\RequirePackage{fix-cm}
\documentclass[3p,times]{article}
\usepackage{placeins}
\usepackage{amssymb}
\usepackage{graphicx,color,ulem}
\usepackage{float}	
\usepackage{amsmath,amssymb}
\usepackage{subfig}
\usepackage{color}
\usepackage{multirow}
\usepackage{algorithmic}
\usepackage{algorithm}
\usepackage{url}
\newcommand{\x}{\mathbf{x}}
\newcommand{\xa}{\mathbf{x}^\textnormal{a}}
\newcommand{\xb}{\mathbf{x}^\textnormal{b}}

\newcommand{\xtin}{\mathbf{x}_0}

\newcommand{\xbara}{\overline{\textbf{x}}^{\rm a}}
\newcommand{\y}{\mathbf{y}}

\newcommand{\nobs}{\textsc{n}_\textnormal{obs}}
\newcommand{\nens}{\textsc{n}_\textnormal{ens}}
\newcommand{\nvar}{\textsc{n}_\textnormal{var}}
\newcommand{\p}{\mathbf{p}}

\newcommand{\PD}{\mathcal{P}}
\newcommand{\Pa}{\mathcal{P}^{\rm a}}
\newcommand{\Pb}{\mathcal{P}^{\rm b}}
\newcommand{\Bini}{\mathbf{B}_0}
\graphicspath{ {Plots/} }
\title{A Hybrid Monte-Carlo Sampling Smoother for Four Dimensional Data Assimilation}
\author{
	Ahmed Attia \and Vishwas Rao \and Adrian Sandu 
}
\date{\today}
\begin{document}
 \thispagestyle{empty}
\setcounter{page}{0}

\begin{Huge}
\begin{center}
%Computer Science Technical Report CSTR-{\tt insert number here} \\
Computational Science Laboratory Technical Report CSL-TR-19-2015\\
\today
\end{center}
\end{Huge}
\vfil
\begin{huge}
\begin{center}
Ahmed Attia, Vishwas Rao, and Adrian Sandu
\end{center}
\end{huge}

\vfil
\begin{huge}
\begin{it}
\begin{center}
``A Hybrid Monte-Carlo Sampling Smoother for Four Dimensional Data Assimilation''
\end{center}
\end{it}
\end{huge}
\vfil

\begin{large}
\begin{center}
Computational Science Laboratory \\
Computer Science Department \\
Virginia Polytechnic Institute and State University \\
Blacksburg, VA 24060 \\
Phone: (540)-231-2193 \\
Fax: (540)-231-6075 \\  
Email: \url{sandu@cs.vt.edu} \\
Web: \url{http://csl.cs.vt.edu}
\end{center}
\end{large}

\vspace*{1cm}

\begin{tabular}{ccc}
\includegraphics[width=2.5in]{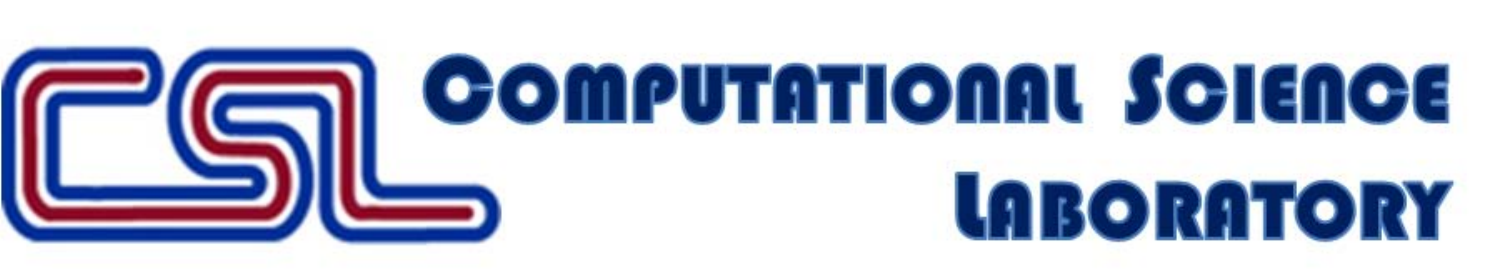}
&\hspace{2.5in}&
\includegraphics[width=2.5in]{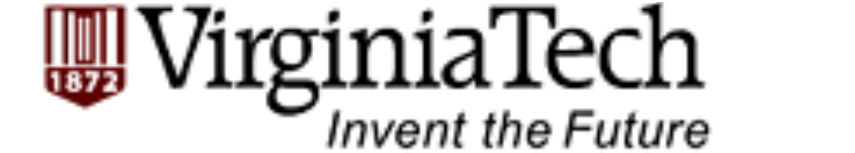} \\
{\bf\em Innovative Computational Solutions} &&\\
\end{tabular}

\newpage

 \maketitle
%\maketitle
%
%
%
%%------------------------------------------
\begin{abstract}
This paper constructs an ensemble-based sampling smoother for four-dimensional data assimilation using a Hybrid/Hamiltonian Monte-Carlo  approach. 
The smoother samples efficiently from the posterior probability density of the solution at the initial time.
Unlike the well-known ensemble Kalman smoother, which is optimal only in the linear Gaussian case, the proposed methodology naturally accommodates non-Gaussian errors and non-linear model dynamics and observation operators. 
Unlike the four-dimensional variational met\-hod, which only finds a mode of the posterior distribution, the smoother provides an estimate of the posterior uncertainty. One can use the ensemble mean as the minimum variance estimate of the state, or can use the ensemble in conjunction with the variational approach to estimate the background errors for subsequent assimilation windows. 
Numerical results demonstrate the advantages of the proposed method compared to the traditional variational and ensemble-based smoothing methods.
%
%----------------------------------<<<<<<<<<<<<< 
%\keywords{Data Assimilation \and Variational Methods \and Ensemble Smoothers \and Markov Chain \and Hybrid Monte Carlo}
%----------------------------------<<<<<<<<<<<<< 
%
\end{abstract}
%%------------------------------------------
%
%
\newpage
\tableofcontents
\newpage\setcounter{page}{1}
%
%%---------------------------------------------------------------- 
\section{Introduction}
\label{sec:Introduction}
Data assimilation (DA) is the process of combining information from predictions made by imperfect models, from noisy observations, and from priors to produce a consistent description of the state of a dynamical system. The application of DA to large scale systems such as the atmosphere is of great practical interest. Two approaches for solving large DA problems have gained widespread acceptance. The first approach, originating from control theory, are variational methods such as the three-dimensional variational (3D-Var) and  four-dimensional variational (4D-Var) strategies \cite{Sandu_2011_subspaceDA}. The variational methods find a maximum a-posteriori (MAP) estimate of the true state of the system. The second approach are the ensemble-based statistical estimation methods. The most successful family of ensemble data assimilation algorithms includes the ensemble Kalman filter (EnKF)~\cite{Evensen_1994} and its variants, the square-root Kalman filters~\cite{Whitaker_2002a},  the ensemble adjustment Kalman filter~\cite{Anderson_2001}, the ensemble transform Kalman filter~\cite{Bishop_2001a}, and efficient implementations of Kalman Filter, such as  ~\cite{Ruiz:2015}, using the Sherman-Morrison formula. All variants of EnKF provide a minimum variance estimate of the state by approximating the expected value of the posterior distribution. Variational and statistical estimation methods yield identical estimates (only) in the case of linear dynamics, linear observations, and Gaussian errors.
  
Both 3D-Var and EnKF are filtering methods that estimate the true state of the system at the specific time instances where observations are available.  For many oceanographic, meteorological, and hydrological applications, it is advantageous to employ smoothing methods such as 4D-Var and EnKS that simultaneously use information from all observations available at different time  points within an assimilation time window. Strong-constraint 4D-Var  updates the state of the system at the initial time of an assimilation window given a background estimate of the initial condition and a set of observations distributed through this interval. The ensemble Kalman smoother, on the other hand, estimates the posterior distributions of the state at time points in the window given all past, present, and future observations (in the assimilation window). 
   
4D-Var requires the derivation and implementation of the tangent linear and the adjoint numerical models, a challenging and effort-intensive task for large-scale models. The 4D-Var algorithm provides a single analysis state, the best posterior estimate of the state of the system. The uncertainty in the estimated state is not inherently provided by the 4D-Var algorithm~\cite{cheng2010hybrid}.
Previous work proposed to quantify the uncertainty in the 4D-Var analysis, by approximating the analysis error covariance 
matrix using an ensemble of simulations~ \cite{gustafsson2014four,hunt04}. These schemes provide an approximation of the analysis error covariance matrix that is inconsistent with the 4D-Var analysis itself because the covariance estimates are usually obtained from independent schemes such as EnKF. Approaches to quantify 4D-Var analysis uncertainty based on subspace error  decompositions \cite{cheng2010hybrid,Sandu_2011_subspaceDA} are statistically consistent but require additional computational effort.

The EnKS is optimal when the observation operators are linear and the errors are Gaussian. However, these assumptions are unlikely to hold for real applications. The analysis ensemble generated by the EnKS allows to find a minimum variance estimate (e.g., ensemble mean), as well as a measure of the analysis uncertainty  (e.g., the analysis error covariance matrix). When the posterior distribution is nearly Gaussian EnKS offers an efficient practical algorithm. However, when the observation operators are nonlinear and the errors are non-Gaussian, the EnKS is not expected to yield good results. 
 
The Markov Chain Monte Carlo (MCMC) family of algorithms provides a powerful foundation to sample from complicated probability distributions.
These algorithms work in general by generating a Markov chain whose stationary distribution is the target probability distribution. MCMC sampling is considered to be the gold standard~\cite{law2012evaluating} in data assimilation. The main practical limitation of MCMC  is the considerable computational cost required to achieve convergence, and to explore the entire state space  in the case of high dimensional state spaces. Scalable and accelerated MCMC algorithms are being continuously developed to improve convergence and space exploration. 
Hybrid Monte Carlo (HMC), also known as Hamiltonian Monte Carlo, is an accelerated MCMC sampling algorithm that reduces the correlation between successive states by using Hamiltonian dynamics to generate proposal states~\cite{duane1987hybrid}. Moreover, HMC targets states with high acceptance probability leading to fast convergence and fast space exploration. 
To the best of our knowledge, HMC was first considered in the context of DA in~\cite{bennett1994open} to solve a nonlinear inverse problem by minimizing 
the residual between the solution and ill-posed boundary conditions. Posterior error statistics are approximated by sampling the nearby states to the optimal 
state after convergence. In~\cite{bennett1994open} a simulated annealing strategy is used during the sampling process where the minimum is obtained at a low
temperature and posterior samples are collected at high temperatures. Solving the weak-constraint 4D-Var problem using a gradient method then using HMC to estimate the analysis error statistics is also discussed  in ~\cite[Chapter 6]{Evensen_2007_book}.
%
%--------------------------------------------------------------------------------------------------------------------------

This work develops a nonlinear non-Gaussian smoo\-ther to solve the four dimensional data assimilation problem.  The new method uses a Hybrid Monte Carlo approach to sample the posterior distribution and is named the HMC smoother. This work extends the sampling filter, proposed in~\cite{Attia_HMCFilter_TR}, to the four dimensional case where time-distributed observations are assimilated at once. HMC smoother provides an ensemble of states sampled from the posterior distribution of the state of the system at the initial time of the  assimilation window. In a practical setting the smoothing step is carried out sequentially over consecutive assimilation windows. The generated ensemble encapsulates the uncertainty in the posterior distribution at the beginning of the current window; when propagated to the beginning of the next assimilation window it provides flow-dependent information about the background error distribution in the next assimilation cycle.

%-------------------------------------added based on HMC_sampling_Filter paper reviews-------------------------------------
%
The use of HMC for solving smoothing problems was presented also in~\cite{alexander2005accelerated}, where a generalized version of HMC is used 
in an attempt to reducing the number of chain steps required to ensure independence of the generated states. The underlying dynamical system in the generalized version of HMC is not Hamiltonian. There are several important differences between this work and ~\cite{alexander2005accelerated}.
In~\cite{alexander2005accelerated} the only source of uncertainty is model error, in form of additive random noise included in the dynamics. Here we work in the strong constraint 4D-Var framework and consider the model to be perfect; the system state is uncertain  due to uncertainties in the initial conditions. While in~\cite{alexander2005accelerated} numerical experiments are carried out with a one-dimensional system, here we experiment with shallow water equations over the sphere, a moderately large multidimensional  nonlinear model relevant for geophysical applications.  Finally, we propose to use higher order symplectic integrators, as tested in~\cite{Attia_HMCFilter_TR}, to efficiently sample from complex posterior distributions that arise when the observation operators and models are highly nonlinear.
%
%--------------------------------------------------------------------------------------------------------------------------
    
The remaining part of the paper is organized as follows. 
Section \ref{sec:Data_Assimilation} reviews the variational and the ensemble approaches for
solving the data assimilation problem. 
The HMC smoother is presented in Section \ref{sec:HMC_and_Sampling_Smoother}. 
Experimental settings and numerical results are discusses in Section \ref{sec:Numerical_Experiments_and_Results}.
Conclusions and future directions are summarized in Section \ref{sec:Conclusion_and_Future_Work}.
  
%----------------------------------<<<<<<<<<<<<<
 
%
%
%%---------------------------------------------------------------- 
\section{Data Assimilation}
\label{sec:Data_Assimilation}

This section reviews the 4D-Var and the EnKS data assimilation schemes.
\subsection{Four-dimensional variational data assimilation}
\label{subsec:4D-VAR_DA}
4D-Var calculates the optimal initial condition for the state of the dynamical system over a specific assimilation time window, based on a background state and using all observations available within this time window \cite{sandu2011chemical}. The background initial state is usually the forecast produced by propagating the previous window analysis through the model dynamics.
To be specific, let the current assimilation window be the time interval $[t_0\,,t_F]$. Given a background state $\xb_0 =\xb[t_0]$, and a set of observations
$\{\y_k=\y[t_k]\}_{k=0,1,\ldots,\nobs}$, available at the discrete time points $\{t_k\}_{k=0,1,\ldots,\nobs} \subset [t_0\,,t_F]$, the 4D-Var analysis is obtained by solving the following optimization problem:
\begin{align}
\label{eqn:4DVAR_Optimization}
 \min_{\xtin}{
   \mathcal{J}(\xtin) } &= % \frac{1}{2} (\xtin-\xb_0)^T \mathbf{B}_0^{-1} (\xtin-\xb_0)  \\
                          \frac{1}{2}  \left\Vert \xtin-\xb_0 \right\Vert ^2 _{\mathbf{B}_0^{-1} }  \\
           % &  \quad   + \frac{1}{2}\sum_{k=0}^{\nobs}{ \Bigl( (\mathcal{H}_k(\x_k)-\y_k)^T \mathbf{R}_k^{-1} (\mathcal{H}_k(\x_k)-\y_k) \Bigr)}
                        &+ \frac{1}{2}\sum_{k=0}^{\nobs}{ \left\Vert \mathcal{H}_k(\x_k)-\y_k \right\Vert^2 _ {\mathbf{R}_k^{-1}}  }\,.  \nonumber  
\end{align}
Here $\mathcal{H}_k$ is the observation operator (generally nonlinear) that maps the model space into the observation space at time point $t_k$. The dimension of observation space $\rm m$ is usually much lower than the dimension of the state space, that is $\rm m \ll \nvar $. $\mathbf{B}_0$ is the background error covariance matrix, and $\mathbf{R}_k$'s are the observation error covariance matrices at each times $t_k;\,k=1,\ldots,\nobs$. 
The background error covariance matrix determines how information from observed areas are extrapolated to unobserved regions or where observations  are sparsely available~\cite{Tremolet_2006}.
The state $\x_k=\x[t_k]$, is produced by propagating the initial state $\xtin$ through the model dynamics from time $t_0$ to point $t_k$
\begin{equation}
\label{eqn:Forward_Model}
 \x_k = \mathcal{M}_{0,k}(\xtin).
\end{equation}
The model solution operator $\mathcal{M}$ represents the discretized partial differential equations that govern the evolution of the dynamical system. Realistic atmospheric and ocean models typically have $\nvar \sim 10^6 - 10^9$ state variables. 

In this work we consider the strong-constraint 4D-Var case which assumes that the numerical model \eqref{eqn:Forward_Model} is perfect. The methodology proposed here can be immediately extended to the weak-constraint 4D-Var framework~\cite{Sasaki_1970}, where model errors are accounted for by adding the corresponding model error terms to equation \eqref{eqn:4DVAR_Optimization} \cite{Tremolet_2006}. 

Perturbations (small errors $\delta \x$ ) of the state of the system evolve according to the tangent linear model:
\begin{equation}
\label{eqn:TLM}
\delta \x_k = \mathbf{M}_{0,k} (\x_0) \cdot \delta \x_0\ ,\ \ t_0 \leq t \leq t_F\,, 
\end{equation}
where $\mathbf{M}_{0,k} = \partial \mathcal{M}_{0,k}/\partial \x[t_0]$ is the Jacobian of the model solution operator.

In strong-constraint 4D-Var the model dynamics act as constraints to the  optimization problem \eqref{eqn:Forward_Model}.  The optimal initial condition obtained by solving the optimization problem \eqref{eqn:4DVAR_Optimization} constrained by the model dynamics \eqref{eqn:Forward_Model}, is referred to as the analysis state $\xa_0=\xa[t_0]$.
The gradient of the cost functional \eqref{eqn:4DVAR_Optimization} is:
%
%\begin{equation}
\begin{align}
\label{eqn:gradient}
\nabla _{\xtin} \mathcal{J}(\xtin) &= {\mathbf{B}^{-1} ( \xtin - \xb_0 )}    \\
& + \sum_{k=0}^{\nobs}{ \mathbf{M}^T_{0,k} \mathbf{H}_k^T
    \mathbf{R}_k^{-1} ( \mathcal{H}_k(\x_k) - \y_k ) } \,,  \nonumber
\end{align}
%\end{equation}
%
where $\mathbf{M}^T_{0,k} $ is the adjoint of the tangent linear model operator \eqref{eqn:TLM}, $\mathbf{H}_{k} = \partial \mathcal{H}_k/\partial \x_k$ is the Jacobian of the observation operator $\mathcal{H}_k$, and $\mathbf{H}_{k}^T$ is its adjoint. 
Gradient-based minimization using \eqref{eqn:gradient} requires the development of the tangent linear and the adjoint models, which is a challenging task for high-dimensional complex models of practical interest. 
The performance of the optimization can be improved by using the second order derivative information and adaptive observations as described in \cite{Cioaca2014377}.
%

%%%%%%%%%%%%%%%%%%%%%%%%%%%%%%%%%%
\subsection{Bayesian interpretation of 4D-Var} 
\label{subsec:Bayesian_4D-VAR_DA}
%%%%%%%%%%%%%%%%%%%%%%%%%%%%%%%%%%
%
The knowledge of the system state at the initial time $t_0$ prior to obtaining new observations is described by the background (prior) probability density $\mathcal{P}^b(\xtin)$. The ``sampling model'' gives the probability distribution of observations conditioned by the initial state \\
$\mathcal{P}(\y_0,\y_1,\ldots,\y_{\nobs}|\xtin)$, 
under the belief that the dynamical model $\x_k = \mathcal{M}_{0,k}(\x_0)$ perfectly represents reality.  From Bayes' theorem:
\begin{subequations}
\label{eqn:Bayes}
\begin{align}
  \Pa(\xtin) &= 
  \PD(\xtin| \y_0,\y_1,\ldots,\y_{\nobs}) \\
&= \frac{ \PD(\y_0,\y_1,\ldots,\y_{\nobs} | \xtin)\, \Pb(\xtin) }{\PD(\y_0,\y_1,\ldots,\y_{\nobs}) }, \label{eqn:Bayes1_basicForm}
%\\
%&\propto \PD(\y_0,\y_1,\ldots,\y_{\nobs}|\xtin) \Pb(\xtin)  \label{eqn:Bayes1_proportionalityForm} \\
% \label{eqn:Bayes1_likelihoodForm}
%&\propto \mathcal{L}(\xtin | \y_0,\y_1,\ldots,\y_{\nobs}) \Pb(\xtin)\,,
\end{align}
\end{subequations}
%
%where \eqref{eqn:Bayes1_likelihoodForm} uses the likelihood function $\mathcal{L}$. 
The posterior (analysis) PDF $\Pa(\xtin)$ is the probability distribution of the initial state after incorporating the new knowledge contained in the observations. The denominator in \eqref{eqn:Bayes1_basicForm} is the marginal density of the observations and acts as a normalization factor.

The background and observations errors are usually assumed to have Gaussian distributions: 
\begin{subequations}
\begin{align}
\label{eqn:Prior_PDF}
  \Pb(\xtin) & \propto  \exp{ \left( -\frac{1}{2} \left\Vert \xtin-\xb_0\right\Vert^2_{\mathbf{B}_0^{-1}} \right) }\,, \\
\label{eqn:Likelihood_PDF}
\PD(\y_k | \x_k) %&= (2\pi)^{-\frac{\rm m}{2}} |\mathbf{R}_k|^{-\frac{1}{2}} \\
    % & \times \exp{ \left( -\frac{1}{2} (\mathcal{H}_k(\x_k)-\y_k)^T \mathbf{R}_k^{-1} (\mathcal{H}_k(\x_k)-\y_k) \right) }  \nonumber\\
    & \propto  \exp{ \left( -\frac{1}{2} \left\Vert \mathcal{H}_k(\x_k)-\y_k\right\Vert^2_{ \mathbf{R}_k^{-1}} \right) },
\end{align}
where $\mathbf{B}_0$ is the background error covariance matrix and $\mathbf{R}_k$'s are the observation error covariance matrices at times $t_k;\,k=1,\ldots,\nobs$.
If the observation errors at different time points are independent, and the model is perfect, the joint sampling model can be written as
\begin{multline}    
\PD(\y_0,\y_1,\ldots,\y_{\nobs} | \xtin) \propto  \\
   \exp{\left( \sum_{k=0}^{\nobs}{ \left( -\frac{1}{2} \left\Vert \mathcal{H}_k(\x_k)-\y_k\right\Vert^2_{ \mathbf{R}_k^{-1}} \right)}\right)}\,.%(\mathcal{H}_k(\x_k)-\y_k)^T \mathbf{R}_k^{-1} (\mathcal{H}_k(\x_k)-\y_k) \right) } \right)}\,. 
\end{multline}
\end{subequations}
Bayes' rule \eqref{eqn:Bayes} yields the following posterior PDF
\begin{subequations}
\label{eqn:posterior_4DVAR}
\begin{align}
        \Pa(\xtin) & \propto    \exp{\Bigl( - \mathcal{J}(\xtin) \Bigr)}, \\ 
\mathcal{J}(\xtin) & =  \frac{1}{2} \left\Vert \xtin-\xb_0 \right\Vert ^2 _{\mathbf{B}_0^{-1}} \label{eqn:4DVAR_Cost_Functional_HMC} 
                     +  \frac{1}{2} \sum_{k=0}^{\nobs}{ \left\Vert  \mathcal{H}_k(\x_k)-\y_k \right\Vert ^2 _{\mathbf{R}_k^{-1}} }\,. 
\end{align}
\end{subequations}
For nonlinear models and nonlinear observation operators the posterior \eqref{eqn:posterior_4DVAR} is not Gaussian. The kernel of the posterior is $\exp(-\mathcal{J}(\xtin))$, where $\mathcal{J}(\xtin)$ is the cost functional of the 4D-Var problem. 4D-Var computes the analysis $\xa_0$ as the minimizer of $\mathcal{J}$. The 4D-Var solution is the MAP  estimate of the initial state under the assumptions that the background and observation errors are Gaussian.  For highly nonlinear observation operators and highly nonlinear models the posterior can have multiple modes. In this case the 4D-Var numerical solution can be trapped in a local minimum of the cost functional. 

The posterior distribution \eqref{eqn:Bayes} contains the complete characterization of the uncertain initial state of the dynamical system. However, calculating the full posterior with high dimensional models is infeasible in practice.  A practical approach is to describe the posterior probability density by an ensemble of states, and to use it to estimate moments of the distribution. Sampling directly from the posterior PDF of the initial condition \eqref{eqn:posterior_4DVAR} acts as a smoother; the ensemble mean provides an estimate of the true state of the system, and the ensemble covariance an estimate of the posterior uncertainty that is totally consistent with the analysis state. In contrast, the current practice to use the analysis obtained from 4D-Var and the covariance obtained from EnKF or EnKS leads to inconsistent representations of uncertainty.
%
%%
%% 

%%%%%%%%%%%%%%%%%%%%%%%%%%%%%%%%%%%%%%%%%%%%%%%
\subsection{Ensemble Kalman filter and smoother}
\label{subsec:EnKS}   
%%%%%%%%%%%%%%%%%%%%%%%%%%%%%%%%%%%%%%%%%%%%%%%
%
Filtering  is the process of calculating the posterior distribution of the uncertain state of a dynamical system at a specific time given observations only available at that time instance. The ensemble Kalman filter~\cite{Evensen_1994} represents probability distributions by samples.  Let $\{\xb_k(e)\}_{e=1,2,\ldots,{\nens}}$ be an ensemble of forecast states at time $t_k$,  and $\y_k=\y[t_k]$ the observation vector (set of measurements) at $t_k$. If the forecast (background) ensemble is represented by the matrix $\mathbf{X}_k^{\rm b}$, whose columns are the forecast ensemble members, then the updated (analysis) ensemble matrix $\mathbf{X}_k^{\rm a}$ at the same time $t_k$ is obtained as~\cite{Evensen_2003}
\begin{equation}
\label{eqn:EnKF_Update}
\mathbf{X}_k^{\rm a} = \mathbf{X}_k^{\rm b} \cdot \mathbf{T}_k\,, 
\end{equation}
where the matrix $\mathbf{T}_k \in \mathbb{R}^{\nens \times \nens}$ is a nonlinear transformation constructed from the forecast ensemble and the observations at time $t_k$~\cite{Evensen_2003}.  Square-root filters~\cite{Whitaker_2002a}, the ensemble transform Kalman filter~\cite{Bishop_2001a}, and the ensemble adjustment Kalman filter~\cite{Anderson_2001} can all be written in the form \eqref{eqn:EnKF_Update} for specific choices of the transformation $\mathbf{T}_k$.

Smoothing is the process of calculating the posterior distribution of the uncertain states of a dynamical system given past, present, and future observations~\cite{briers2010smoothing}. There are three approaches to smoothing: fixed-point,  fixed-interval, and fixed-lag smoothing~\cite{briers2010smoothing}, with the second and third ones being the most popular~\cite{ravela2007fast}. Any scheme that can be used to solve any of the three smoothing problems can also be employed as a single  fixed-interval smoothing scheme~\cite{briers2010smoothing}.
The ensemble smoother (ES) was introduced in~\cite{Evensen_1996} as a linear variance minimization algorithm. The ensemble Kalman smoother (EnKS)~\cite{Evensen_2000} employs an ensemble of states to describe distributions and obtains the posterior using the Kalman filter updates equations.  EnKS is optimal in case of linear dynamics, Gaussian errors, and large number of ensemble members ~\cite{Evensen_2007_book}. 

To construct EnKS~\cite{Evensen_2000} the EnKF update equations \eqref{eqn:EnKF_Update} are used repeatedly to develop a fixed-lag, and a fixed-interval smoothing algorithms. A fixed-point smo\-other can be written as~\cite{Evensen_2003}
\begin{equation}
\label{eqn:EnKS_FixedPoint_Update}
\mathbf{X}^{\rm s}_0 = \mathbf{X}^{\rm a}_0\cdot \prod_{k=0}^{\nobs} {\mathbf{T}_k}\,.
\end{equation}
The update equation \eqref{eqn:EnKS_FixedPoint_Update} is used recursively for fixed-interval smoothing, where smoothed ensembles are obtained at specified set of times, and they are conditioned only on observations available at later times in the interval. Ravela and McLaughlin~\cite{ravela2007fast} presented efficient, fast versions of the fixed-interval and the fixed-lag EnKS. 
The fast fixed-interval smoother has a computational cost that scales linearly with respect to the length of the interval. In this work, we use the fast fixed-interval EnKS~\cite{ravela2007fast}, with a single smoothing point (fixed-point smoother) chosen at the beginning of the time interval.

EnKS computes the minimum variance estimate of the state. This is not expected to be very accurate if the observations are highly nonlinear or if the Gaussianity assumptions are severely violated. As shown in Section \ref{sec:Numerical_Experiments_and_Results}, the HMC sampling smoother proposed herein is capable of handling nonlinear observation and model operators, and consequently produces posterior estimates that are more useful than the EnKS ensemble, and contain more information than the 4D-Var MAP analysis. 
Hybrid methods such as \cite{Ruiz:2015a} make use of optimization over ensembles using the trust-region framework.

%%%%%%%%%%%%%%%%%%%%%%%%%%%%%%%%%%%%%%%%%%%%%%%
\section{The hybrid Monte-Carlo sampling smoother}
\label{sec:HMC_and_Sampling_Smoother}   
%%%%%%%%%%%%%%%%%%%%%%%%%%%%%%%%%%%%%%%%%%%%%%%
    %  
The most popular and successful class of sampling algorithms is the Markov chain Monte-Carlo (MCMC) \cite{neal1993probabilistic}, first introduced by Metropolis \textit{et al.}~\cite{metropolis1953equation}.
MCMC algorithms sample from  a general probability distribution $\PD(\x)$ by building a Markov chain whose invariant distribution is $\PD(\x)$. MCMC algorithms have an advantage of not requiring the normalization of target distributions. However, traditional MCMC samplers are often considered impractical for large dimensional problems due to the following drawbacks: The Markov chain may take a very long time to reach stationarity. A large number of (burn-in) states are generated and discarded before starting the sampling process, in order to guarantee that the collected samples are obtained from the true target PDF. The samples should be independent, however the Markov chain is not completely memoryless; in order to achieve independence of sampled states, the sampler usually drops some intermediate states between each selected state. Another drawback of most of Monte-Carlo sampling methods is the curse of dimensionality~\cite{neal1993probabilistic}: as the dimension of the state spaces grows, the number of sample members needed to represent the probability distribution, grows rapidly. The number of samples required to efficiently represent the probability distribution can be controlled if the sampler surveys sufficiently fast the entire state space. The sampler can become trapped in a high-probability mode of a multi-modal distribution, and fail to represent the other probability modes. 
   
%%%%%%%%%%%%%%%%%%%%%%%%%%%%%%%%%%%%%%%%%%%%%%%
\subsection{Hybrid Monte-Carlo}
\label{subsec:Hybrid_Monte-Carlo}
%%%%%%%%%%%%%%%%%%%%%%%%%%%%%%%%%%%%%%%%%%%%%%%

Hybrid/Hamiltonian Monte-Carlo (HMC)~\cite{duane1987hybrid} follows an auxiliary-variable approach in order to alleviate the limitations of the traditional MCMC algorithms. 
  
The phase space of a Hamiltonian dynamical system consists of points $(\p,\x) \in \mathbb{R}^{2\,\nvar}$,  where $\x \in \mathbb{R}^{\nvar}$ is the position variable, and $\p \in \mathbb{R}^{\nvar}$ is the momentum variable. The Hamiltonian dynamics is governed by the set of ordinary differential equations (ODEs):
\begin{subequations}
\label{eqn:Hamiltonian}
\begin{equation}
 \label{eqn:Hamiltonian_equations}
 \begin{split}
 \frac{d\x}{dt} &= \nabla_\p\, H(\p,\x)\,, \\
 \frac{d\p}{dt} &= - \nabla_\x\, H(\p,\x)\,,
 \end{split}
\end{equation} 
where the Hamiltonian function $H$ describes the total energy of the system 
\begin{equation}
\label{eqn:Hamiltonian_function}
H(\p,\x)  =  \frac{1}{2} \, \p^T \mathbf{M}^{-1} \p + \mathcal{J}(\x)\,.
\end{equation}  
\end{subequations}
The first term of the Hamiltonian \eqref{eqn:Hamiltonian_function} represents the potential energy of the system, while the second term corresponds to the kinetic energy. The exact (analytic) flow of the Hamiltonian system \eqref{eqn:Hamiltonian_equations} advances the solution in time from $t=0$ to $t=T$:
\begin{equation}
\label{eqn:Hamiltonian_flow}
\Phi_T:\mathbb{R}^{2\nvar} \rightarrow \mathbb{R}^{2\nvar}; \quad \Phi_T\bigl(\p[0],\x[0]\bigr)=\bigl(\p[T],\x[T]\bigr)\,.
\end{equation}
This flow cannot be calculated exactly in practice, and has to be approximated by an equivalent numerical solution using a time reversible symplectic integrator~\cite{sanz2014Markov,sanz1994numerical}. The most common symplectic integrator is leapfrog (St{\"o}rmer--Verlet)  \cite{sanz2014Markov,sanz1994numerical}. One step of the position Verlet algorithm advances the solution of the Hamiltonian equations \eqref{eqn:Hamiltonian_equations}  from time $t_j$ to time $t_{j+1} = t_j + h$ using:
 \begin{subequations}
 \label{eqn:Verlet}
 \begin{eqnarray}
  \x_{j+1/2} &=&  \x_j + \frac{h}{2}\, \mathbf{M}^{-1}\, \p_j\,,   \\
  \p_{j+1}   &=&  \p_j - h \,\nabla_\x \mathcal{J}(\x_{j+1/2})  \,,   \\
  \x_{j+1}   &=&  \x_{j+1/2} + \frac{h}{2}\, \mathbf{M}^{-1} \,\p_{j+1}. 
 \end{eqnarray}
 \end{subequations} 
The optimal time step size $h$ must satisfy $h \propto {(1/{\nvar})}^{1/4}$ ~\cite{beskos2013optimal}, and careful empirical tuning of the step size is usually required for good performance~\cite{Attia_HMCFilter_TR}. Several other symplectic integrators with more stages and higher accuracy than Verlet have also been developed~\cite{blanes2014numerical}.  An infinite dimensional time integrator was also introduced in~\cite{beskos2011hybrid}.

For practical considerations it is advisable to split the interval $[0,T]$ where the Hamiltonian system evolves into $m$ smaller sub steps of length $h=T/m$.  The flow of the numerical solution obtained by the symplectic integrator will be denoted by $\tilde{\Phi}_T$ and is an approximation of the exact flow  $\Phi_T$.
 
The key idea in HMC sampling is to add an auxiliary variable $\p$ to the target variable $\x$ and sample from the joint probability
distribution of $(\x,\p)$. The auxiliary variable is chosen such that the sampling procedure from the joint distribution is much faster
than sampling from the marginal distribution  of the target variable. In HMC sampling the target and the auxiliary variables are thought of as the position and momentum components of a Hamiltonian system, respectively. The Hamiltonian dynamics of the system serves as a transition kernel to the Markov chain.

The kernel of the stationary probability distribution of the Hamiltonian system \eqref{eqn:Hamiltonian} is ~\cite{neal1993probabilistic}
 \begin{subequations}
 \label{eqn:Canonical_Pdf}
 \begin{align} 
 \exp{( - H(\p,\x) )}  &=  \exp{\left( -\frac{1}{2} \p^T \mathbf{M}^{-1} \p - \mathcal{J}(\x) \right)}   \\
&= \exp{\left( -\frac{1}{2} \p^T \mathbf{M}^{-1} \p \right)}\cdot \pi(\x)\,,
 \end{align}
 \end{subequations}
where $\pi(\x)=\exp{\Bigl(-\mathcal{J}(\x)\Bigr)}$ is the probability distribution of the position variable. The joint probability distribution of the state  $(\p,\x)$ in the phase space $\mathbb{R}^{2\,\nvar}$ is the product of the marginal distributions of both the position and the momentum. This simply means that the two variables $\x\, \text{and}\, \p$ are independent~\cite{sanz2014Markov}. Independence of both position and momentum makes it possible to sample from the marginal distribution of each variable by sampling from their joint distribution. The marginal PDF of the momentum variable is a Gaussian distribution with zero mean and covariance matrix $\mathbf{M}$ (also known as the mass matrix), i.e., $\p \sim \mathcal{N}(0,\mathbf{M})$. 

Let $\x \sim \pi(\x)$ be a random variable that is the target of the MCMC sampling algorithm. View $\x$ as the position variable in the Hamiltonian system, and add the momentum $\p$ as an auxiliary variable. The symplectic integrator is used to propose a state that is either accepted or rejected using an acceptance/rejection rule based on the loss of energy. Algorithm \ref{alg:HMC_Sampling} \cite{sanz2014Markov} summarizes the HMC steps to sample from the probability distribution $\pi(\x)$.
 %
 %----------------------------------------------HMC-Sampling Algorithm-----------------------------------------------
 \begin{algorithm}[htpb]
 \begin{algorithmic}[1]
\STATE Initialize the Markov chain. Preferably $(\p_0,\x_0)$ should have high probability w.r.t. the target distribution.
\STATE At each step $k$ of the Markov chain draw the random auxiliary variable $\p_k \sim \mathcal{N}(0,\mathbf{M})$.
\STATE Use a symplectic numerical integrator (e.g. position Verlet) to advance the current state $(\p_k,\x_k )$
by a time increment $T$ to obtain a \textit{proposal} state $( \p^* , \x^* )$: 
\begin{equation}
    ( \p^* , \x^* ) = \tilde{\Phi}_T\bigl(( \p_k , \x_k )\bigr).
\end{equation}
\STATE Use the Hamiltonian \eqref{eqn:Hamiltonian_function} to approximate the loss of energy $\Delta H$.
\STATE Calculate the acceptance probability:
\begin{equation}
\label{eqn:acceptance_probability}
a^{(k)} = 1 \wedge e^{-\Delta H}.
\end{equation}
\STATE Discard both $\p^*$ and $\p_k$.
\STATE \textit{(Acceptance/Rejection)} Draw a uniform random variable $u^{(k)}\sim \mathcal{U}(0,1)$:
\begin{enumerate}
\item[i-  ] If $a^{(k)} > u^{(k)}$ accept the proposal as the next sample: $\x_{k+1} := \x^*$;
\item[ii- ] If $a^{(k)} \leq u^{(k)}$ reject the proposal and continue with the current state: $\x_{k+1} := \x_k$.
\end{enumerate}
\STATE Repeat steps $2$ to $7$ until sufficiently many distinct samples are drawn.
 \end{algorithmic}
 \caption{HMCMC Sampling~\cite{sanz2014Markov}.}
 \label{alg:HMC_Sampling}
 \end{algorithm}
 %-------------------------------------------------------------------------------------------------------------------
 %
The loss of energy between the current and the proposed state is usually calculated as the difference between the Hamiltonians at the current and the proposed states:
 \begin{equation}
 \label{eqn:loss_of_Energy_Verlet}
  \Delta H =  H( \p^* , \x^* ) - H( \p_k, \x_k )\,.
 \end{equation} 
This equation \eqref{eqn:loss_of_Energy_Verlet} is valid for the Verlet \eqref{eqn:Verlet}, two-stage, three-stage, and four-stage symplectic integrators~\cite{blanes2014numerical,sanz2014Markov}. See~\cite{Attia_HMCFilter_TR} for details on different symplectic integrators and corresponding expressions for energy.
The length of the Hamiltonian trajectory $T$ and the number of steps $m$ are parameters to be tuned by the user~\cite{neal2011mcmc}. Another user-tunable parameter is the mass matrix $\mathbf{M}$, a symmetric positive definite matrix that represents the covariance of the momentum variable. The choice of the mass matrix does not alter the fact that HMC sampling Algorithm \ref{alg:HMC_Sampling} converges to the stationary distribution $\pi(\x)$. However, a good choice of $\mathbf{M}$ can considerably improve sampling efficiency. One popular and simple choice is to take $\mathbf{M}$ a constant multiple of the identity matrix. Ideally, if the variances of the target distribution $\pi(\x)$ are known (or can be approximated), the diagonal of $\mathbf{M}$ should be chosen as the corresponding precisions (reciprocals of these variances)~\cite{neal2011mcmc}. We found that this choice results in a very fast convergence of the chain to stationarity. 

HMC sampling Algorithm \ref{alg:HMC_Sampling} tends to explore the state space faster than traditional MCMC, and  the acceptance probability of all generated states is close to one. Several enhancements to the HMC sampling, such as parallel tempering~\cite{earl2005parallel,swendsen1986replica},
have been proposed to guarantee that the algorithm escapes local modes of high probability.
    
%%%%%%%%%%%%%%%%%%%%%%%%%%%%%%%%%%%%%%%%%%%%%%%
\subsection{Sampling smoother algorithm}
\label{subsec:Sampling_Smoother}
%%%%%%%%%%%%%%%%%%%%%%%%%%%%%%%%%%%%%%%%%%%%%%%
%
We now present the HMC smoother (smoothing by sampling) that simultaneously accounts for all observations available within a specific assimilation window to obtain posterior estimates of the initial system state.

Consider the assimilation window $[t_0,t_F]$ with a set of observations available at times $t_0,t_1,\ldots,t_{\nobs}$ inside the window, where $t_{\nobs} \equiv t_F$. Under the assumptions discussed in Section \ref{subsec:Bayesian_4D-VAR_DA} the posterior (analysis) probability distribution of the initial state $\x_0$ takes the form \eqref{eqn:posterior_4DVAR}.
We seek to sample from this posterior distribution using the HMC approach. For this we set the potential energy term in \eqref{eqn:Hamiltonian_function} to be the 4D-Var cost functional \eqref{eqn:4DVAR_Cost_Functional_HMC}. Consequently the  target probability distribution $\pi(\x)$ coincides with the 4D-Var posterior distribution \eqref{eqn:posterior_4DVAR}, i.e., $\pi(\x) = \exp(-\mathcal{J}(\x))$.

The smoother works sequentially \textit{over consecutive assimilation windows} by applying the forecast and analysis (sampling) steps in succession. In the forecast step each state of the analysis ensemble is propagated in time to the end of the previous assimilation window (the beginning of the current window). The result of the forecast step is a forecast ensemble $\mathbf{X}^{\rm b}$ (or just a single background state $\xb_0$) at the beginning of the current time window, i.e. at $t_0$. One can just propagate the analysis state (e.g. the mean of the analysis ensemble) to obtain the current background state  $\xb_0$. However, propagating the full analysis ensemble makes it possible to build an ensemble-based (flow-dependent) background error covariance matrix at the beginning of the current window. This background error covariance matrix includes the errors of the day, and can considerably enhance the quality of the analyses generated by a data assimilation scheme.  We will assume herein that the full forecast ensemble is generated in the forecast step. In the analysis step, the HMC sampling strategy summarized in Algorithm \ref{alg:HMC_Sampling} is applied to obtain the analysis ensemble at the initial time of the current assimilation window. 
 
The HMC sampling smoother is detailed in Algorithm \ref{alg:Sampling_Smoother_Algorithm}.
\begin{algorithm}[H]
\begin{algorithmic}[1]
\STATE \textbf{Analysis step:} 
Given the background state and observations, draw an ensemble of initial  states from the posterior distribution \eqref{eqn:posterior_4DVAR} as follows:
\begin{enumerate}

 \item[i- ]  Calculate an ensemble-based forecast error covariance matrix $\mathbf{B}_0^{\rm ens}$, and use it together with a fixed (modeled) matrix to construct the background error covariance matrix $\mathbf{B}_0$~\cite{Attia_HMCFilter_TR}.
(This step can be omitted by using the modeled background error covariance matrix; however, the use of forecast ensemble is expected to improve the quality of the generated analysis ensemble.)
\item[ii-] Build the mass matrix $\mathbf{M}$ as a diagonal matrix such that $\mathrm{diag\left(\mathbf{M}\right)} = \mathrm{diag\left(\mathbf{B^{-1}_0}\right)}$.
 \item[iii-] Initialize the Markov chain to the best estimate of the current state available, e.g., the background state $\xtin^{\rm{b}}$,
 or a suboptimal 4D-Var solution. A good choice speeds up the convergence of the chain.
 \item[iv- ] Follow the steps in Algorithm \ref{alg:HMC_Sampling} to generate the chain and select ensemble members after the chain reaches stationarity. Dropping a small number (5 to 10) steps between each selected states helps to ensure the independence of the generated ensemble members.
\end{enumerate}
    \STATE \textbf{Forecast step:}
Propagate each member of the analysis ensemble, using the full forward model, to the end of the current assimilation window
(beginning the next assimilation window).
\end{algorithmic}
\caption{The Proposed Sampling Smoother}
 \label{alg:Sampling_Smoother_Algorithm}
 \end{algorithm}

The generated ensemble of states $\{ \xa_0(e)\}_{e=1,2,\ldots ,\nens}$, samples the posterior PDF  $\Pa(\xtin)$, and can be used to calculate the best estimate of the initial condition of the system (e.g., the mean ($\xbara_0$) of the ensemble), and to estimate the analysis error covariance matrix $\mathbf{A}_0$:
 \begin{subequations}
 \label{eqn:Ensemble_mean_and_covariance}   
 \begin{eqnarray}
 \label{eqn:Ensemble_mean}   
\xbara_0&=&  (\nens)^{-1}\, \sum_{e=1}^{\nens}{\xa_0(e) } \,, \\
\nonumber
\Delta\mathbf{X}_0^{\rm a}   &=&  [\xa_0(1)- \xbara_0, \ldots,  \xa_0(\nens)- \xbara_0] \,, \\
 \label{eqn:Ensemble_covariance}   
\mathbf{A}_0   &=&  (\nens-1)^{-1}\, \left( \Delta\mathbf{X}_0^{\rm a} \,  \left( \Delta\mathbf{X}_0^{\rm a} \right)^T \right) \,.
 \end{eqnarray}
 \end{subequations}
The forecast and analysis steps are repeated sequentially on subsequent assimilation windows. The propagated ensemble can be used to estimate the analysis covariance at the final time using \eqref{eqn:Ensemble_covariance}, such as to provide ``flow-dependent'' information for the background error covariance matrix used in the subsequent assimilation interval. 
    
%%%%%%%%%%%%%%%%%%%%%%%%%%%%%%%%%%%%%%%%%%%%%%%
\section{Numerical Experiments}
\label{sec:Numerical_Experiments_and_Results}  
%%%%%%%%%%%%%%%%%%%%%%%%%%%%%%%%%%%%%%%%%%%%%%%

The proposed HMC sampling smoother is tested against the EnKS and the 4D-Var schemes on two numerical models. We first illustrate the distinctive features of the HMC smoother with a simple one-dimensional model with a nonlinear observation operator and a bimodal posterior distribution. Next, we employ the shallow water on the sphere to test the sampling smoother on a problem relevant to geophysics, and to compare its performance against the conventional 4D-Var scheme.

%%%%%%%%%%%%%%%%%%%%%%%%%%%%%%%%%%%%%%%%%%%%%%%
\subsection{A one-dimensional model}
\label{subsec:oneD_Model}
%%%%%%%%%%%%%%%%%%%%%%%%%%%%%%%%%%%%%%%%%%%%%%%

Consider the following model:
\begin{subequations}    
\label{eqn:oneD_Potential_Model}
\begin{align}
  \frac{d\x}{dt} &= -\frac{dV(\x)}{d\x}\,,  \\
  V(\x) &= (\x+1)^2\, (\x-1)^2\,, \label{eqn:oneD_Potential}
\end{align}
\end{subequations}
that describes the position of a particle over the entire real line moving under the effect of the potential field \eqref{eqn:oneD_Potential}. This model is similar to the one used in~\cite{alexander2005accelerated}. The potential field has two local minima at $\pm 1$, which are expected to act as attractors for the particle. We set the reference initial condition to $\x_0^{\rm true}=-0.15$,  and chose the background initial condition to be $\xb_0 = 0.1$.  Note that the true and the background initial conditions lie in the basins of attraction of different equilibria. The background errors are assumed to be normally distributed with zero mean and standard deviation $\sigma_{\xtin}=\sqrt{2}$.

\begin{figure}%[htpb]
\centering
    \includegraphics[width=0.90\linewidth,height=3.8cm]{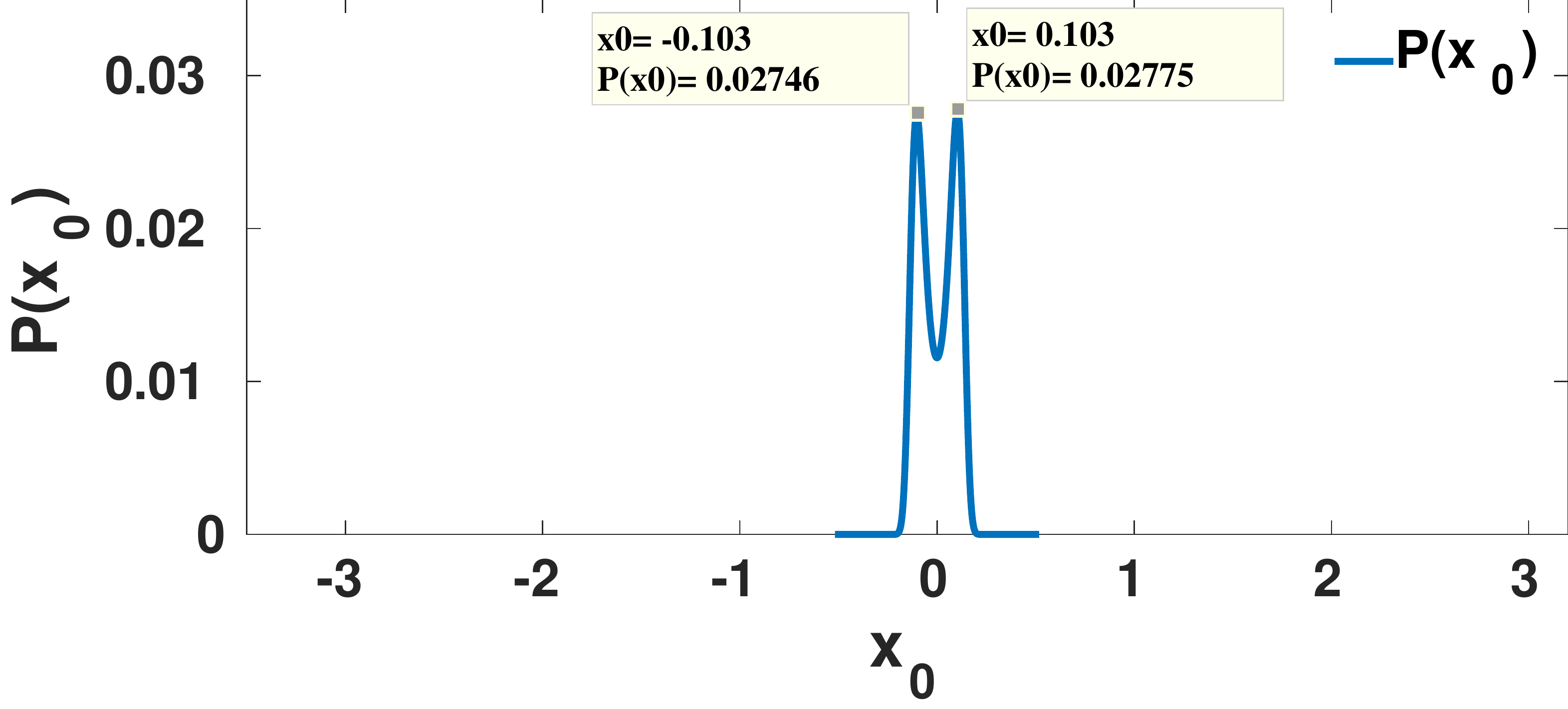}   
\caption{The non-normalized kernel of the posterior distribution \eqref{eqn:OneD_Potential_Posterior_Kernel}. 
         The right peak is slightly higher than than the left one as a result of the prior being a Gaussian distribution centered around $\xb=0.1$ with small standard deviation ($\sigma_{\xtin}=\sqrt{2}$)
         Given the current settings, the right peak occurs at $\xtin=0.103$ with $P(\xtin)=0.02775$ while the left peak occurs at $\xtin =-0.103$ with $P(\xtin)=0.02746$.
         } 
\label{fig:oneD_Potential_Posterior_Kernel}
\end{figure}

Synthetic observations are obtained from the reference solution by applying the quadratic observation operator 
\begin{equation}
\label{eqn:Observation_Operator}
\mathcal{H}(\x_k) = (\x_k)^2\,.
\end{equation}
Observation errors are assumed to be Gaussian with zero mean and standard deviation $\sigma_{\rm obs} = 0.05$. The simulation time window is $[t_0,t_F]=[0,0.12]$ (units), with equally spaced $12$ observation points. The posterior distribution of the initial state reads
%
%\begin{subequations}
%\label{eqn:OneD_Potential_Posterior_Kernel}
\begin{align}
\label{eqn:OneD_Potential_Posterior_Kernel}
 \Pa(\x_0) \propto  &
 %\exp{ \left(-\mathcal{J}(\x)  \right) }\,,  \\
 %\mathcal{J}(\x) &= \frac{1}{2} \left(  \left(\frac{\x_0 - 0.1}{1.41} \right)^2 +
 %   \sum_{k=1}^{\nobs=12}{\left(\frac{\x_k^2 - \y_k}{0.05} \right)^2 }   \right) \,,
 \exp \left(  -\frac{1}{2} \left(\frac{\x_0 - 0.1}{1.41} \right)^2 \right. \\
 \nonumber
 & \qquad \qquad \left. -\frac{1}{2}
    \sum_{k=1}^{\nobs=12}{\left(\frac{\x_k^2 - \y_k}{0.05} \right)^2 }   \right) \,,
\end{align}
%\end{subequations}
%
where $\x_k$ is obtained by propagating $\x_0$ forward in time, from $t_0=0$ to $t_k = k \times 0.01$ (units), using the model \eqref{eqn:oneD_Potential_Model}. The non-normalized  posterior density \eqref{eqn:OneD_Potential_Posterior_Kernel} is illustrated in Figure \ref{fig:oneD_Potential_Posterior_Kernel}.
Traditional assimilation methods, like 4D-Var and EnKS, are expected to have difficulties capturing the bimodal nature of the posterior distribution of the initial condition. 
Since the prior PDF is a Gaussian centered around the background state $\xb=0.1$ with standard deviation $\sigma_{\xtin}=\sqrt{2}$, the right peak in Figure \ref{fig:oneD_Potential_Posterior_Kernel} is slightly taller than the left peak. 
                 With Gaussian background prior centered around one of the peaks, smaller standard deviation would damp the other peak.
Capturing only that right peak completely misses the true solution, which is negative. Numerical results presented below show that the proposed HMC smoother is capable of building a representative ensemble from  the bimodal posterior distribution. The analysis ensemble can then be used to draw more useful conclusions (e.g. statistics) than what can be obtained from  analysis results obtained by the traditional methods.
\subsubsection{Numerical results with the one-dimensional model}    
HMC smoothing was carried out to collect an ensemble of $100$ members from the posterior \eqref{eqn:OneD_Potential_Posterior_Kernel}. We tested several symplectic integrators~\cite{Attia_HMCFilter_TR}, and found that all show similar behavior. We chose the position Verlet symplectic integrator due to its minimal computational cost for all our experiments. The Hamiltonian system step size is empirically tuned to $T=0.1$, with step length $h = 0.01$, and number of steps $m=10$. The number of burn-in steps is chosen to be $20$ (for this simple model we already know that the forecast state $0.1$ lies in the support of the posterior  and the  burn-in steps could be omitted; in general one can incorporate convergence tests to shorten the number of burn-in steps and ensure that the collected samples are from the target distribution). Four states are dropped between consecutive selected states at stationarity to guarantee the independence of the samples. 
    
\begin{figure*}
    \subfloat[HMC sampling smoother.,height=3.8cm]{\label{fig:HMC_Histogram_Poterior_Ensemble}%
 \includegraphics[width=0.45\linewidth,height=3.8cm]{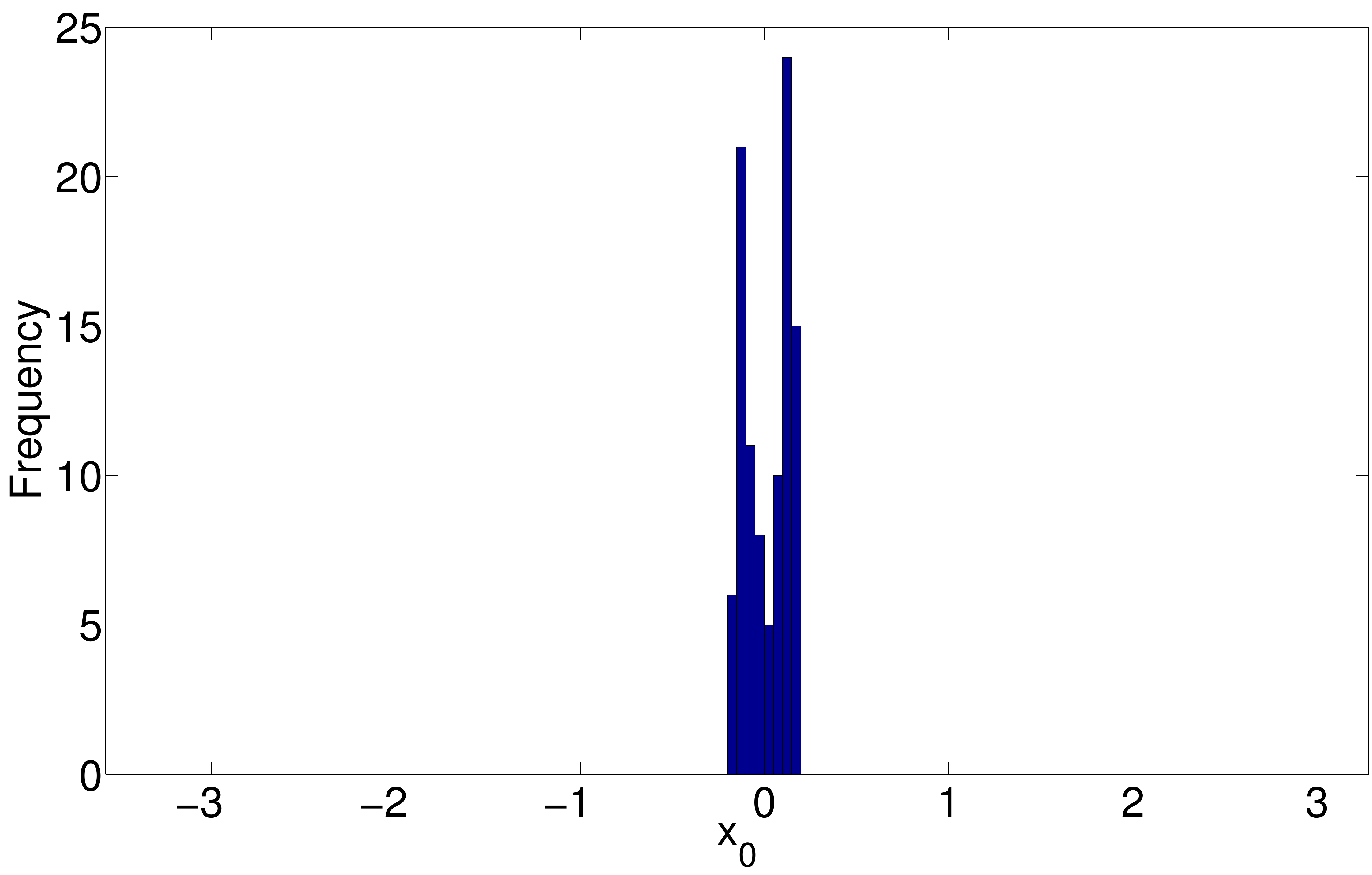}   
}%
    \hfill
    \subfloat[Ensemble Kalman smoother.]{\label{fig:EnKS_Histogram_Poterior_Ensemble}%
 \includegraphics[width=0.45\linewidth,height=3.8cm]{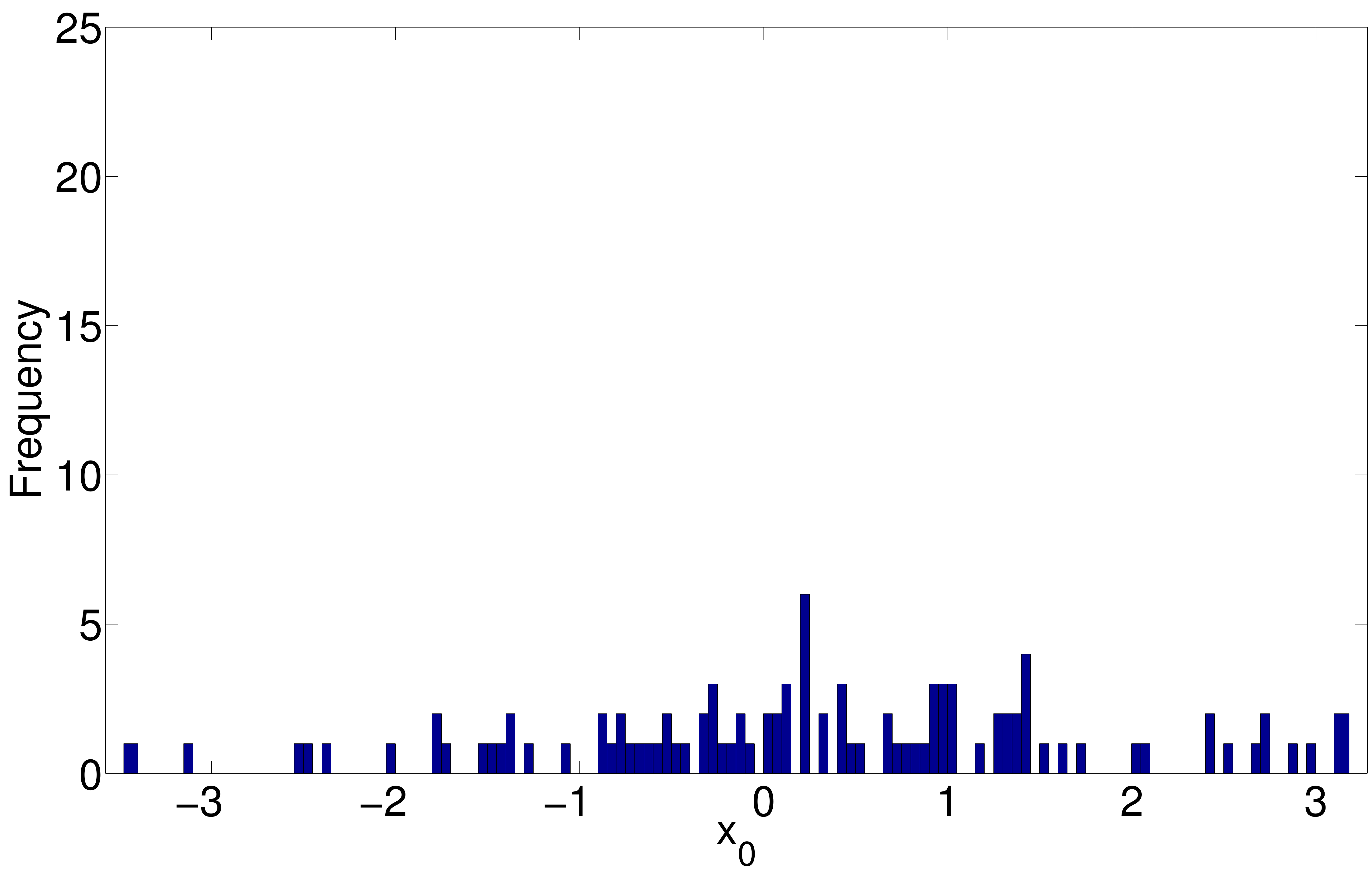}   
}
    \caption{Histograms of the analysis ensembles generated by HMC smoother, and EnKS. The number of ensemble members generated by each smoother is $100$. For the HMC smoother the step of the symplectic integrator~\eqref{eqn:Verlet} is $T=0.1$,  with $h=0.01$ and $m=10$.
}
    \label{fig:Histogram_Poterior_Ensemble}
    \end{figure*}

The histograms of the analysis ensembles obtained with the HMC smoother and EnKS are shown in Figure \ref{fig:Histogram_Poterior_Ensemble}. The HMC smoother generates an analysis ensemble that matches the kernel shown in Figure \ref{fig:oneD_Potential_Posterior_Kernel}, but EnKS fails to generate an accurate analysis ensemble. The most likely state seems to be located in the correct place.
A single analysis state (best estimate) in this case might be misleading. One needs to consider more than one analysis with certain probability to give better description of the true state of the system in case of multi-modal systems.

\begin{figure}[H]
    \centering
  \includegraphics[width=0.95\linewidth,height=3.8cm]{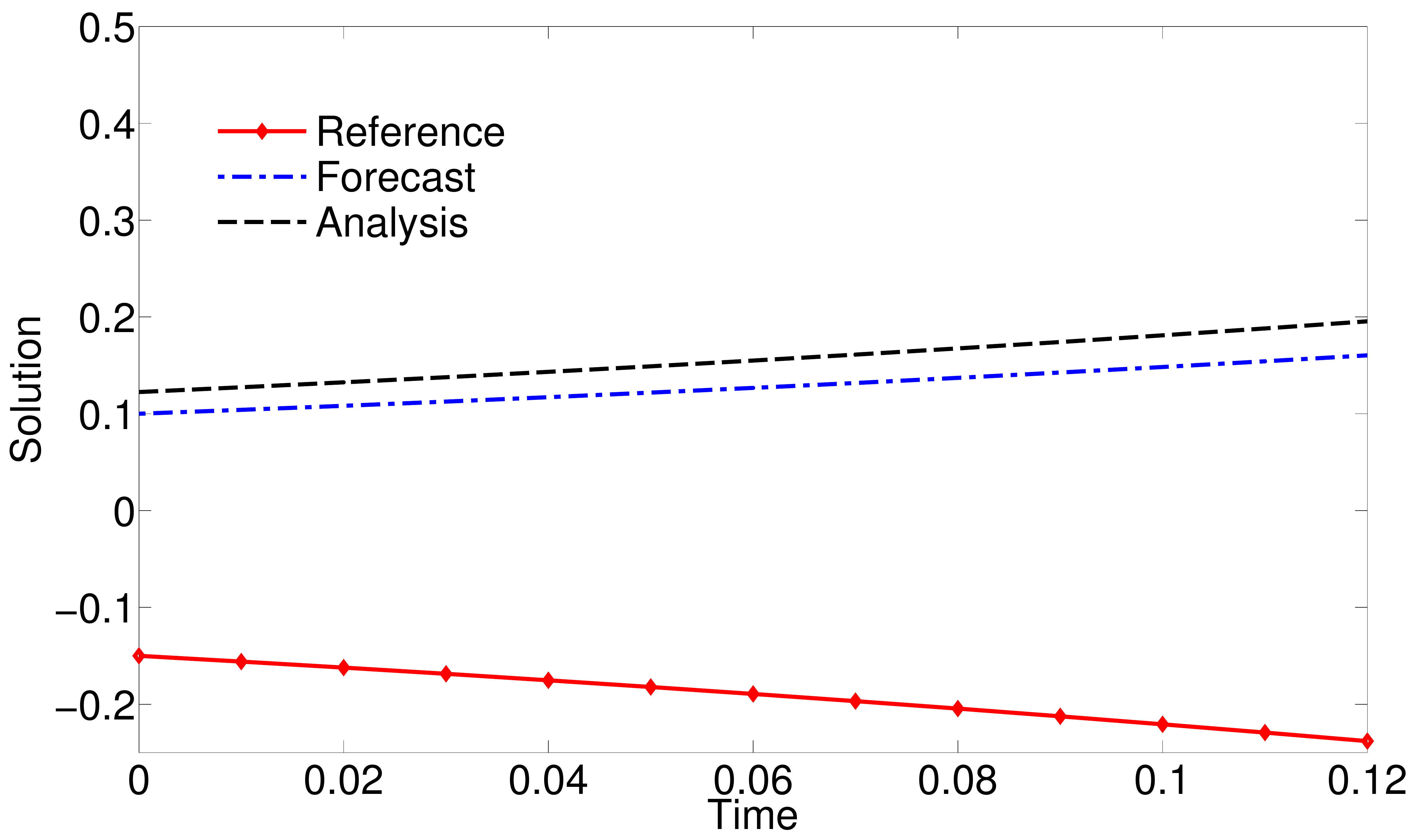}   
    \caption{Data assimilation results using 4D-Var, together with the forecast and reference trajectories plotted over the assimilation window.} 
    \label{fig:4D_Var_Trajectories}
    \end{figure}

The 4D-Var algorithm is expected to be trapped in a local minimum of the posterior distribution. Since the background state is closer to $+1$ than to $-1$, and since the observations \eqref{eqn:Observation_Operator} are insensitive to the sign of solution, we expect the analysis to follow the behavior of the true solution but with the opposite sign. This is confirmed by results in Figure \ref{fig:4D_Var_Trajectories}.
    %

%%%%%%%%%%%%%%%%%%%%%%%%%%%%%%%%%%%%%%%%%%%%%%%
\subsection{Shallow water model on the sphere}
\label{subsec:shallow_water_model}
%%%%%%%%%%%%%%%%%%%%%%%%%%%%%%%%%%%%%%%%%%%%%%%
%
The shallow water equations have been used extensively as a simplified model of the atmosphere that contains the essential wave propagation mechanisms found in general circulation models (GCMs)\cite{Amik:2007}. The shallow water equations in spherical coordinates are  \cite{Navon19911311}:
\begin{subequations}
\label{eqn:SWE_Equations}
\begin{align}
 \frac{\partial u}{\partial t} &+ \frac{1}{a\cos \theta} \left( u \frac{\partial u}{\partial \lambda} 
 + v \cos \theta \frac{\partial u}{\partial \theta} \right)   \\
 &- \left(f + \frac{u \tan \theta}{a} \right) v 
 + \frac{g}{a \cos \theta} \frac{\partial h} {\partial \lambda} = 0, \nonumber \\
  \frac{\partial v}{\partial t} &+ \frac{1}{a\cos \theta} \left( u \frac{\partial v}{\partial \lambda} 
  + v \cos \theta \frac{\partial v}{\partial \theta} \right) \\
  &+ \left(f + \frac{u \tan \theta}{a} \right) u 
  + \frac{g}{a} \frac{\partial h} {\partial \theta} = 0, \nonumber\\
 \frac{\partial h}{\partial t} &+ \frac{1}{a \cos \theta} \left(\frac{\partial\left(hu\right)}{\partial \lambda} 
+ \frac{\partial{\left(hv \cos \theta \right)}}{\partial \theta} \right) = 0\,.
\end{align}
\end{subequations}
Here $f$ is the Coriolis parameter, given by $f = 2 \Omega \sin \theta$, where  $\Omega$ is the angular speed of the rotation of the Earth. In addition, $h$ represents the height of the homogeneous atmosphere, $u$ and $v$ are the zonal and meridional wind components, respectively.  The latitudinal and longitudinal directions are respectively denoted by $\theta$ and $\lambda$.  The radius of the Earth is  denoted by $a$  and $g$ is the acceleration due to gravity. The space discretization is performed using the  unstaggered Turkel-Zwas scheme \cite{navon1987application}. The discretization has nlon=72 nodes in longitudinal direction and  nlat=36 nodes in the latitudinal  direction. The semi-discretization in space leads to a system of ordinary differential equations: 
\begin{equation}
\label{eqn:SWE_FWD_model}
 \x^{\prime} = f(t,\x), \quad \x(t_0) = \x_0; \quad  t_0=0,~t_F=9 \,\textnormal{(hours)}.
\end{equation}
The vector $\x \in \mathbb{R}^n$ with $n=3\times{\rm nlat}\times{\rm nlon}$ contains discrete versions of the zonal wind, meridional wind, and the height variables. We perform the time integration using a $5^{\rm th}$ order Runge-Kutta method. This time-integrator is part of the MATLODE suite~\cite{MATLODE_Tony}, which also has sensitivity analysis capabilities.
%

%%%%%%%%%%%%%%%%%%%%%%%%%%%%%%%%%%
\subsubsection{Observations and background information}
\label{subsec:observations_and_background_info}
%%%%%%%%%%%%%%%%%%%%%%%%%%%%%%%%%%
%    
A reference initial condition is used to generate a reference trajectory. Synthetic linear observations are created from the reference trajectory  by adding Gaussian noise with zero mean  for each of the three components. 
The observation error standard deviation for height component is set to $1.5\%$ of the average magnitude of the reference height component in the reference initial condition. The observation error standard deviation for wind components is set to $10\%$ of the average magnitude of the reference wind component in the initial condition. 
%Since, we assume Gaussian errors here, linearity of the observation operator guarantees that the posterior distribution is Gaussian. 
%This is useful to compare the results of the HMC smoother against the analysis obtained by the 4D-Var scheme.
The initial background state is created by perturbing the reference initial condition with  a Gaussian error drawn from the distribution $\mathcal{N}(0,\mathbf{B}_0)$, with a modeled background error covariance matrix. 
The background error covariance $\mathbf{B}_0$ is modeled  as follows:
\begin{itemize}
   \item Start with a diagonal background error covariance matrix. The standard deviation of the background errors for the height component is $2\%$ of the average magnitude of the reference height component in the 
reference initial condition. The standard deviation of the background errors for the wind components is $15\%$ of the average magnitude of the reference wind component in the reference initial condition. 
   \item Synthetic initial ensemble is created by adding zero-mean 
   Gaussian noise to the reference initial condition with covariances set to the initial (diagonal) background 
   error covariance matrix. Apply the ensemble Kalman filter for $48$ cycles with observations obtained each hour. 
                            The uncertainties in observations are fixed to $1.5\%$, and $10\%$ for the height and wind components respectively.
                            The synthetic observations are obtained by adding Gaussian noise with zero mean and standard deviation equal to the uncertainty level multiplied by the 
                            average magnitude of the corresponding component (height and wind) in the initial condition. 
   \item Decorrelate the ensemble-based covariances using a decorrelation matrix $\rho$ with decorrelation distance $L=1000\,km$.
   \item Calculate $\mathbf{B}_0$ by averaging the ensemble covariances over the last $6$ hours with one matrix per hour. 
                                  In this version the background noise levels are no longer $2\%$ and $15\%$.
\end{itemize}
This method of creating a synthetic initial background error covariance matrix is empirical, but we found that the resulting background error covariance matrix performs well for several algorithms including 4D-Var. Enhancing the quality of this background error covariance matrix can be done by making use of the ensembles generated by the sampling smoother, an idea that we will investigate in future work. 

Data assimilation experiments with this model were conducted for three consecutive assimilation windows. 
The time interval of the first assimilation window is $[0,6]$ hours, the second window is $[6,14]$ hours,  and the third is $[14,20]$ hours.
The short first window can be regarded as a spin-off period for the data assimilation system. Hourly (synthetic) observations are available each of the three windows,  with a total of $6$ observation times in the first window, and $8$ observation times in each of the last two windows.  

Two experiments were conducted. In the first one the background error covariance matrix $\mathbf{B}_0$ is kept fixed for each of the three windows. In the second experiment $\mathbf{B}_0$ is updated with information from the generated ensemble according to the following expression:
\begin{equation}
\label{eqn:Linear_Update_of_B0}
\mathbf{B}_0^{\rm hybrid} = \gamma \times \mathbf{B}_0^{\rm modeled} + (1-\gamma)  \times \mathbf{B}_0^{\rm ensemble}\, ,
\end{equation}
where $\mathbf{B}_0^{\rm hybrid}$ is the updated version of $\mathbf{B}_0$, and $\mathbf{B}_0^{\rm modeled}$ is the fixed version used in the first  experiment. The scalar weight $\gamma$ is a number in the interval $[0,1]$. Selecting $\gamma=1$ ignores the error-of-the-day, while $\gamma=0$ forces the use of only the flow-dependent background error covariance matrix obtained from the ensemble, possibly leading to a singular covariance matrix. In our experiments we chose $\gamma=0.75$. 

The error metric used to compare analyses against the reference solution is the root mean squared error (RMSE):
\begin{equation}
\label{eqn:RMSE_Formula}
   \mathbf{RMSE} = \sqrt{\frac{1}{\nvar} \sum_{i=1}^{\nvar}{(\x_i - \x_i^{\rm true})^2} } \, , 
\end{equation}
where $\x^{\rm true}$ is the reference state of the system. The RMSE is calculated hourly along the trajectory over each assimilation window.

%%%%%%%%%%%%%%%%%%%%%%%%%%%%%%%%%%
\subsubsection{Numerical results with the shallow water on the sphere model}
%%%%%%%%%%%%%%%%%%%%%%%%%%%%%%%%%%

The numerical optimization step in 4D-Var is carried out using the the LBFGS routine implemented in the Poblano optimization toolbox~\cite{dunlavy2010poblano}. 
Here the optimization process is stopped when the norm of the gradient is $1e-10$ or when the relative function tolerance hits $1e-6$. The optimization process takes at least $45$ iterations of LBFGS to converge for the experiment considered here; (see Table \ref{table:SWE_Computational_Cost}) . 

The  HMC sampling smoother is used to generate $100$ ensemble members in each assimilation window. The symplectic integrator used is Verlet \eqref{eqn:Verlet} with an empirically tuned step length of $T_*=0.1$, with $h_*=0.01$, and $m=10$. 
A practically useful recipe is to perturb the step length of the symplectic integrator, a procedure that guarantees that the results obtained are not  contingent on that specific selection of step settings~\cite{blanes2014numerical,neal2011mcmc}. 
The step length $h$ is perturbed with uniform random noise: $h:=(1+r)\times h_*$, $r\sim \mathcal{U}(-0.2,0.2)$.
It is important to notice that the step $h$ is pertubed only once at the beginning of the Hamiltonian trajectory and kept fixed for all the $m$ steps.
                 This actually means that the length of the Hamiltonian trajectory $T$ is perturbed for each proposal state while keeping the number of steps $m$ constant such that at each run the step size $h$ scales accordingly

The number of burn-in steps is set to $30$. We noticed that the HMC smoother converges to the posterior in much fewer steps ($5-10$). Four generated states are discarded between each selected state in the ensemble to guarantee independence of the generated ensemble members.

    \begin{figure}[H]
    \centering
  \subfloat[Experiment with constant $\mathbf{B}_0$]{%
 \includegraphics[width=0.95\linewidth]{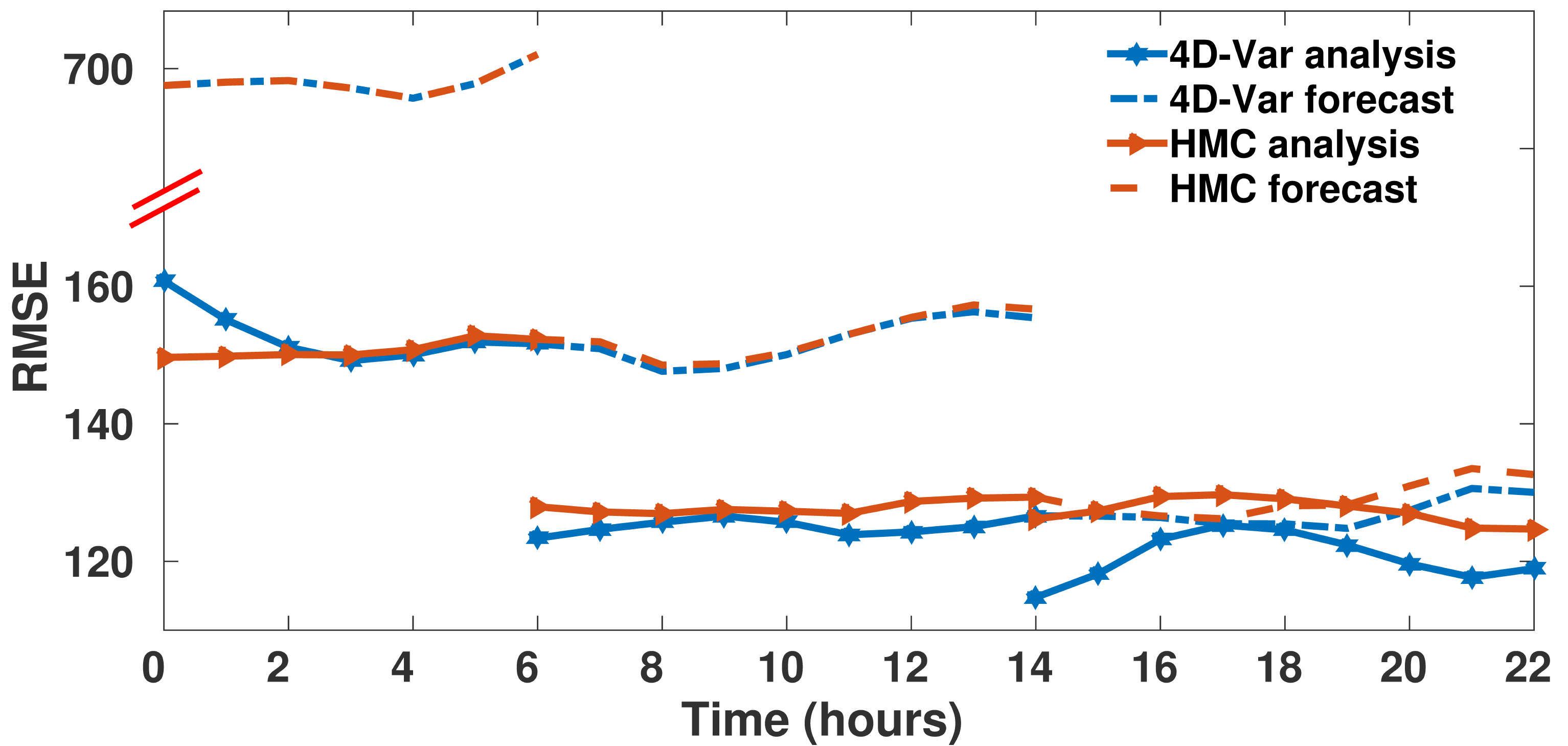}   
 \label{fig:SWE_RMSE_Same_B0}}
  \quad %\hfillht
  \subfloat[Experiment with hybrid $\mathbf{B}_0$ using \eqref{eqn:Linear_Update_of_B0}]{%
 \includegraphics[width=0.95\linewidth]{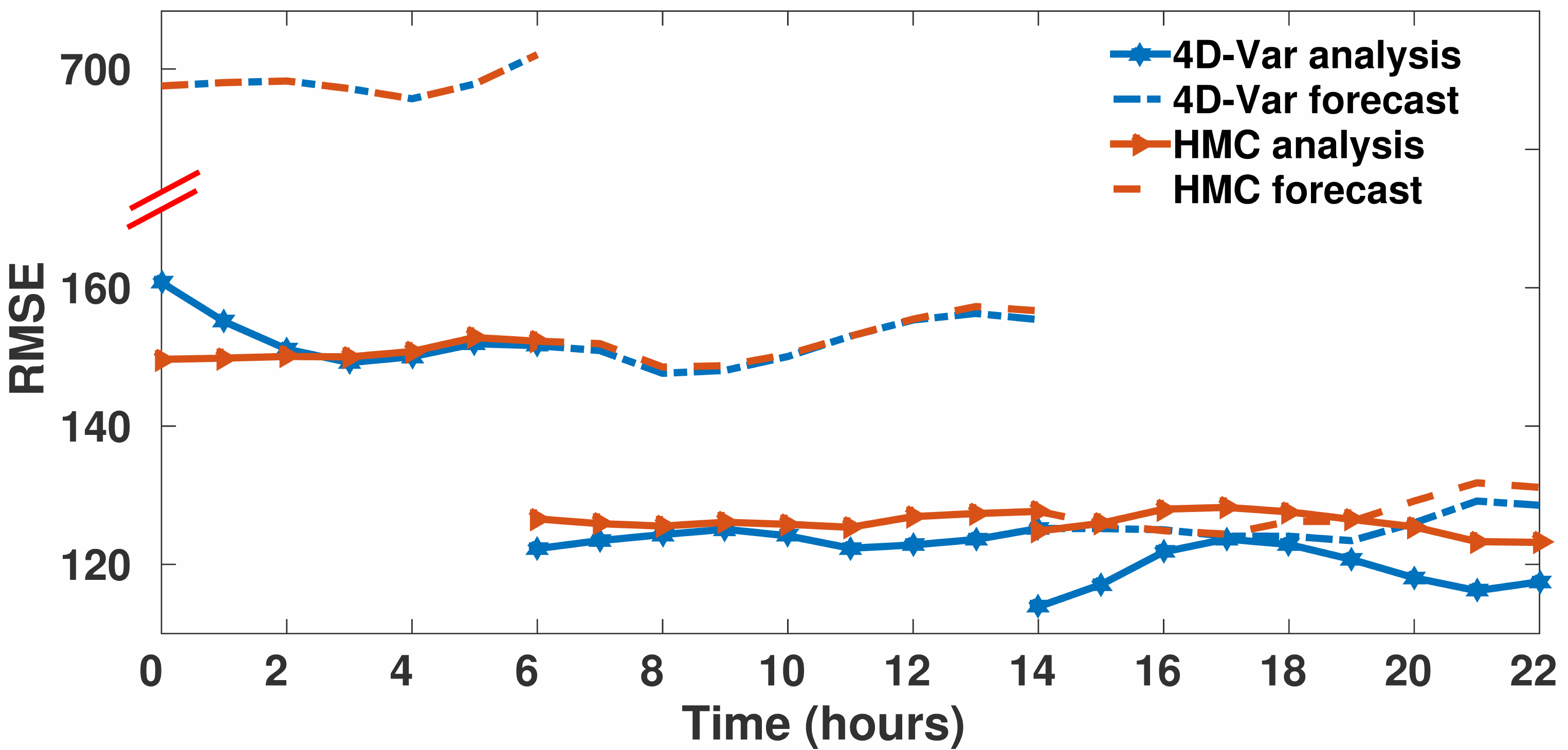}   
 \label{fig:SWE_RMSE_Update_B0}}
  %
    % caption and label of the whole figure
    \caption{Data assimilation results for two scenarios using linear observations are shown. The first panel \ref{fig:SWE_RMSE_Same_B0} shows RMSE with $B_0$ being fixed.
 The second panel \ref{fig:SWE_RMSE_Update_B0} shows RMSE with $B_0$ being updated using \eqref{eqn:Linear_Update_of_B0}.
 The symplectic integrator used in both cases is Verlet \eqref{eqn:Verlet} with step $T=0.1$, where $h=0.01$, and $m=10$. 
 The number of dropped states between selected samples is $4$.} 
    \label{fig:SWE_RMSE}
    \end{figure}    

\begin{figure}[H]
    \centering
          \includegraphics[width=0.95\linewidth]{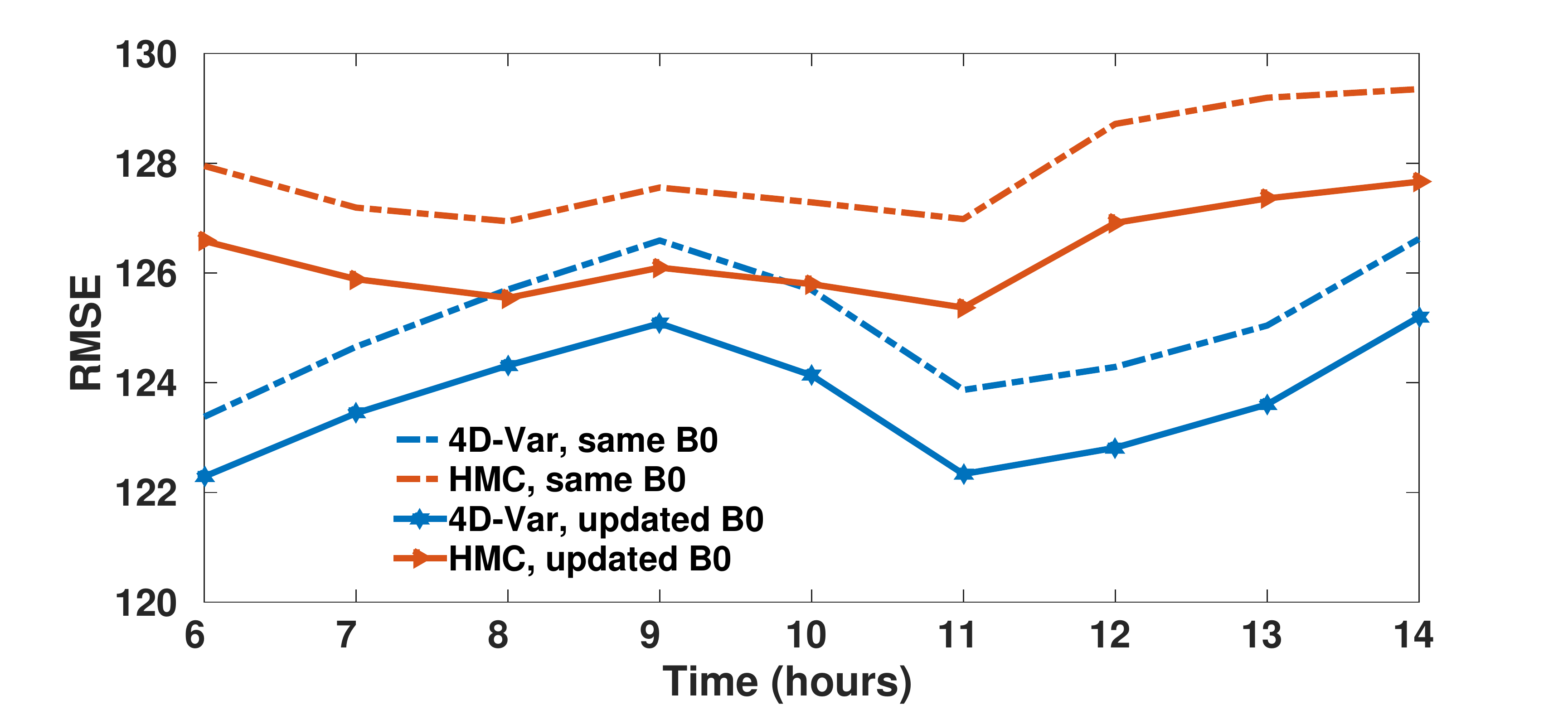}   
    \caption{ Same as results in Figure \ref{fig:SWE_RMSE} with only RMS errors over the second window displayed for two scenarios. 
              RMS errors obtained while keeping $\mathbf{B}_0$ fixed are plotted as dotted lines. 
              RMS errors obtained with $\mathbf{B}_0$ being updated at the beginning of the assimilation window are plotted as dotted lines.  
              } 
    \label{fig:SWE_RMSE_Window_2_Zoomed}
    \end{figure}

\begin{figure*}[htpb]
    \centering
  \subfloat[Reference solution at the initial time, H component]{%
 \includegraphics[width=0.30\linewidth,height=1.90cm]{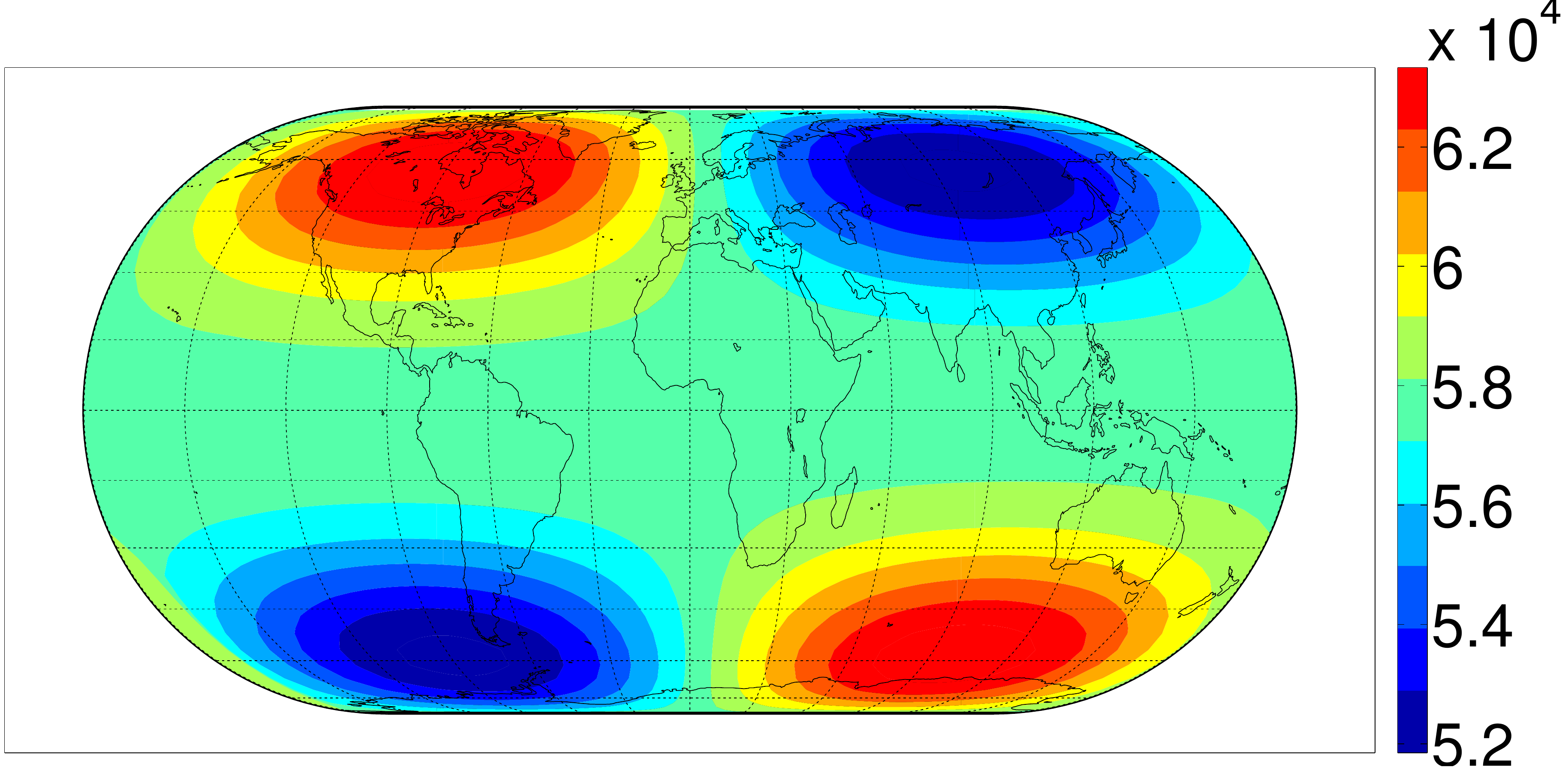}   
 \label{fig:SWE_IC_Ref_H_window_1}}
  \hfill
  \subfloat[Reference solution at the initial time, U component]{%
 \includegraphics[width=0.30\linewidth,height=1.90cm]{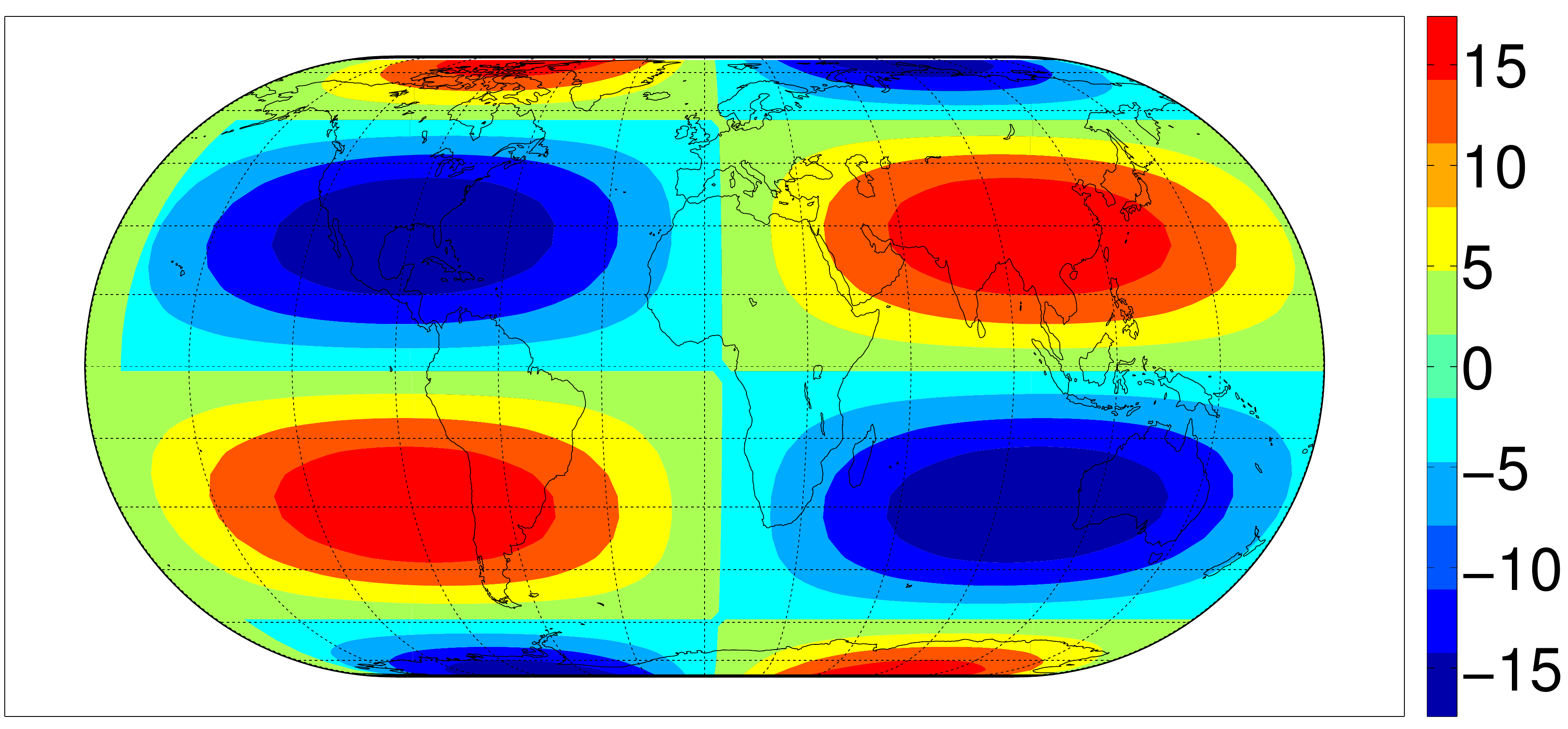}   
 \label{fig:SWE_IC_Ref_U_window_1}}
  \hfill
  \subfloat[Reference solution at the initial time, V component]{%
 \includegraphics[width=0.30\linewidth,height=1.90cm]{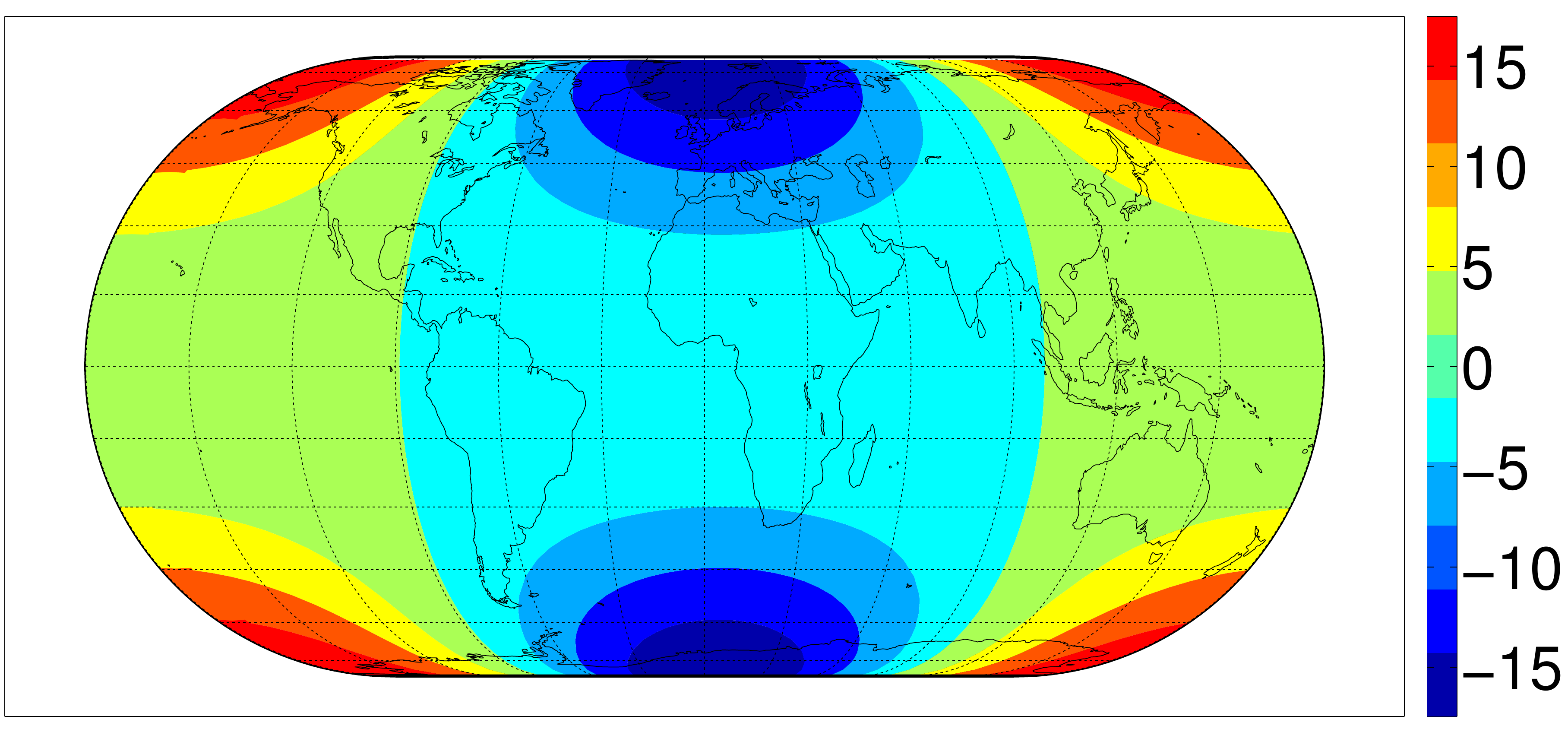}   
 \label{fig:SWE_IC_Ref_V_window_1}}
  \hfill
  \subfloat[Background solution at the initial time, H component]{%
 \includegraphics[width=0.30\linewidth,height=1.90cm]{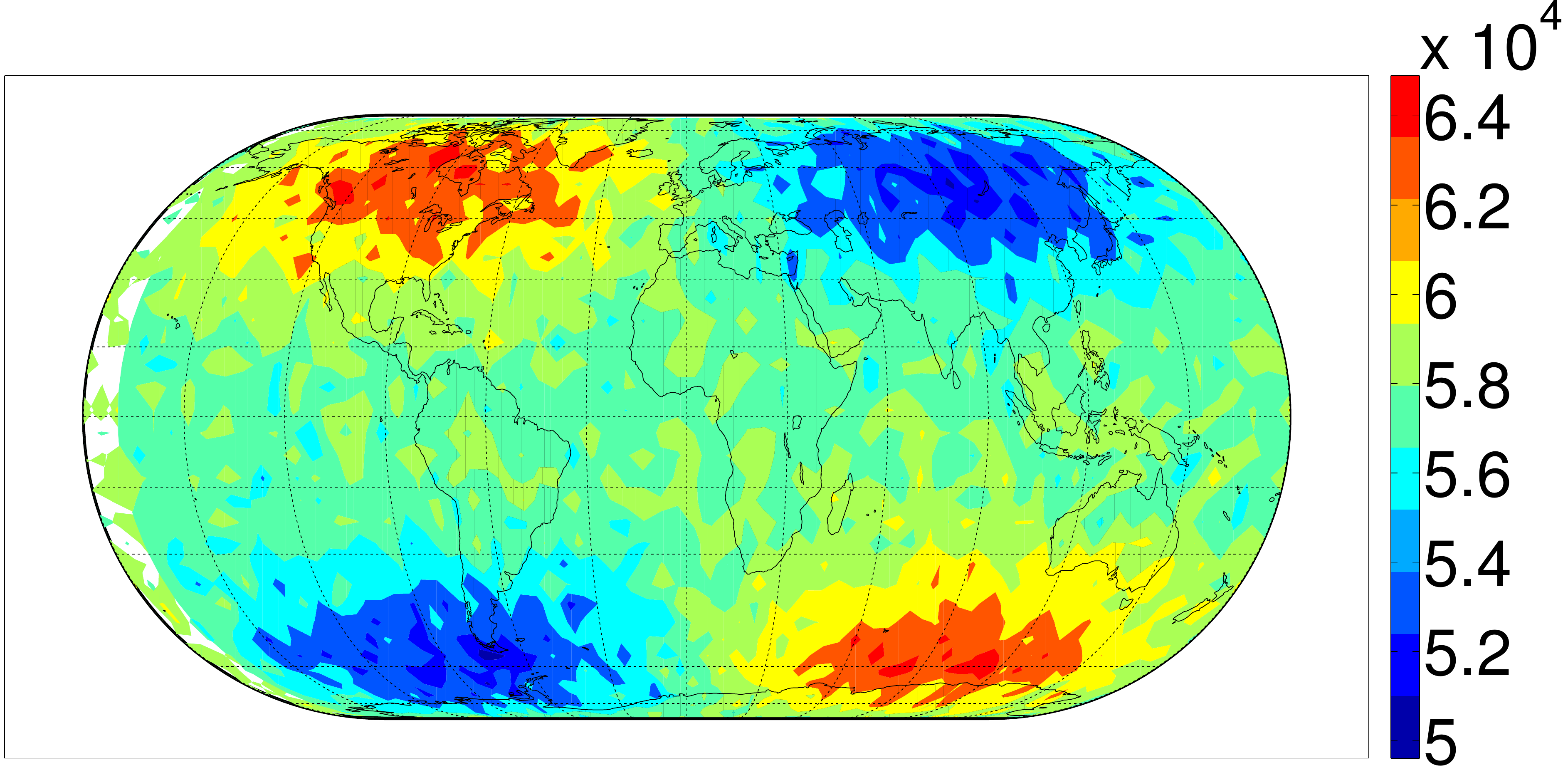}   
 \label{fig:SWE_Forecast_H_window_1}}
  \hfill  
  \subfloat[Background solution at the initial time, U component]{%
 \includegraphics[width=0.30\linewidth,height=1.90cm]{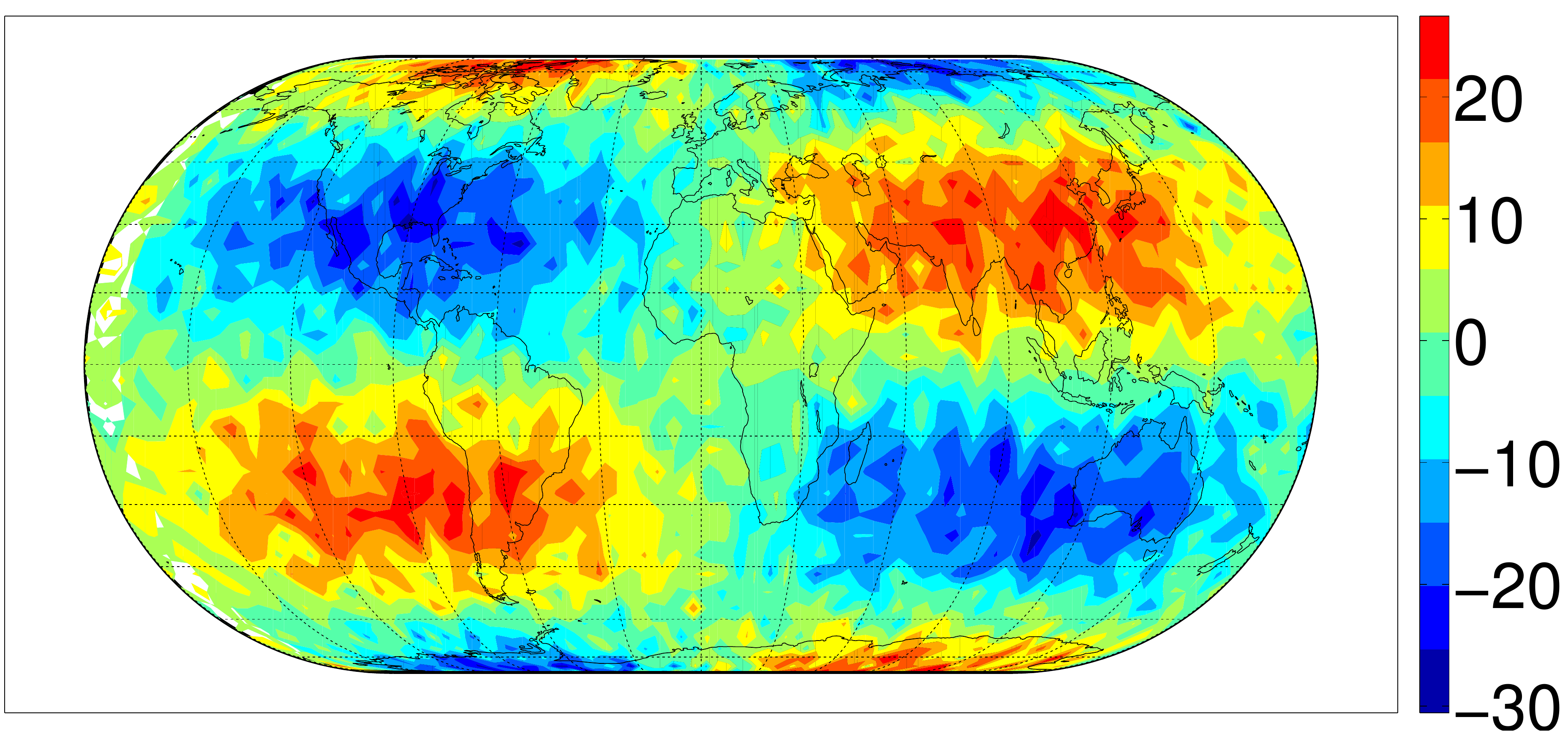}   
 \label{fig:SWE_Forecast_U_window_1}}
  \hfill  
  \subfloat[Background solution at the initial time, V component]{%
 \includegraphics[width=0.30\linewidth,height=1.90cm]{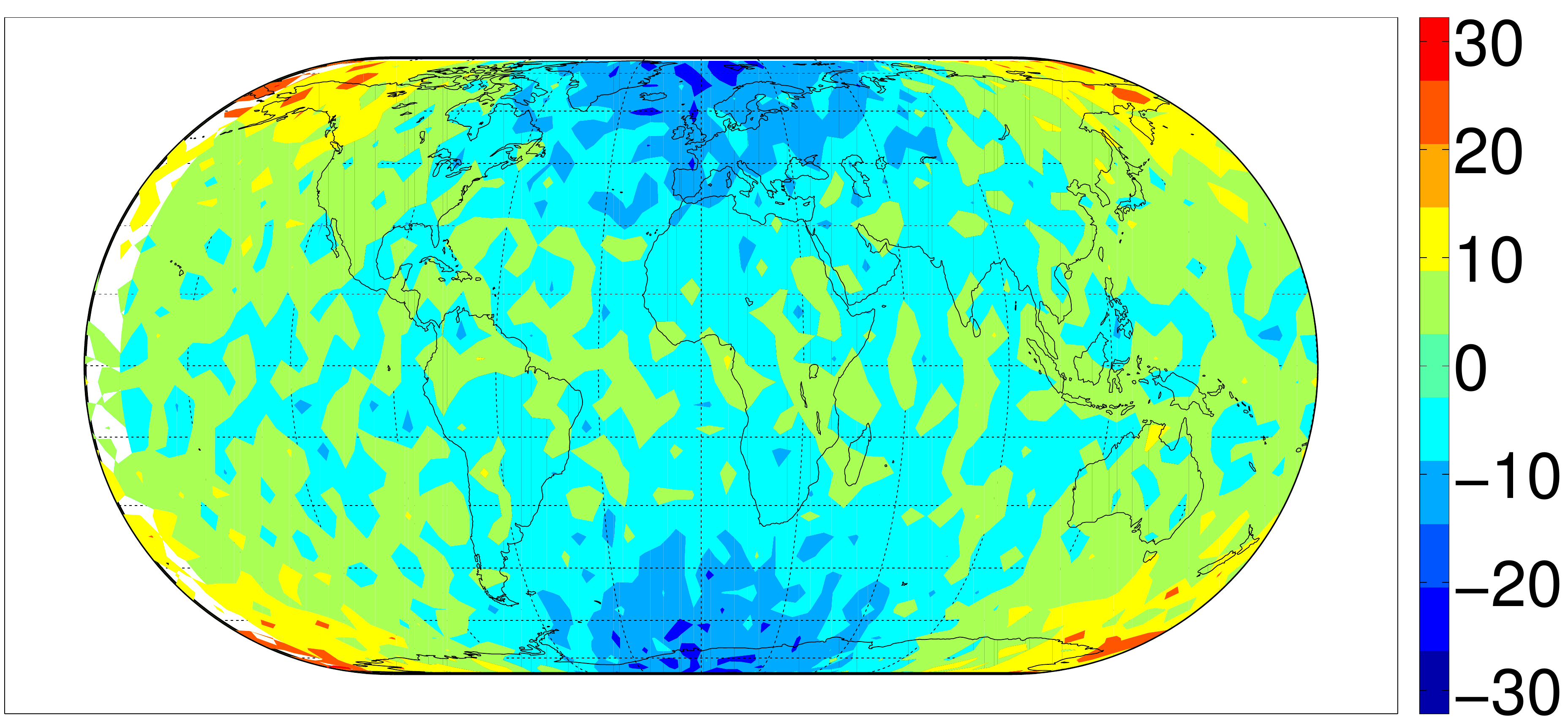}   
 \label{fig:SWE_Forecast_V_window_1}}
  \hfill
  \subfloat[4D-Var analysis at the initial time, H component]{%
 \includegraphics[width=0.30\linewidth,height=1.90cm]{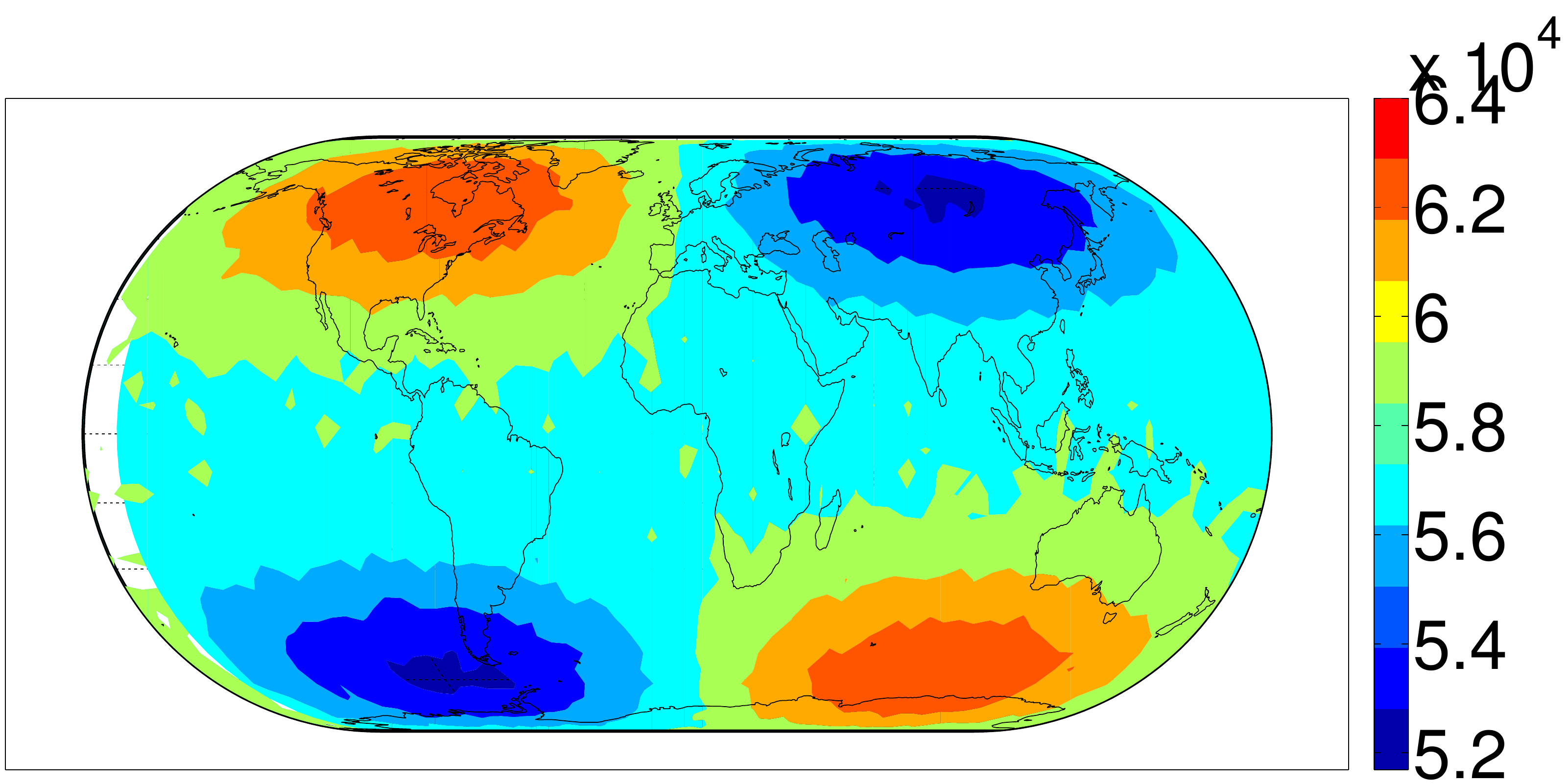}   
 \label{fig:SWE_IC_4DVAR_H_window_1}}
  \hfill
  \subfloat[4D-Var analysis at the initial time, U component]{%
 \includegraphics[width=0.30\linewidth,height=1.90cm]{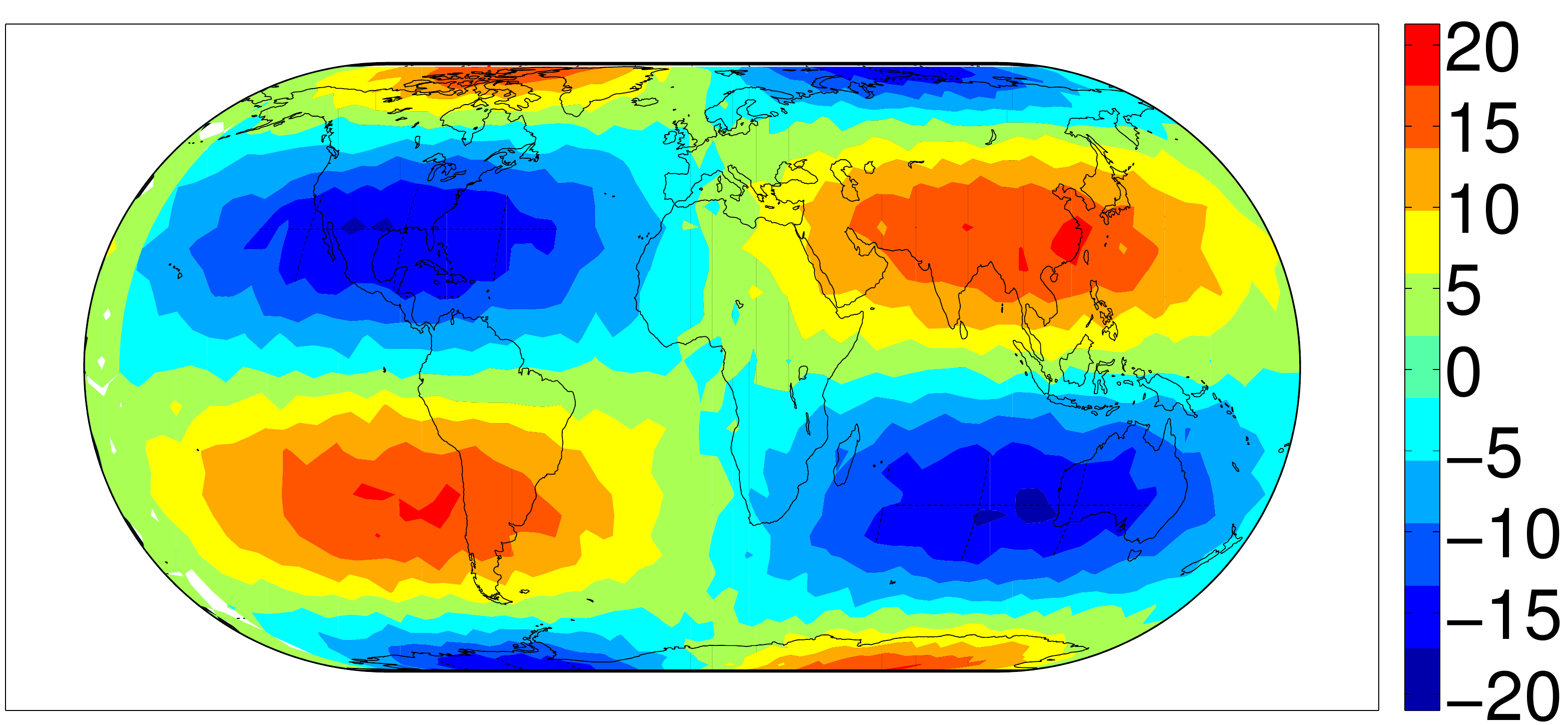}   
 \label{fig:SWE_IC_4DVAR_U_window_1}}
  \hfill
  \subfloat[4D-Var analysis at the initial time, V component]{%
 \includegraphics[width=0.30\linewidth,height=1.90cm]{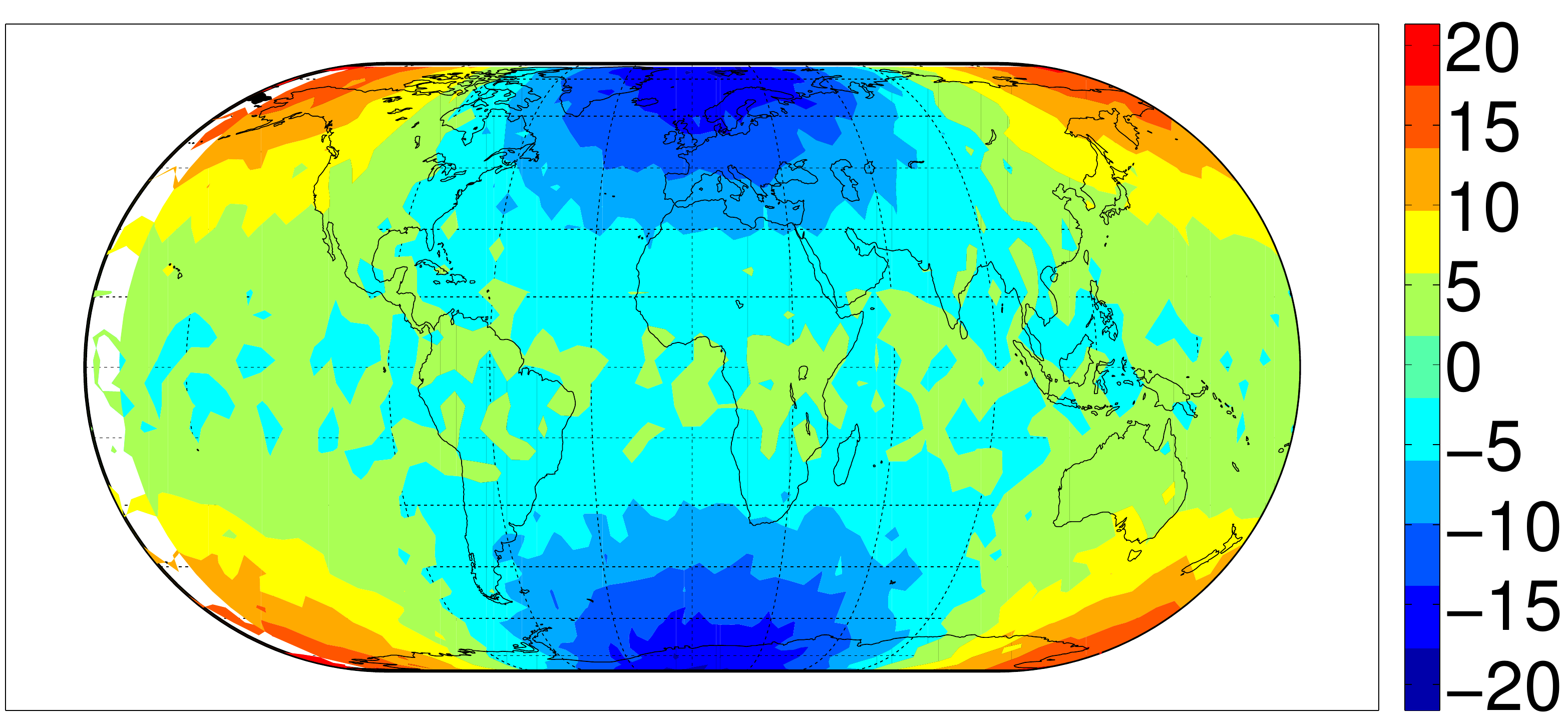}   
 \label{fig:SWE_IC_4DVAR_V_window_1}}
  \hfill
  \subfloat[HMC smoother analysis at the initial time, H component]{%
 \includegraphics[width=0.30\linewidth,height=1.90cm]{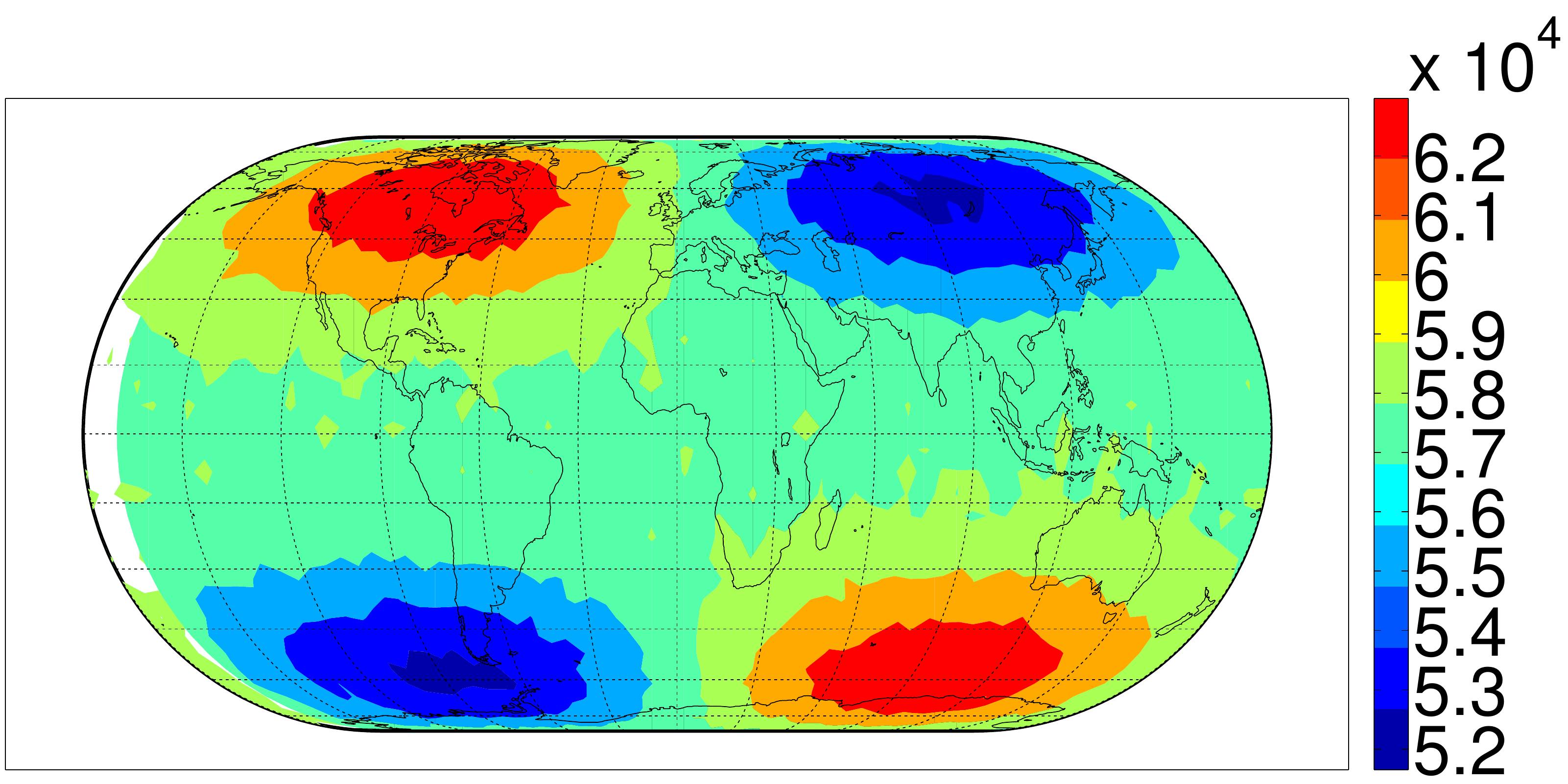}   
 \label{fig:SWE_IC_HMC_H_window_1}}
  \hfill  
  \subfloat[HMC smoother analysis at the initial time, U component]{%
 \includegraphics[width=0.30\linewidth,height=1.90cm]{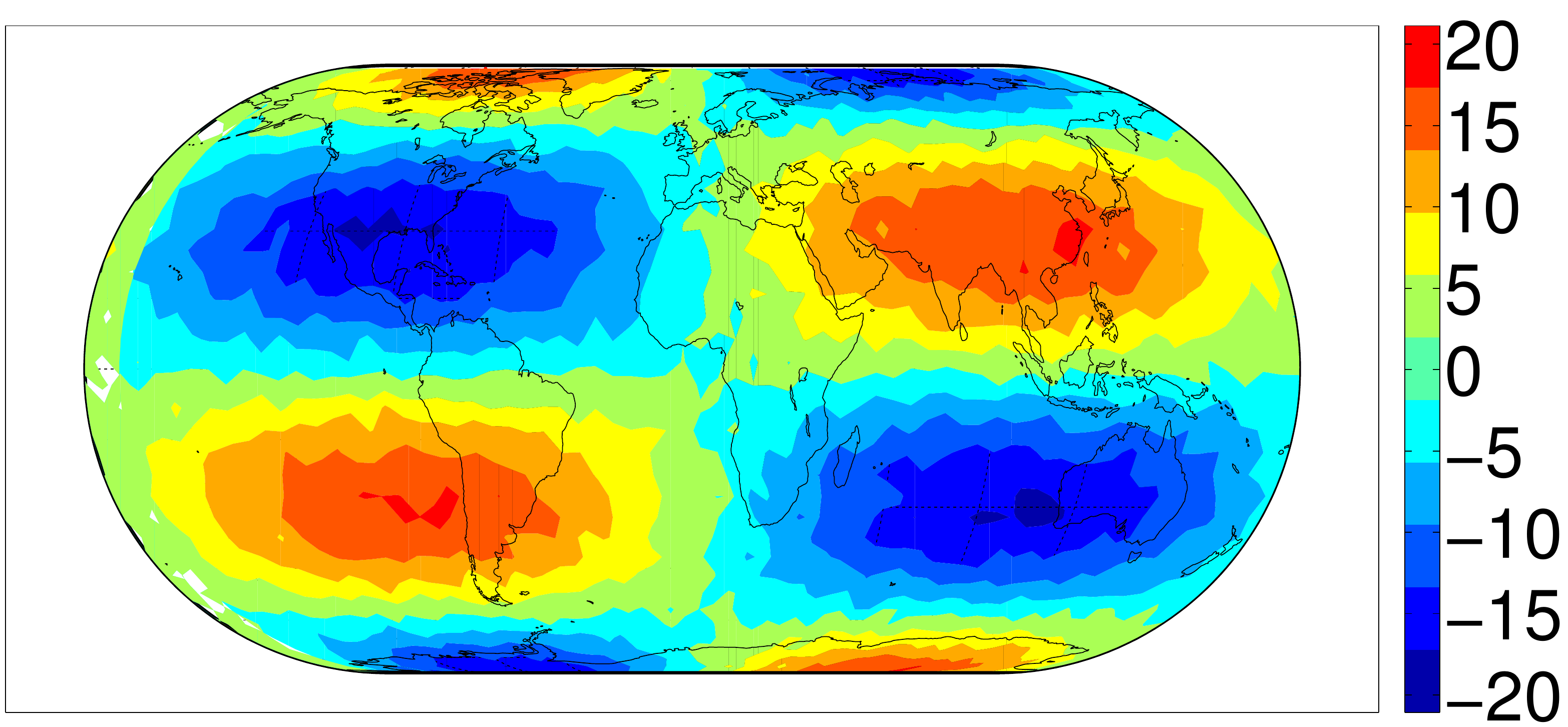}   
 \label{fig:SWE_IC_HMC_U_window_1}}
  \hfill  
  \subfloat[HMC smoother analysis at the initial time, V component]{%
 \includegraphics[width=0.30\linewidth,height=1.90cm]{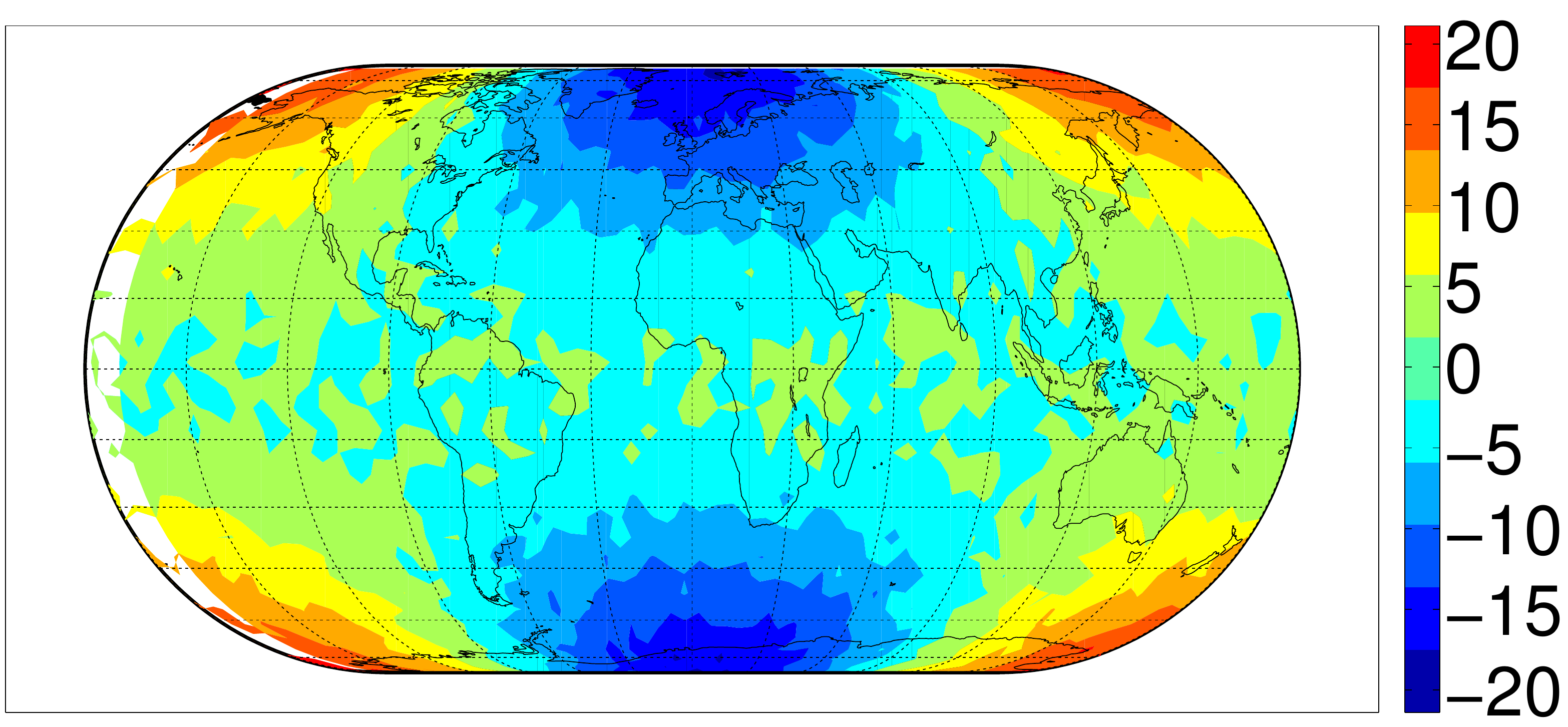}   
 \label{fig:SWE_IC_HMC_V_window_1}}
  %
    % caption and label of the whole figure
    \caption{Four dimensional data assimilation results with linear observations. The initial condition solutions at the beginning of the first window are shown. The data assimilation scheme and the state components are indicated under each panel. The assimilation window length is $6$ hours, with hourly observations. The background error covariance matrix $\mathbf{B}_0$ is kept fixed.
 } 
    \label{fig:SWE_IC_Solutions_Window_1}
    \end{figure*}

The RMS errors for both HMC smoother and 4D-Var over the three assimilation windows are shown in Figure \ref{fig:SWE_RMSE}. Figure \ref{fig:SWE_RMSE_Same_B0} reports the case where the background error covariance matrix $\mathbf{B}_0$ is kept fixed, and Figure \ref{fig:SWE_RMSE_Update_B0} shows the case where $\mathbf{B}_0$ is updated, at the beginning of each assimilation window, according to equation \eqref{eqn:Linear_Update_of_B0}. 
The quality of the analyses in both cases is very similar, however in the second case a slight reduction in RMSE is noticed along the entire trajectory. This is appreciated by Figure \ref{fig:SWE_RMSE_Window_2_Zoomed} zooming onto the RMSE results over the second assimilation window.
    The HMC smoother can sample efficiently from the posterior distribution and resulting analysis competes in accuracy with that obtained using 4D-Var.
    Figure \ref{fig:SWE_IC_Solutions_Window_1} shows the three components at the beginning of the first window for the reference solution, the background state, 4D-Var analysis, and HMC smoother analysis (ensemble mean). Results shown in The analysis recovered from  the noisy background by both 4D-Var and HMC smoother are almost identical.

\begin{figure*}[htpb!]
    \centering
  \subfloat[Reference solution at the initial time, H component]{%
 \includegraphics[width=0.30\linewidth,height=1.90cm]{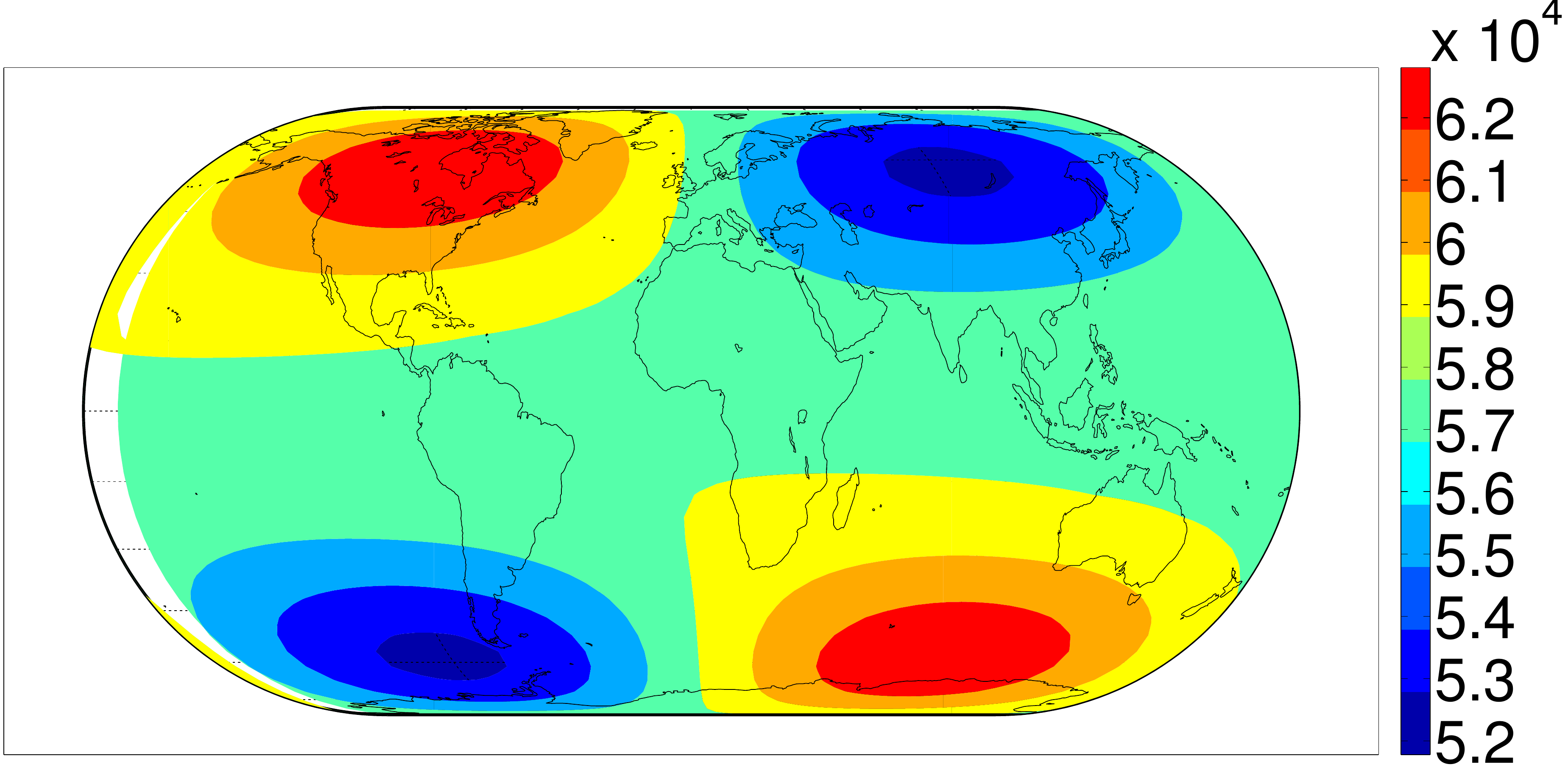}   
 \label{fig:SWE_IC_Ref_H_Same_B0_window_2}}
  \hfill
  \subfloat[Reference solution at the initial time, U component]{%
 \includegraphics[width=0.30\linewidth,height=1.90cm]{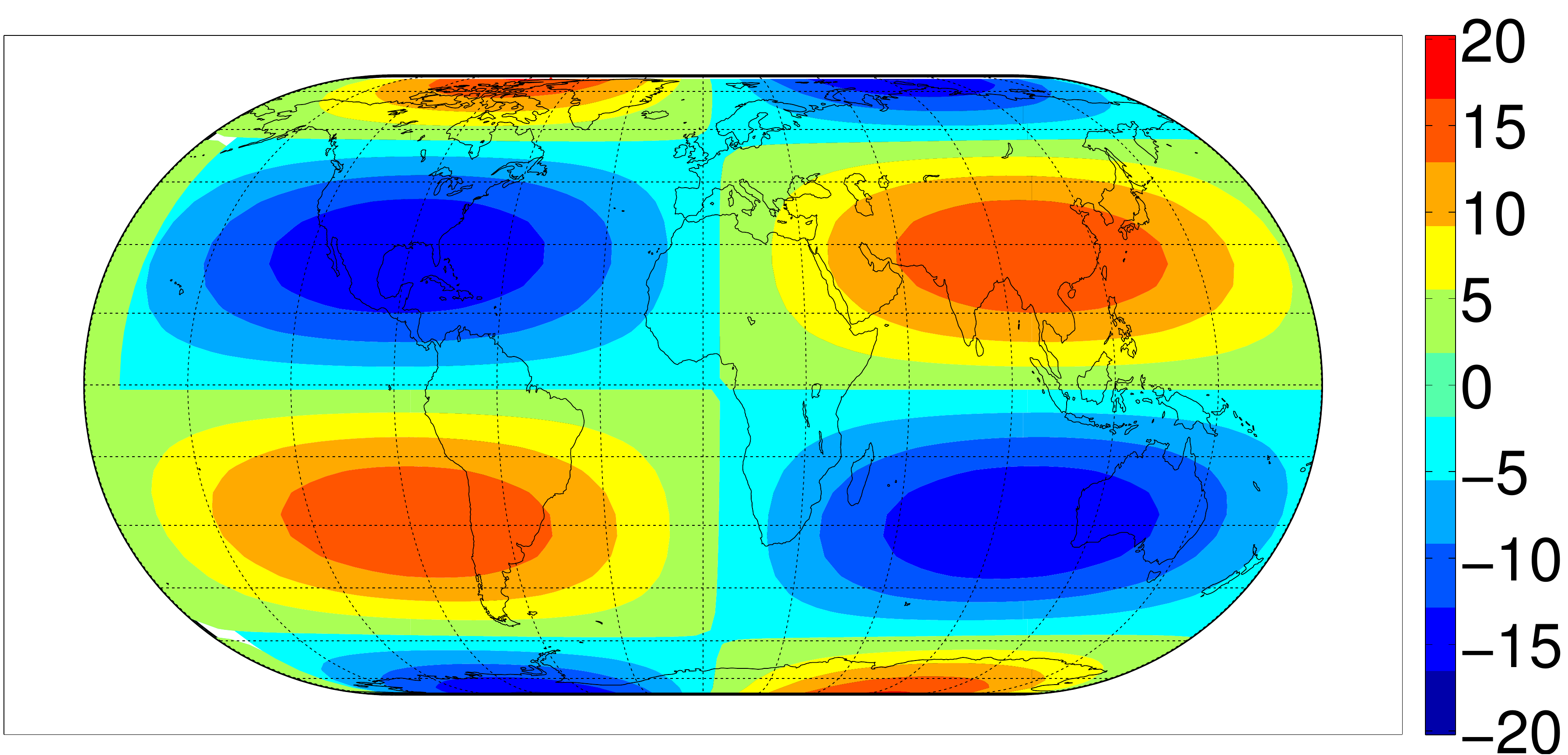}   
 \label{fig:SWE_IC_Ref_U_Same_B0_window_2}}
  \hfill
  \subfloat[Reference solution at the initial time, V component]{%
 \includegraphics[width=0.30\linewidth,height=1.90cm]{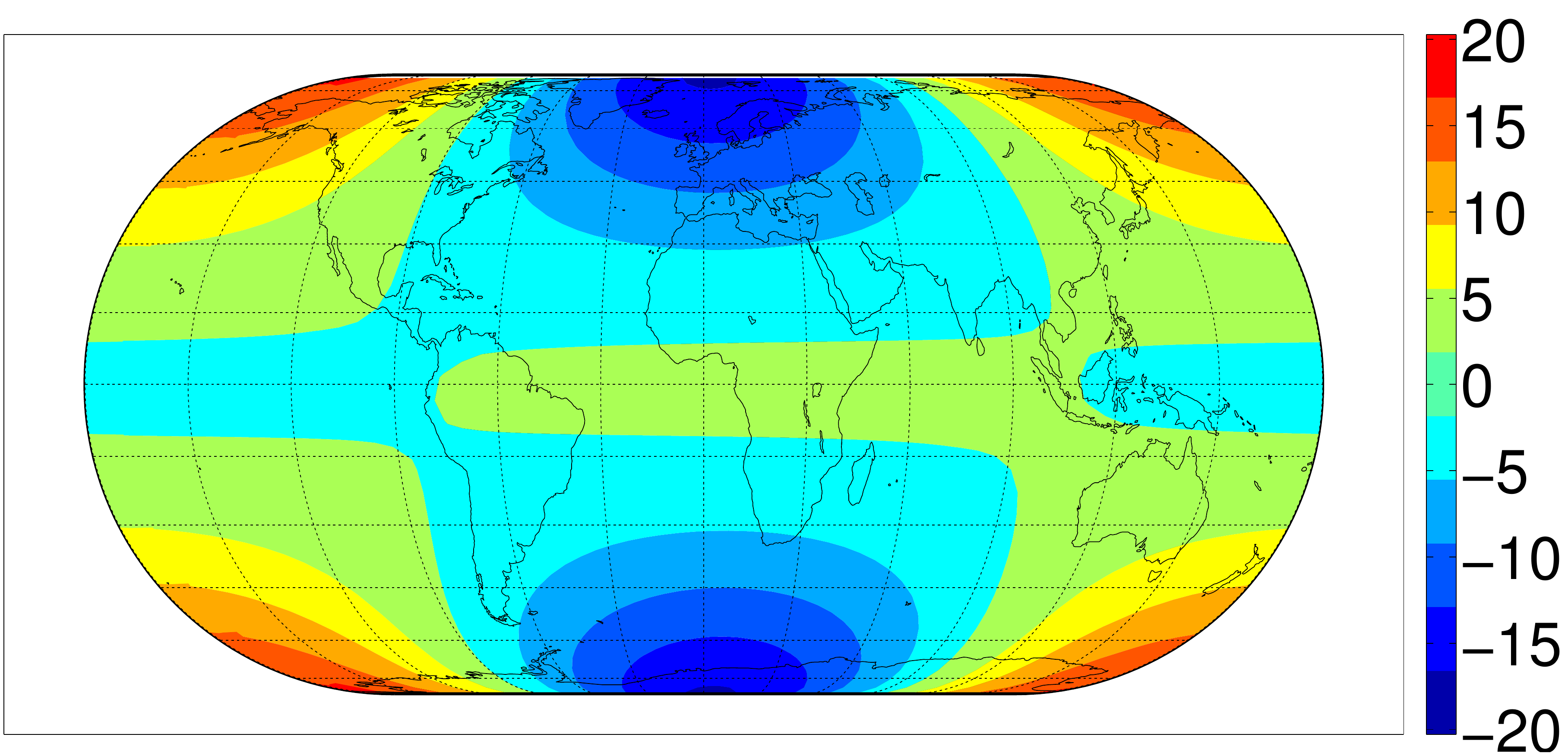}   
 \label{fig:SWE_IC_Ref_V_Same_B0_window_2}}
  \hfill
  \subfloat[HMC smoother analysis at the initial time, H component]{%
 \includegraphics[width=0.30\linewidth,height=1.90cm]{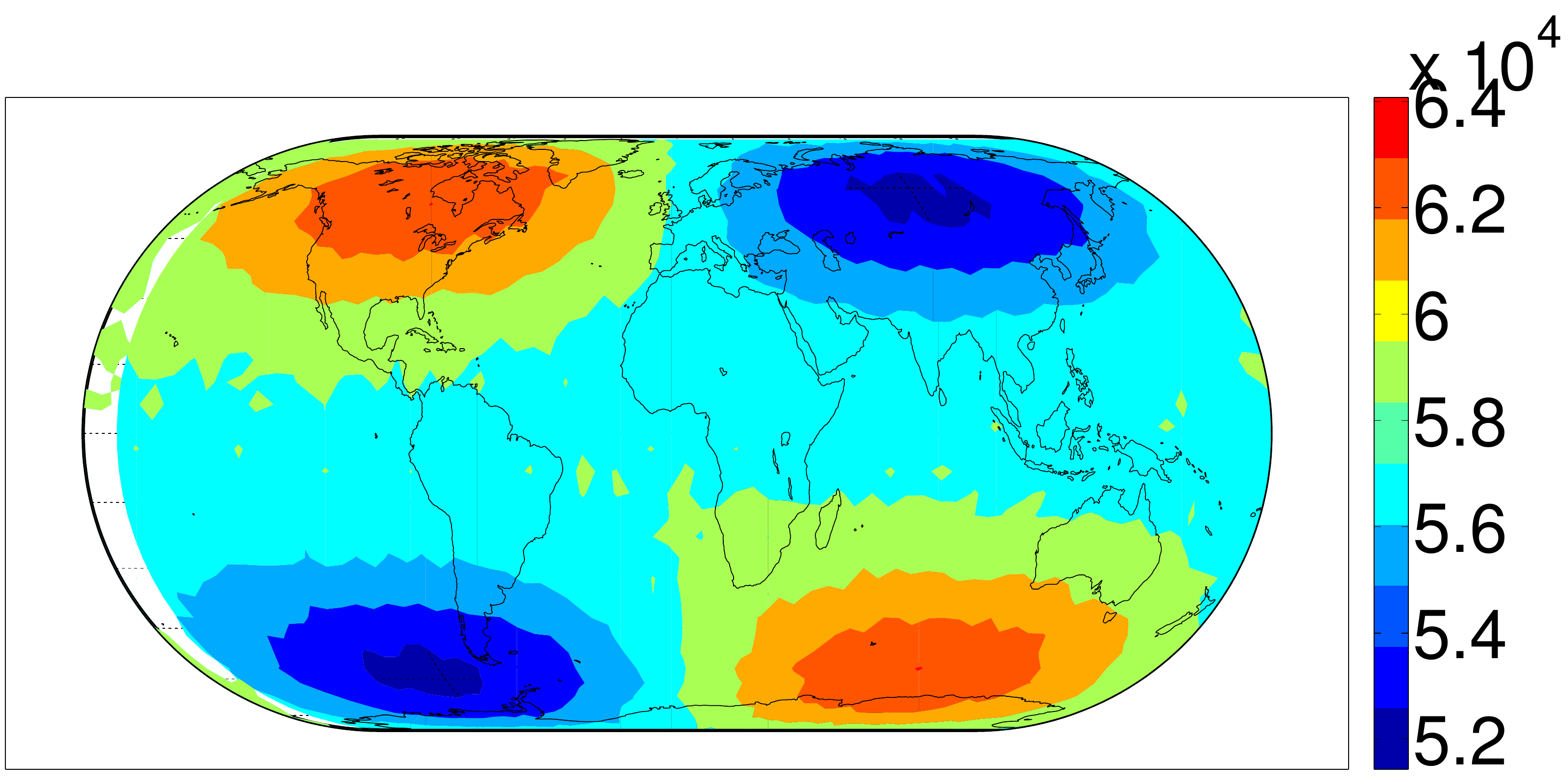}   
 \label{fig:SWE_IC_HMC_H_Same_B0_window_2}}
  \hfill  
  \subfloat[HMC smoother analysis at the initial time, U component]{%
 \includegraphics[width=0.30\linewidth,height=1.90cm]{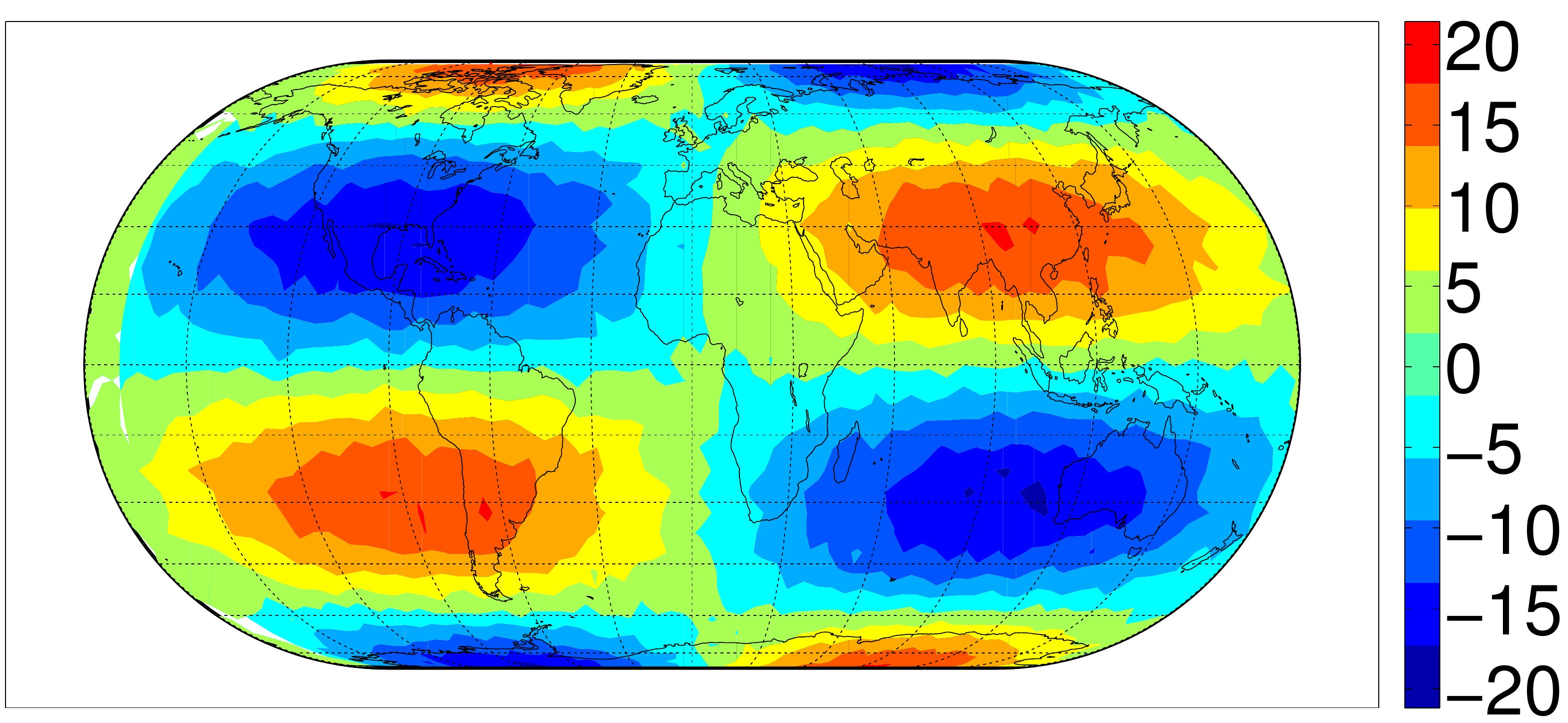}   
 \label{fig:SWE_IC_HMC_U_Same_B0_window_2}}
  \hfill  
  \subfloat[HMC smoother analysis at the initial time, V component]{%
 \includegraphics[width=0.30\linewidth,height=1.90cm]{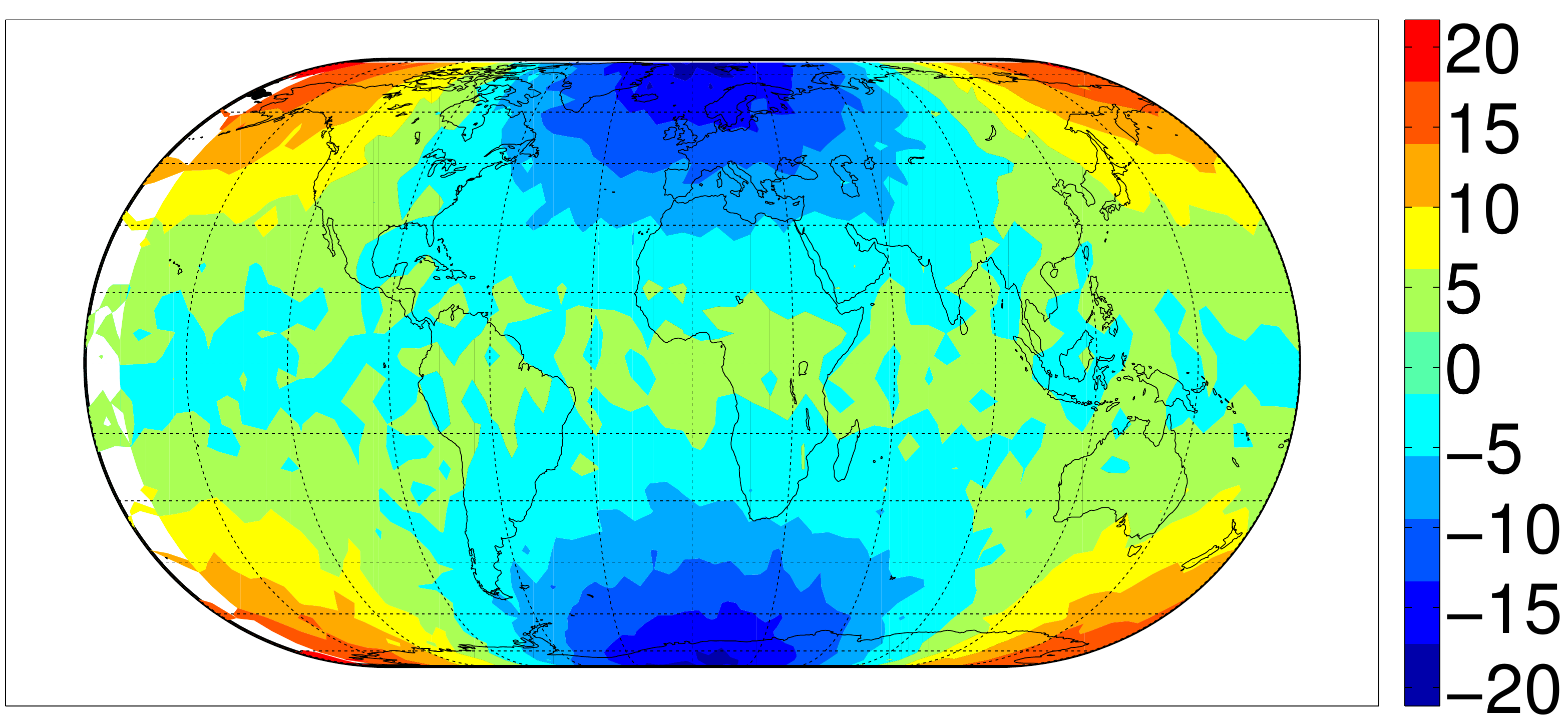}   
 \label{fig:SWE_IC_HMC_v_Same_B0_window_2}}
  \hfill
  \subfloat[4D-Var analysis at the initial time, H component]{%
 \includegraphics[width=0.30\linewidth,height=1.90cm]{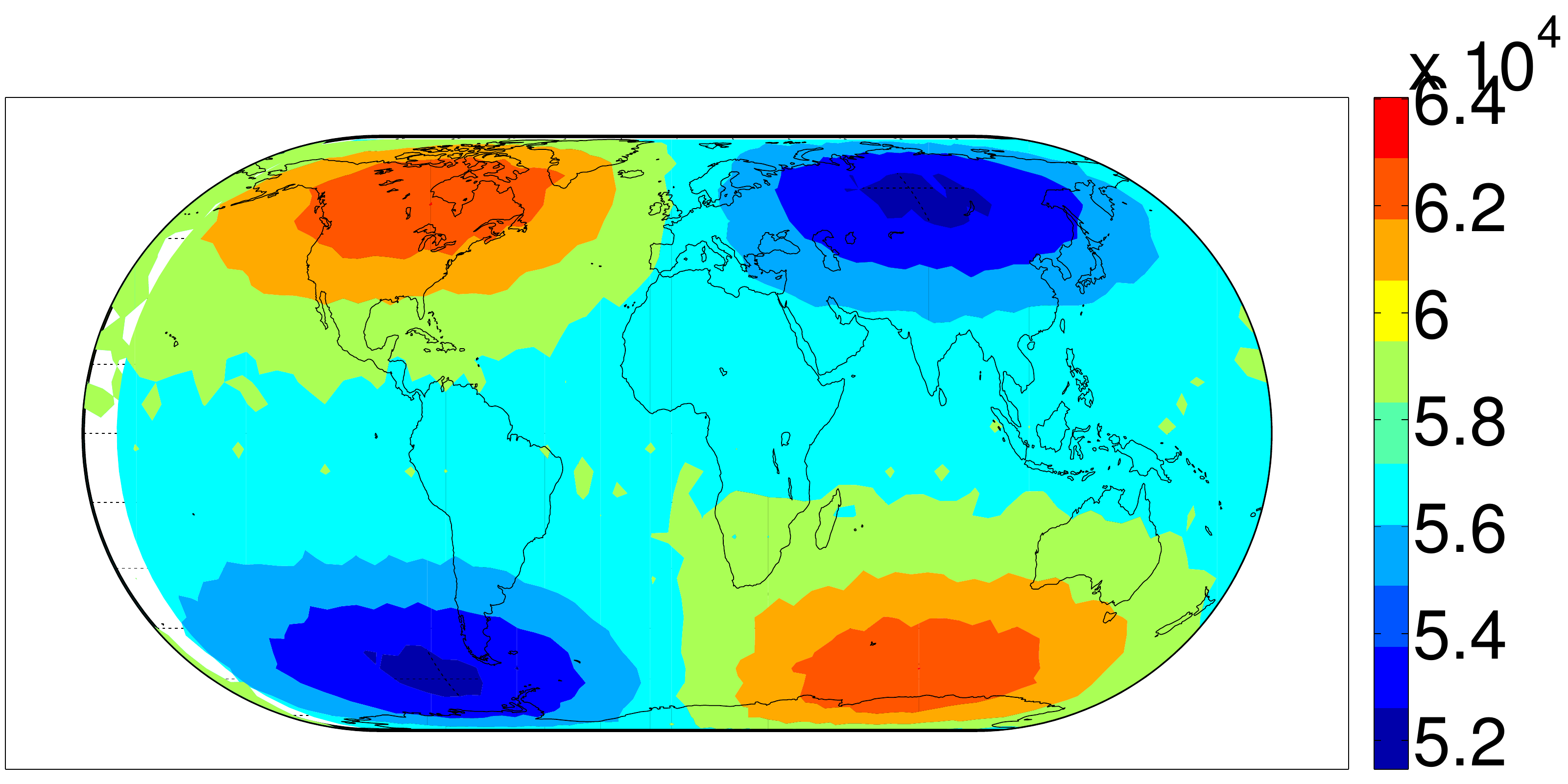}   
 \label{fig:SWE_IC_4DVAR_H_Same_B0_window_2}}
  \hfill
  \subfloat[4D-Var analysis at the initial time, U component]{%
 \includegraphics[width=0.30\linewidth,height=1.90cm]{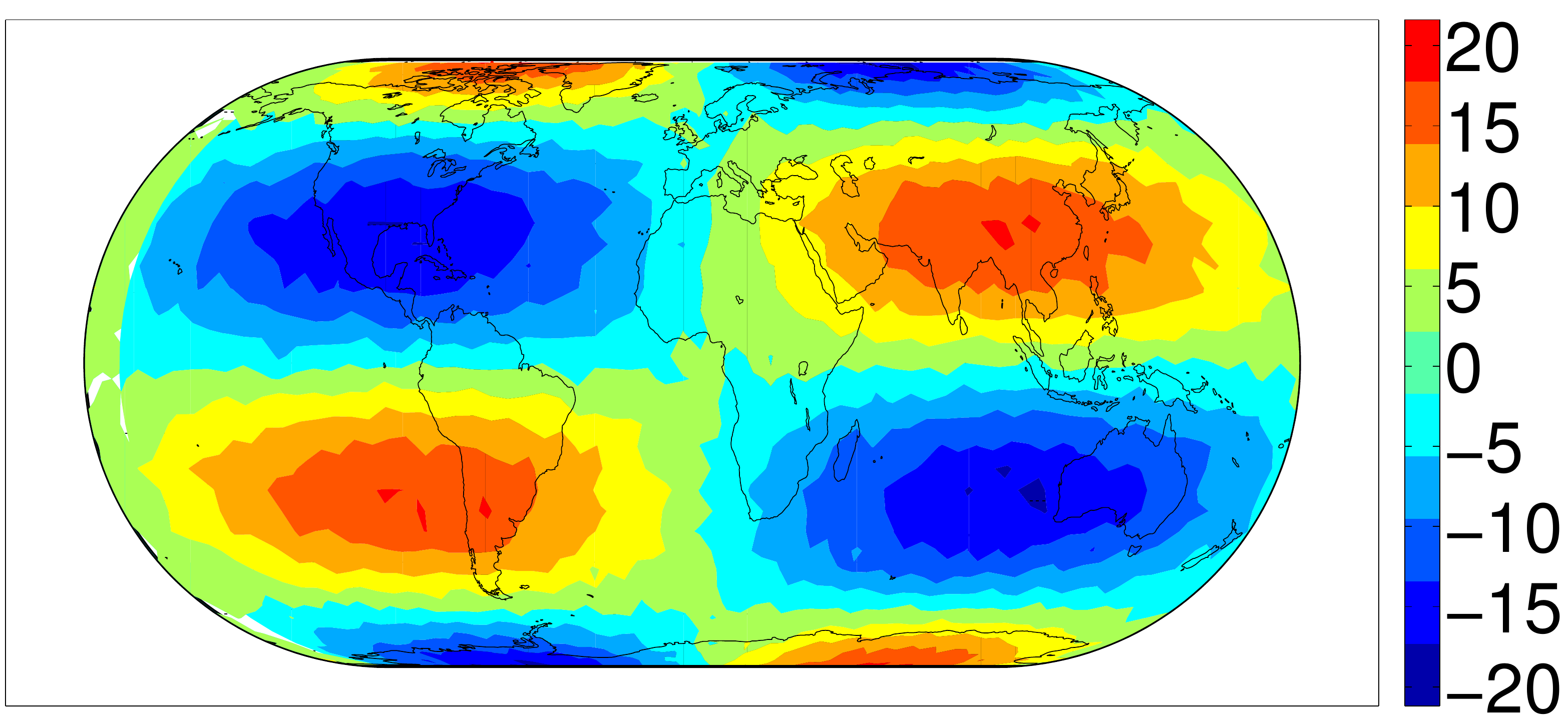}   
 \label{fig:SWE_IC_4DVAR_U_Same_B0_window_2}}
  \hfill  
  \subfloat[4D-Var analysis at the initial time, V component]{%
 \includegraphics[width=0.30\linewidth,height=1.90cm]{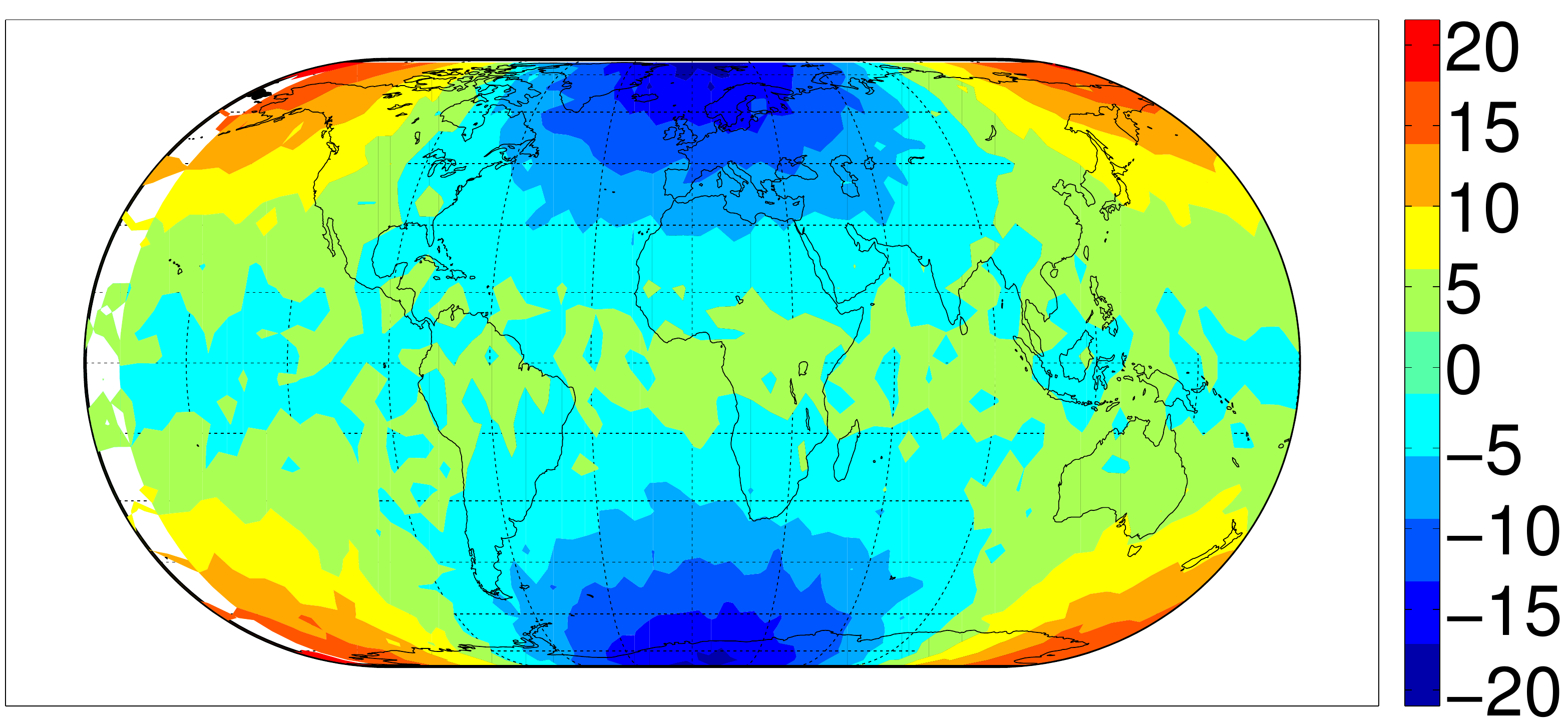}   
 \label{fig:SWE_IC_4DVAR_V_Same_B0_window_2}}
  %
    % caption and label of the whole figure
    \caption{Four dimensional data assimilation results with linear observations. The initial condition solutions at the beginning of the second window are shown.  The data assimilation scheme and the state components are indicated under each panel. The assimilation window length is $8$ hours, with hourly observations.
The background error covariance matrix $\mathbf{B}_0$ is not updated.} 
    \label{fig:SWE_IC_Solutions_Same_B0_Window_2}
    \end{figure*}

\begin{figure*}[htpb!]
    \centering
  \subfloat[Reference solution at the initial time, H component]{%
 \includegraphics[width=0.30\linewidth,height=1.90cm]{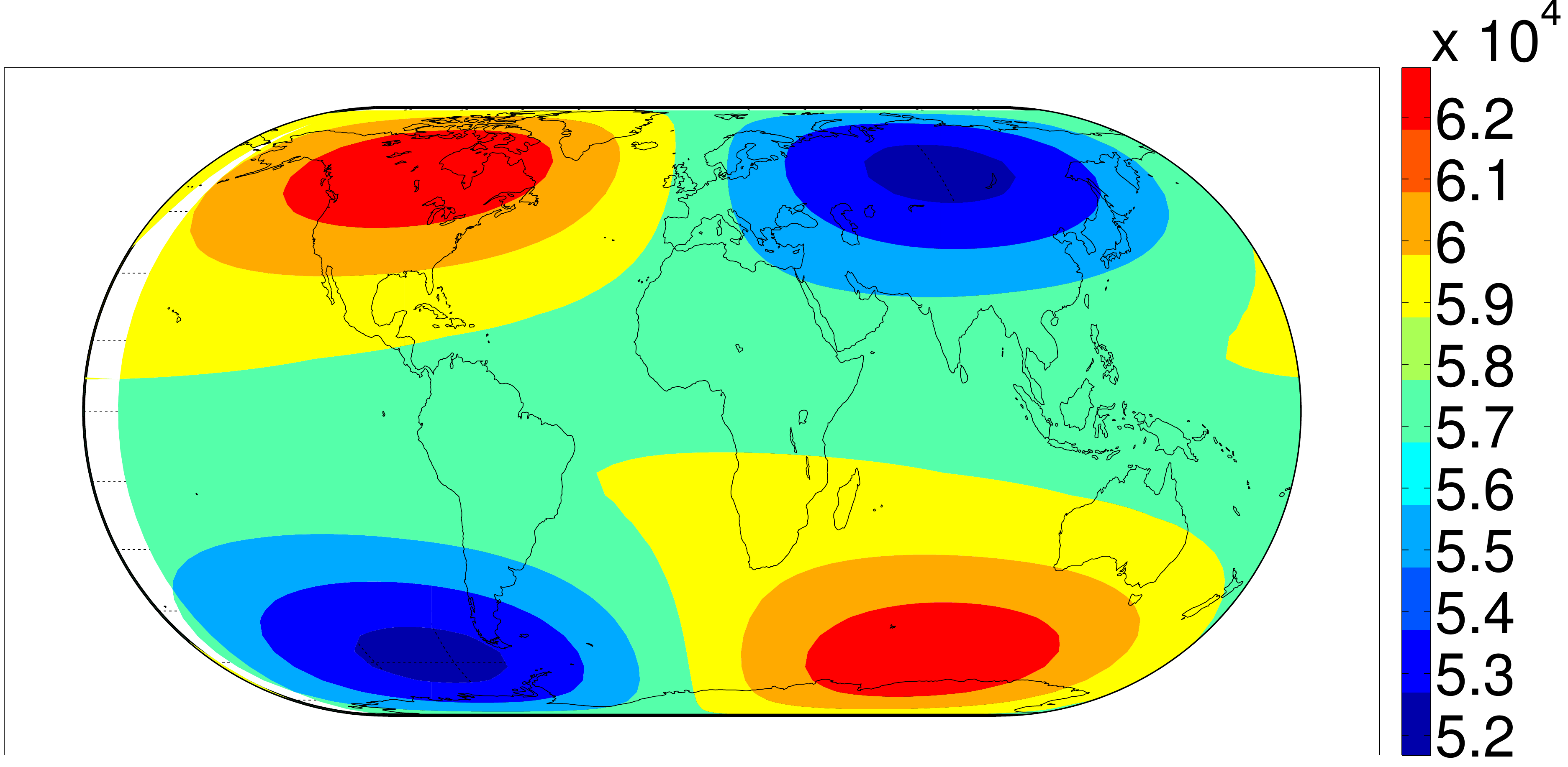}   
 \label{fig:SWE_IC_Ref_H_Same_B0_window_3}}
  \hfill
  \subfloat[Reference solution at the initial time, U component]{%
 \includegraphics[width=0.30\linewidth,height=1.90cm]{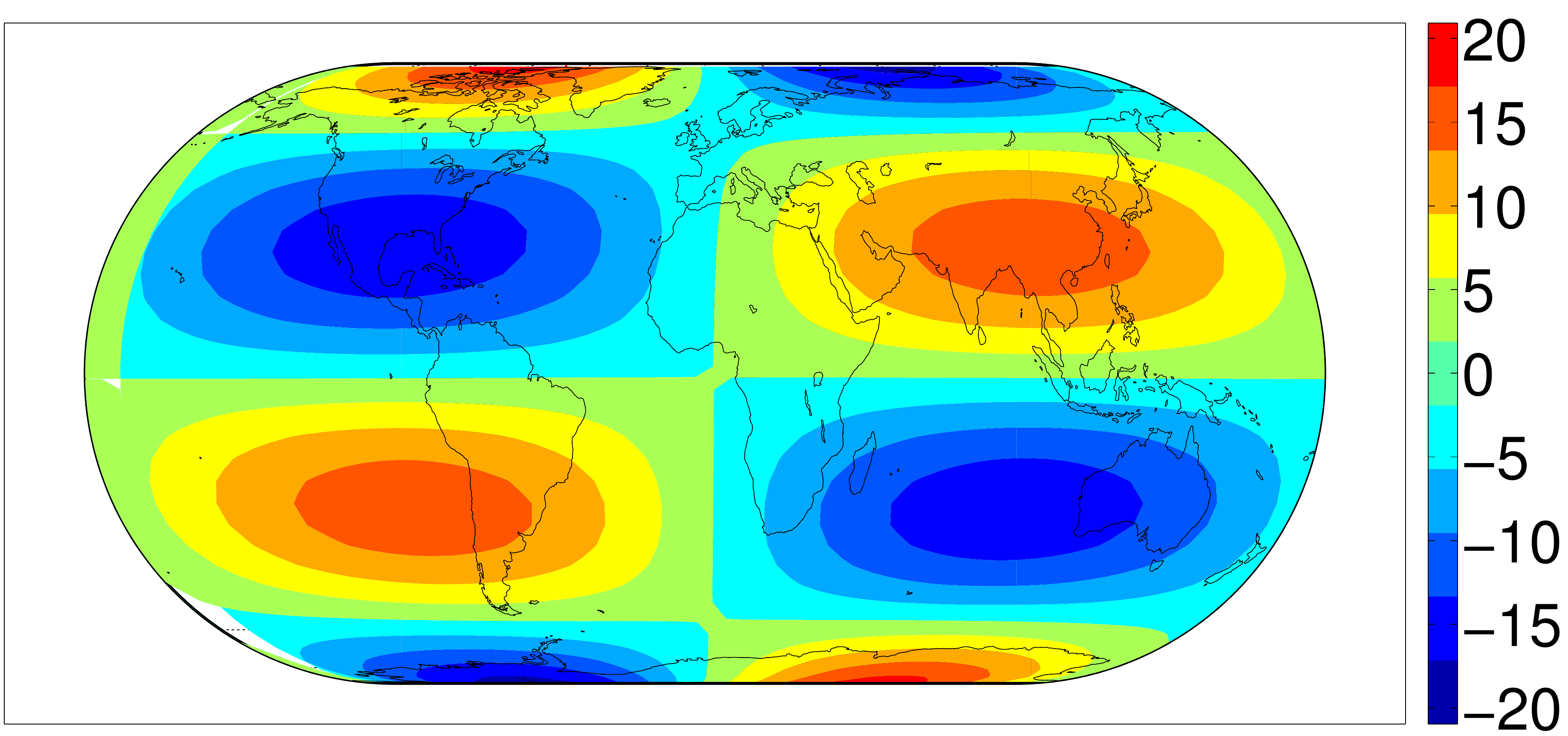}   
 \label{fig:SWE_IC_Ref_U_Same_B0_window_3}}
  \hfill
  \subfloat[Reference solution at the initial time, V component]{%
 \includegraphics[width=0.30\linewidth,height=1.90cm]{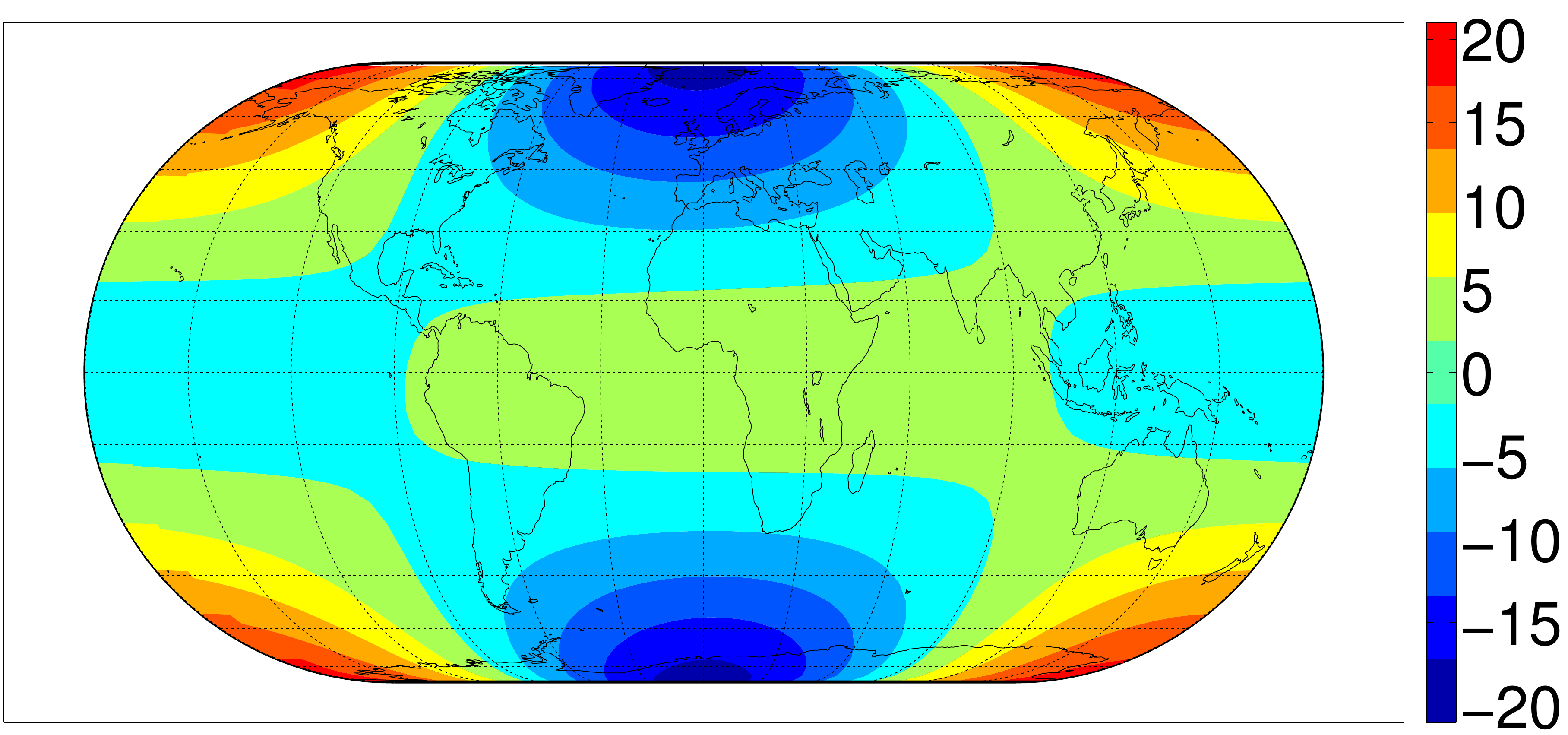}   
 \label{fig:SWE_IC_Ref_V_Same_B0_window_3}}
  \hfill
  \subfloat[HMC smoother analysis at the initial time, H component]{%
 \includegraphics[width=0.30\linewidth,height=1.90cm]{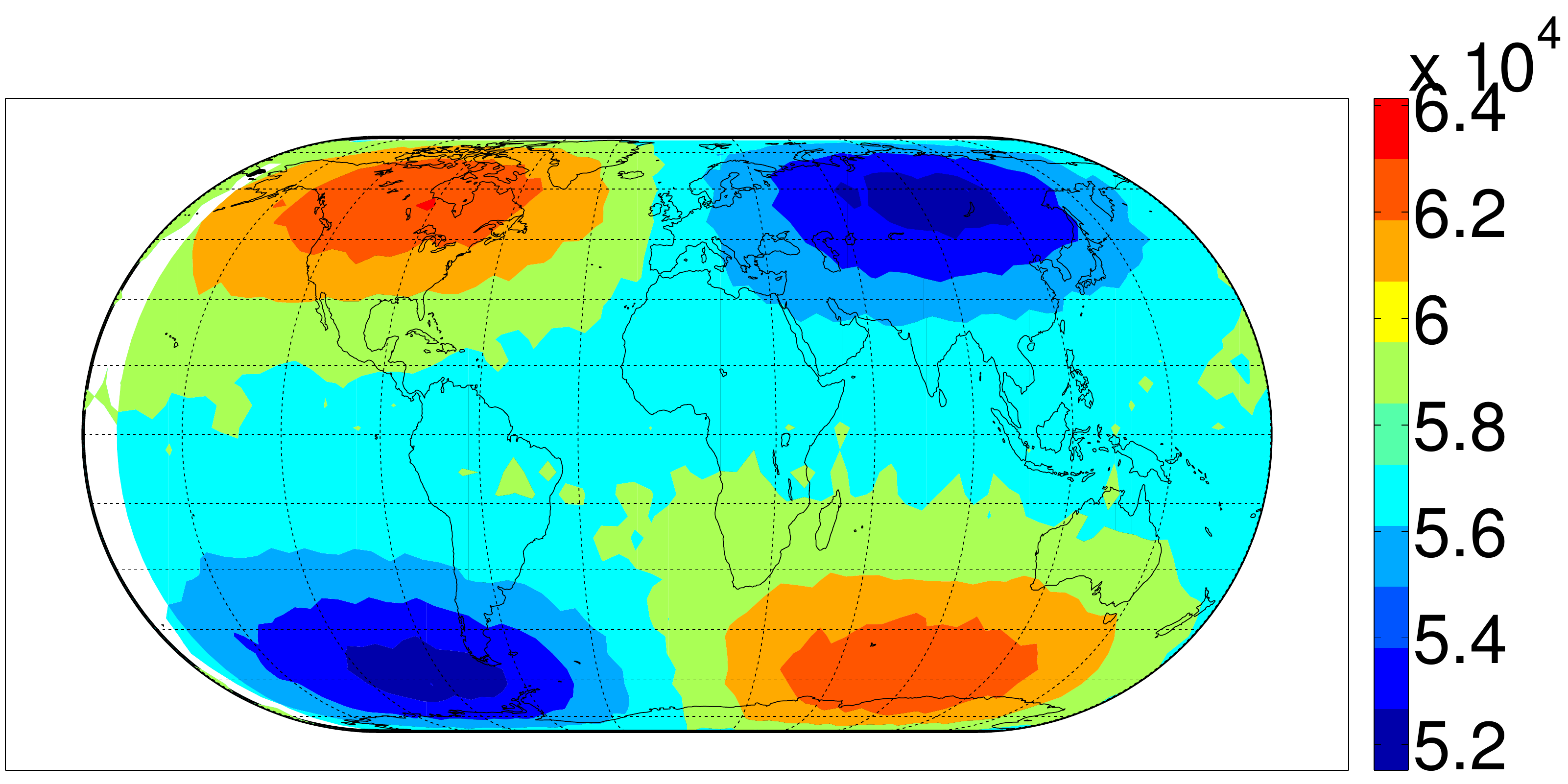}   
 \label{fig:SWE_IC_HMC_H_Same_B0_window_3}}
  \hfill 
  \subfloat[HMC smoother analysis at the initial time, U component]{%
 \includegraphics[width=0.30\linewidth,height=1.90cm]{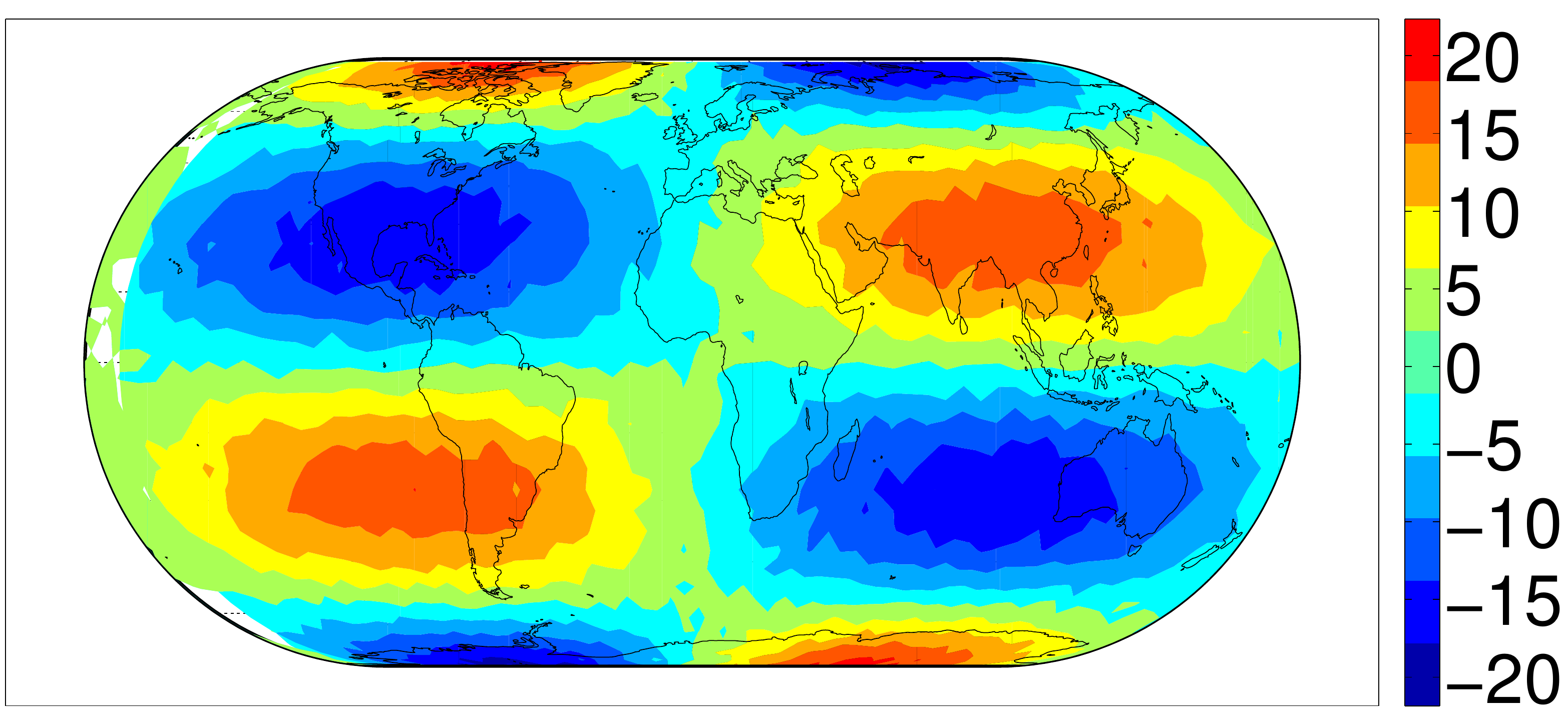}   
 \label{fig:SWE_IC_HMC_U_Same_B0_window_3}}
  \hfill 
  \subfloat[HMC smoother analysis at the initial time, V component]{%
 \includegraphics[width=0.30\linewidth,height=1.90cm]{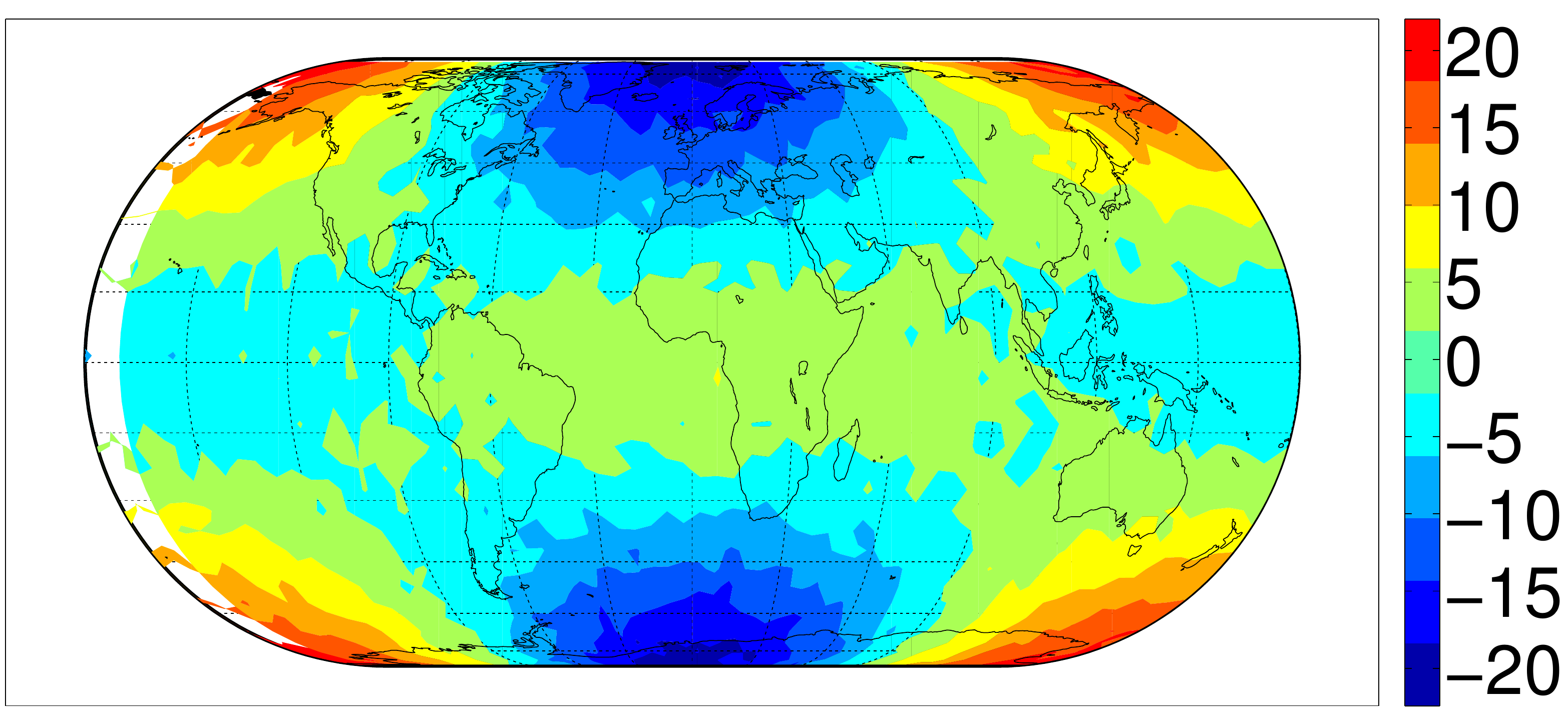}   
 \label{fig:SWE_IC_HMC_v_Same_B0_window_3}}
  \hfill
  \subfloat[4D-Var analysis at the initial time, H component]{%
 \includegraphics[width=0.30\linewidth,height=1.90cm]{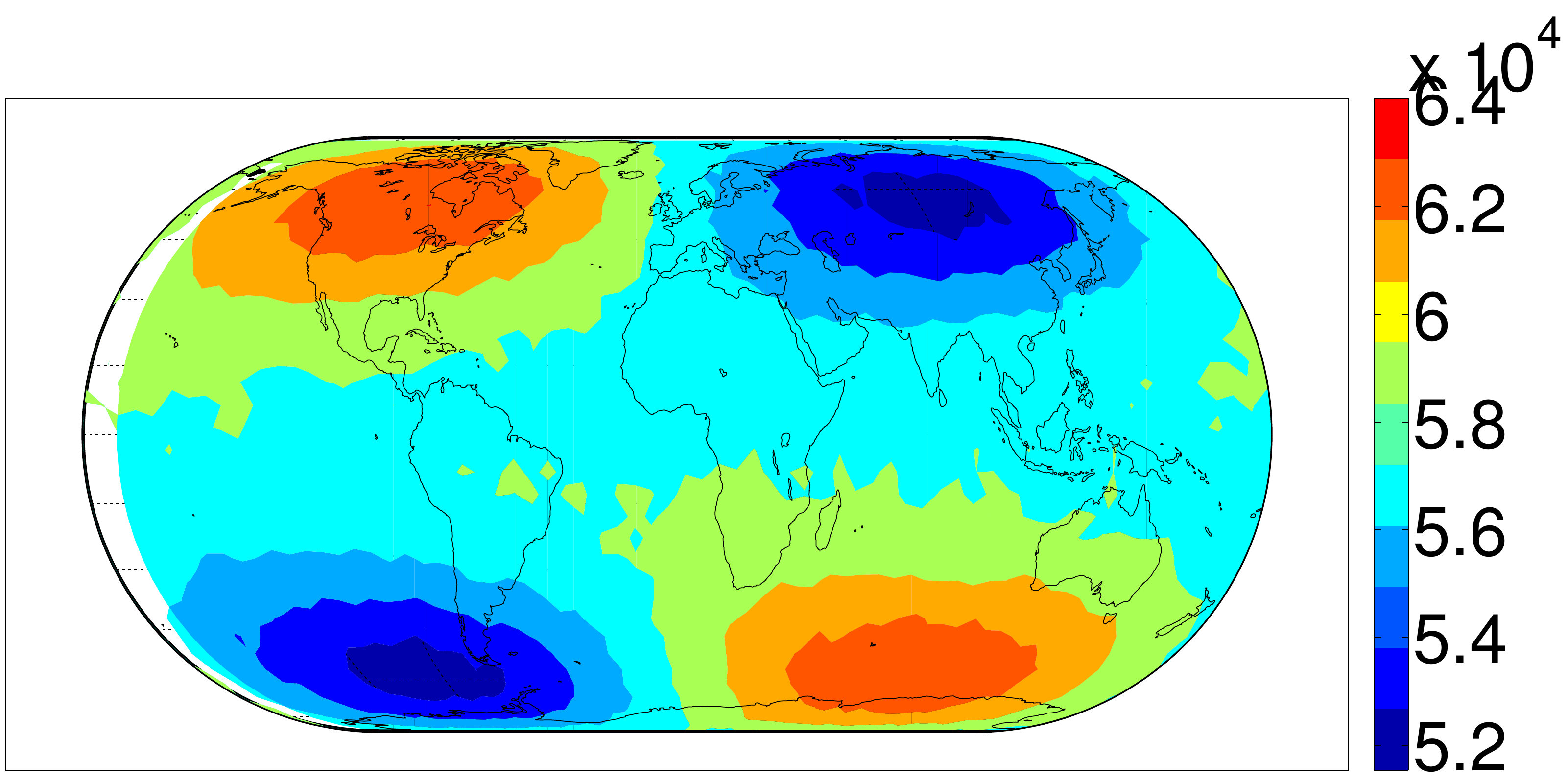}   
 \label{fig:SWE_IC_4DVAR_H_Same_B0_window_3}}
  \hfill
  \subfloat[4D-Var analysis at the initial time, U component]{%
 \includegraphics[width=0.30\linewidth,height=1.90cm]{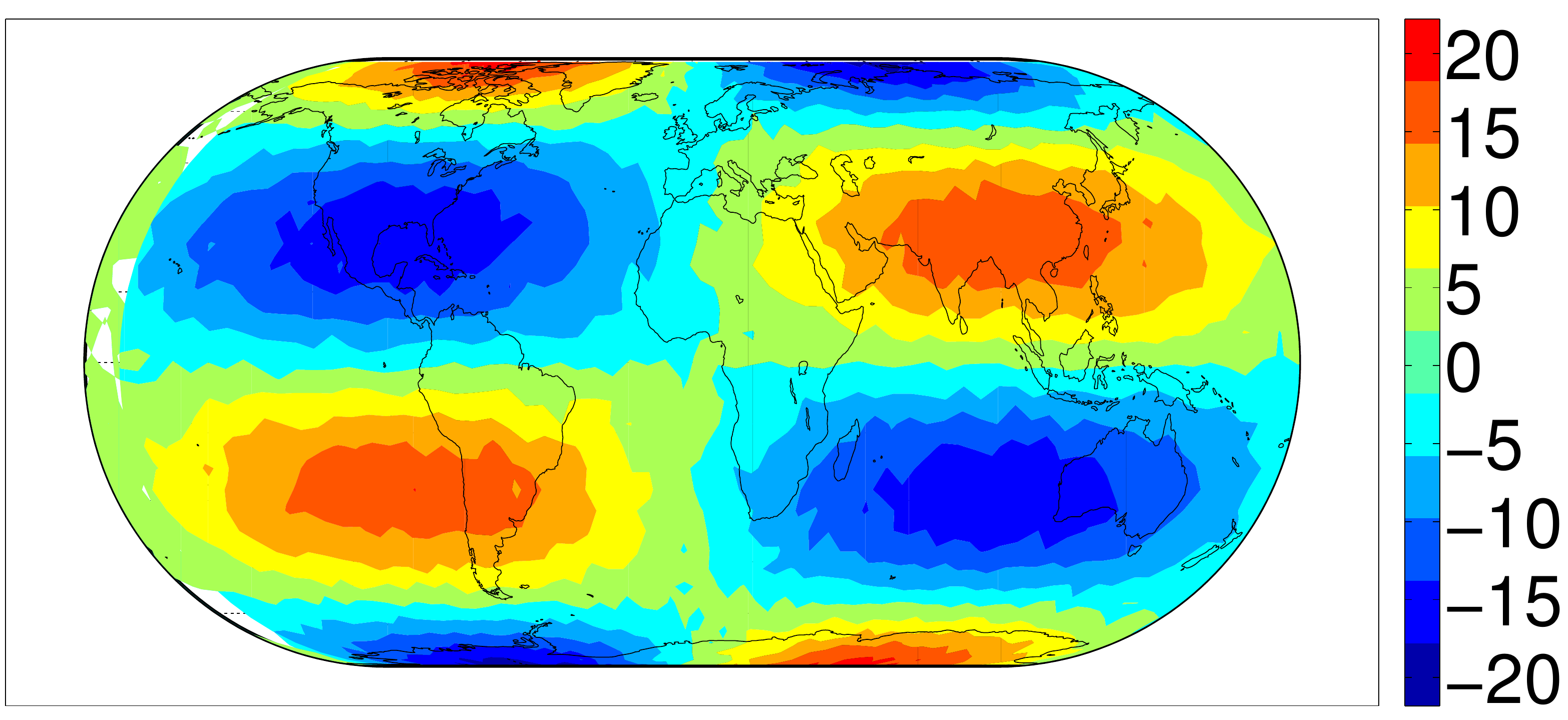}   
 \label{fig:SWE_IC_4DVAR_U_Same_B0_window_3}}
  \hfill  
  \subfloat[4D-Var analysis at the initial time, V component]{%
 \includegraphics[width=0.30\linewidth,height=1.90cm]{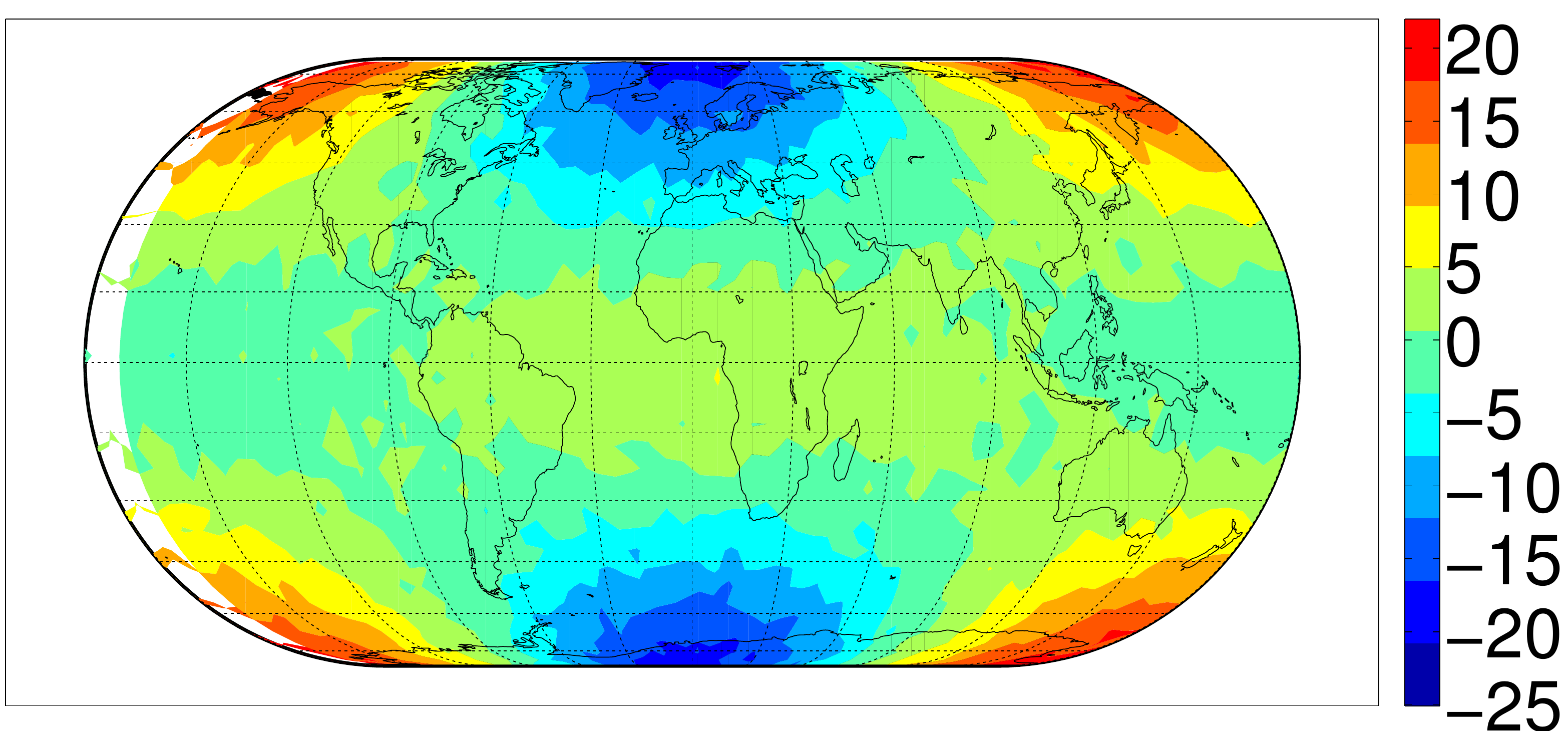}   
 \label{fig:SWE_IC_4DVAR_V_Same_B0_window_3}}
  %
    % caption and label of the whole figure
    \caption{Four dimensional data assimilation results with linear observations. The initial condition solutions at the beginning of the third window are shown. The data assimilation scheme and the state components are indicated under each panel. The assimilation window length is $8$ hours, with hourly observations.
The background error covariance matrix $\mathbf{B}_0$ is kept fixed.} 
    \label{fig:SWE_IC_Solutions_Same_B0_Window_3}
    \end{figure*}

The assimilation results obtained over the next two windows with $\mathbf{B}_0$ kept fixed are shown in Figures \ref{fig:SWE_IC_Solutions_Same_B0_Window_2}  and \ref{fig:SWE_IC_Solutions_Same_B0_Window_3}. The performance of the two schemes, 4D-Var, and the HMC smoother is quite similar, and the HMC smoother analysis competes with the 4D-Var analysis.

\begin{figure*}[htpb!]
    \centering
  \subfloat[Reference solution at the initial time, H component]{%
 \includegraphics[width=0.30\linewidth,height=1.90cm]{SWE_Reference_IC_H_Window_2_Same_B0}   
 \label{fig:SWE_IC_Ref_H_Update_B0_window_2}}
  \hfill
  \subfloat[Reference solution at the initial time, U component]{%
 \includegraphics[width=0.30\linewidth,height=1.90cm]{SWE_Reference_IC_U_Window_2_Same_B0}   
 \label{fig:SWE_IC_Ref_U_Update_B0_window_2}}
  \hfill
  \subfloat[Reference solution at the initial time, V component]{%
 \includegraphics[width=0.30\linewidth,height=1.90cm]{SWE_Reference_IC_V_Window_2_Same_B0}   
 \label{fig:SWE_IC_Ref_V_Update_B0_window_2}}
  \hfill
  \subfloat[HMC smoother analysis at the initial time, H component]{%
 \includegraphics[width=0.30\linewidth,height=1.90cm]{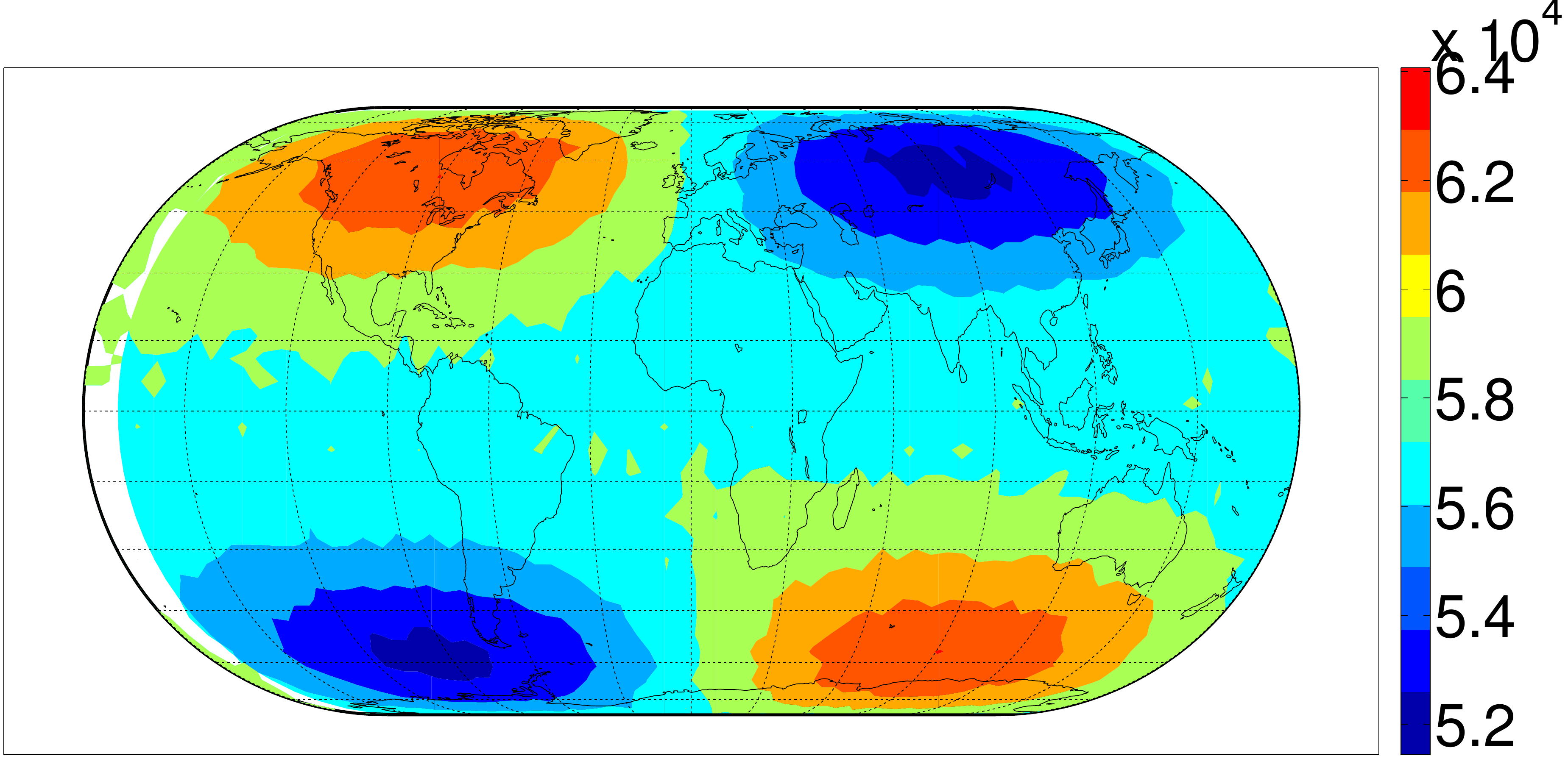}   
 \label{fig:SWE_IC_HMC_H_Update_B0_window_2}}
  \hfill  
  \subfloat[HMC smoother analysis at the initial time, U component]{%
 \includegraphics[width=0.30\linewidth,height=1.90cm]{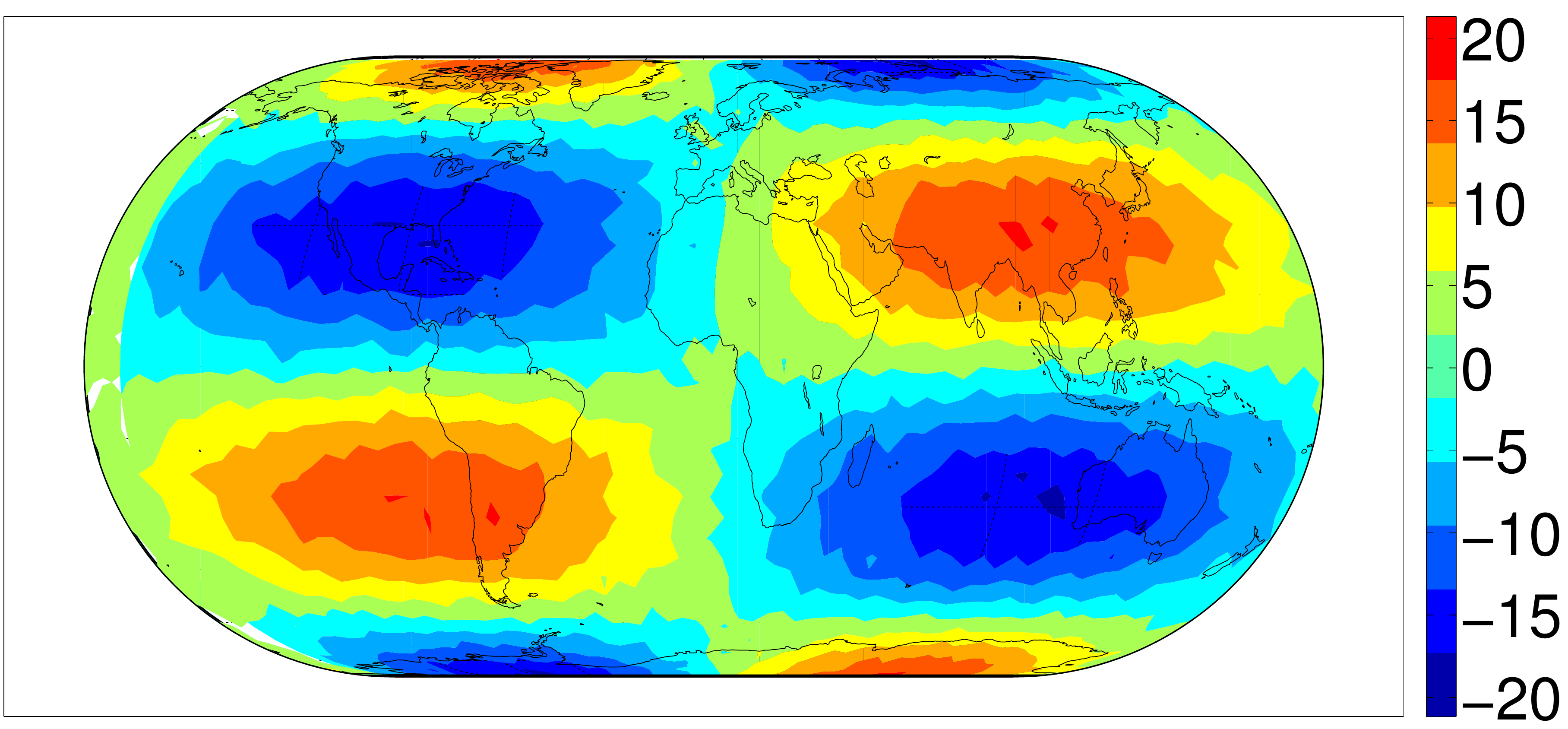}   
 \label{fig:SWE_IC_HMC_U_Update_B0_window_2}}
  \hfill  
  \subfloat[HMC smoother analysis at the initial time, V component]{%
 \includegraphics[width=0.30\linewidth,height=1.90cm]{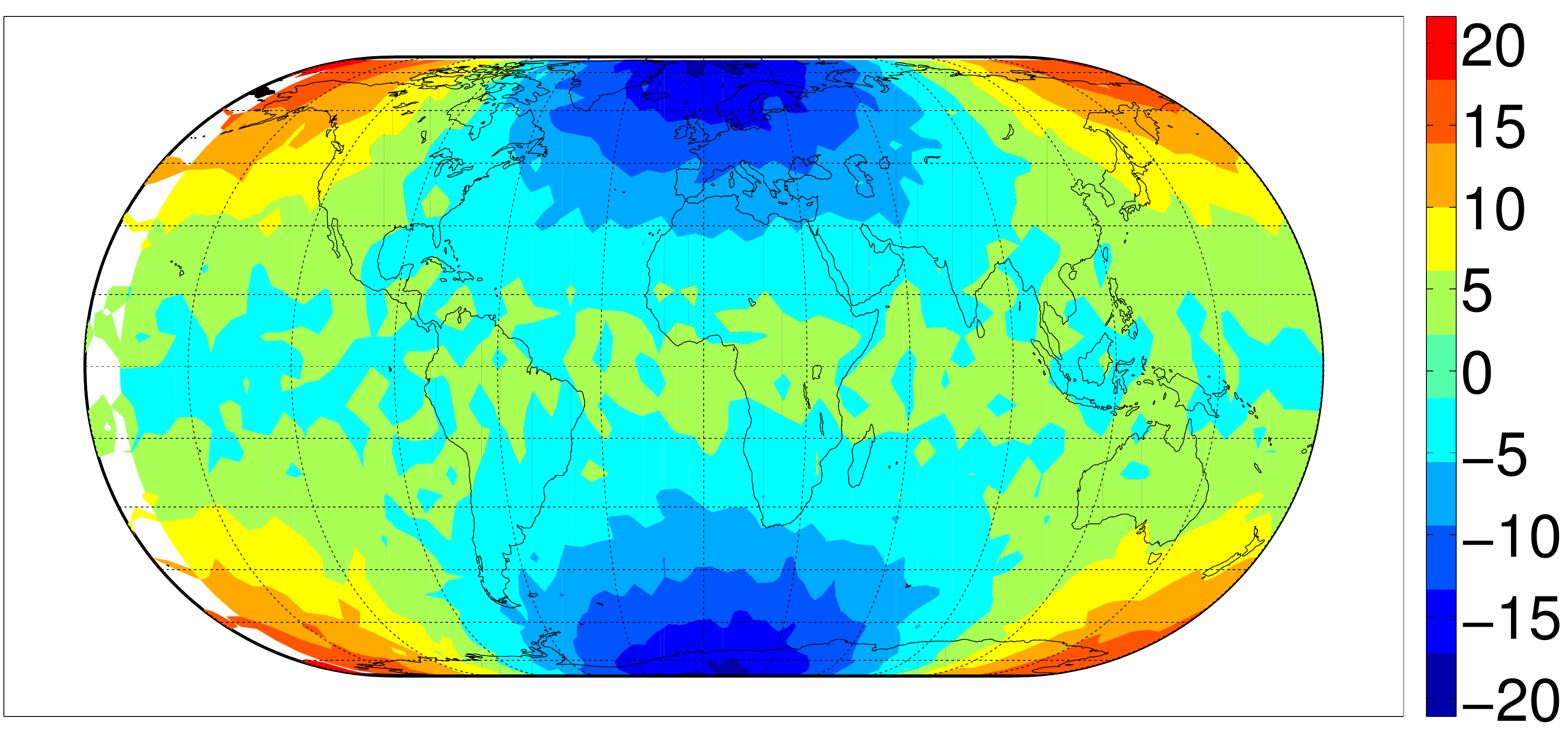}   
 \label{fig:SWE_IC_HMC_v_Update_B0_window_2}}
  \hfill
  \subfloat[4D-Var analysis at the initial time, H component]{%
 \includegraphics[width=0.30\linewidth,height=1.90cm]{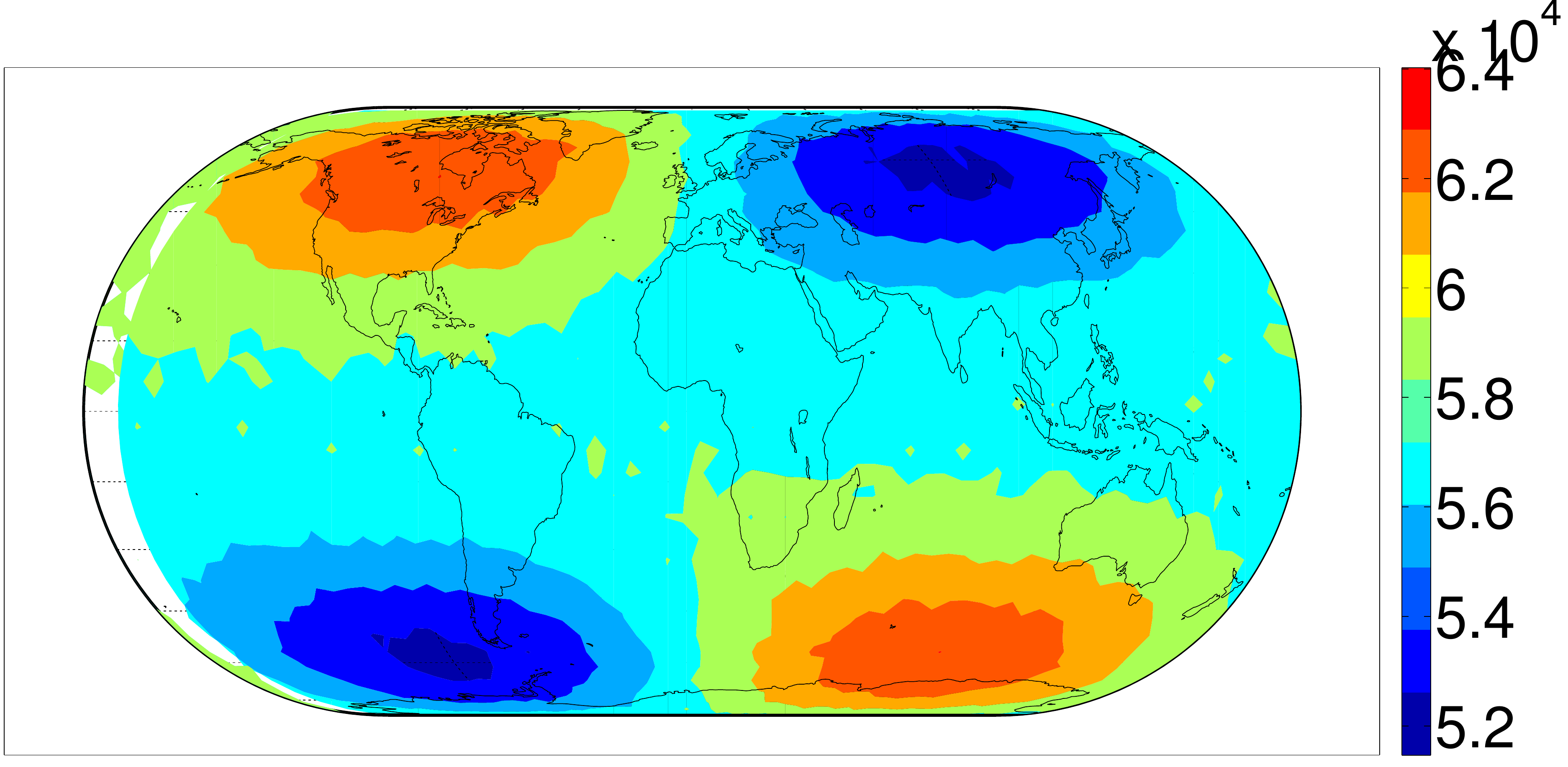}   
 \label{fig:SWE_IC_4DVAR_H_Update_B0_window_2}}
  \hfill
  \subfloat[4D-Var analysis at the initial time, U component]{%
 \includegraphics[width=0.30\linewidth,height=1.90cm]{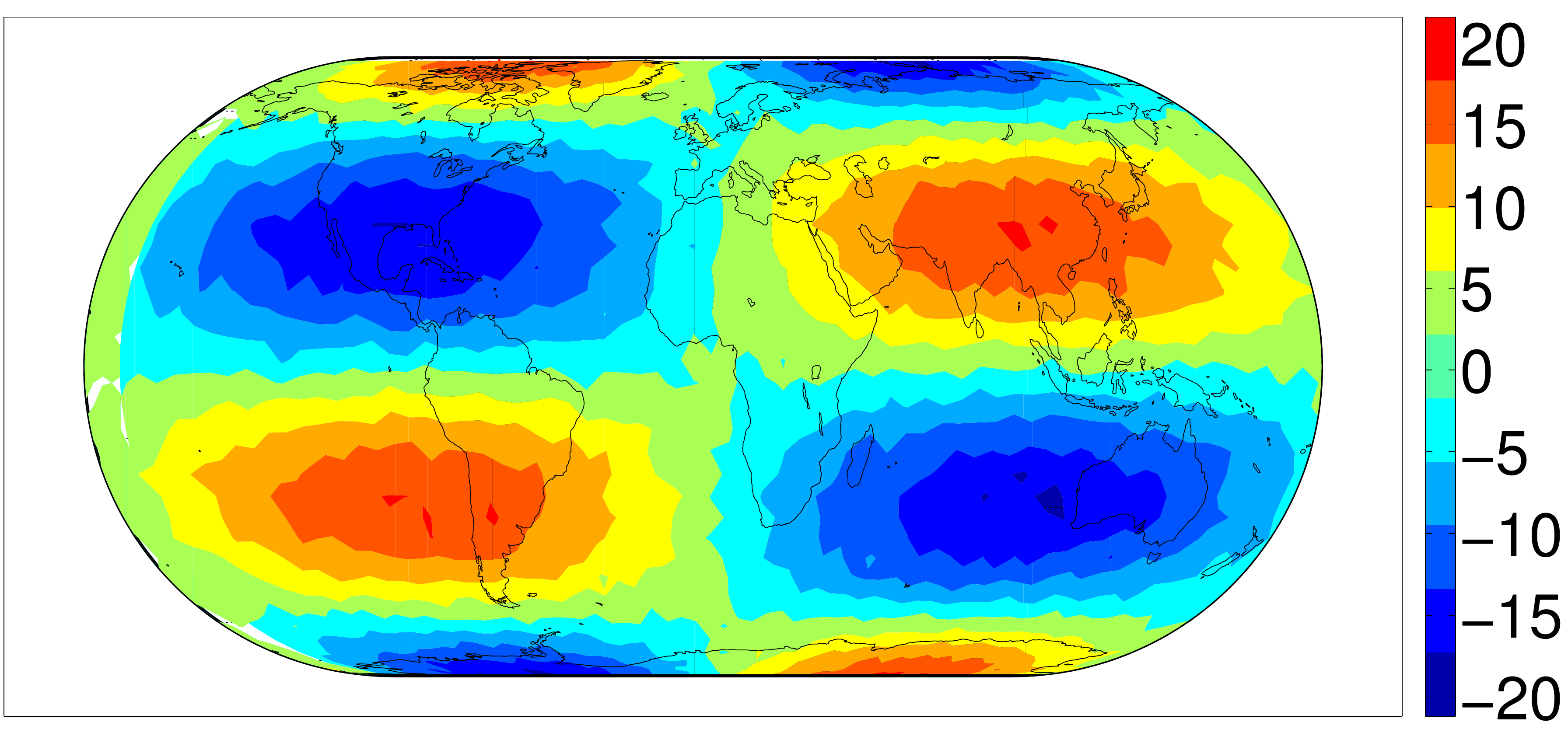}   
 \label{fig:SWE_IC_4DVAR_U_Update_B0_window_2}}
  \hfill  
  \subfloat[4D-Var analysis at the initial time, V component]{%
 \includegraphics[width=0.30\linewidth,height=1.90cm]{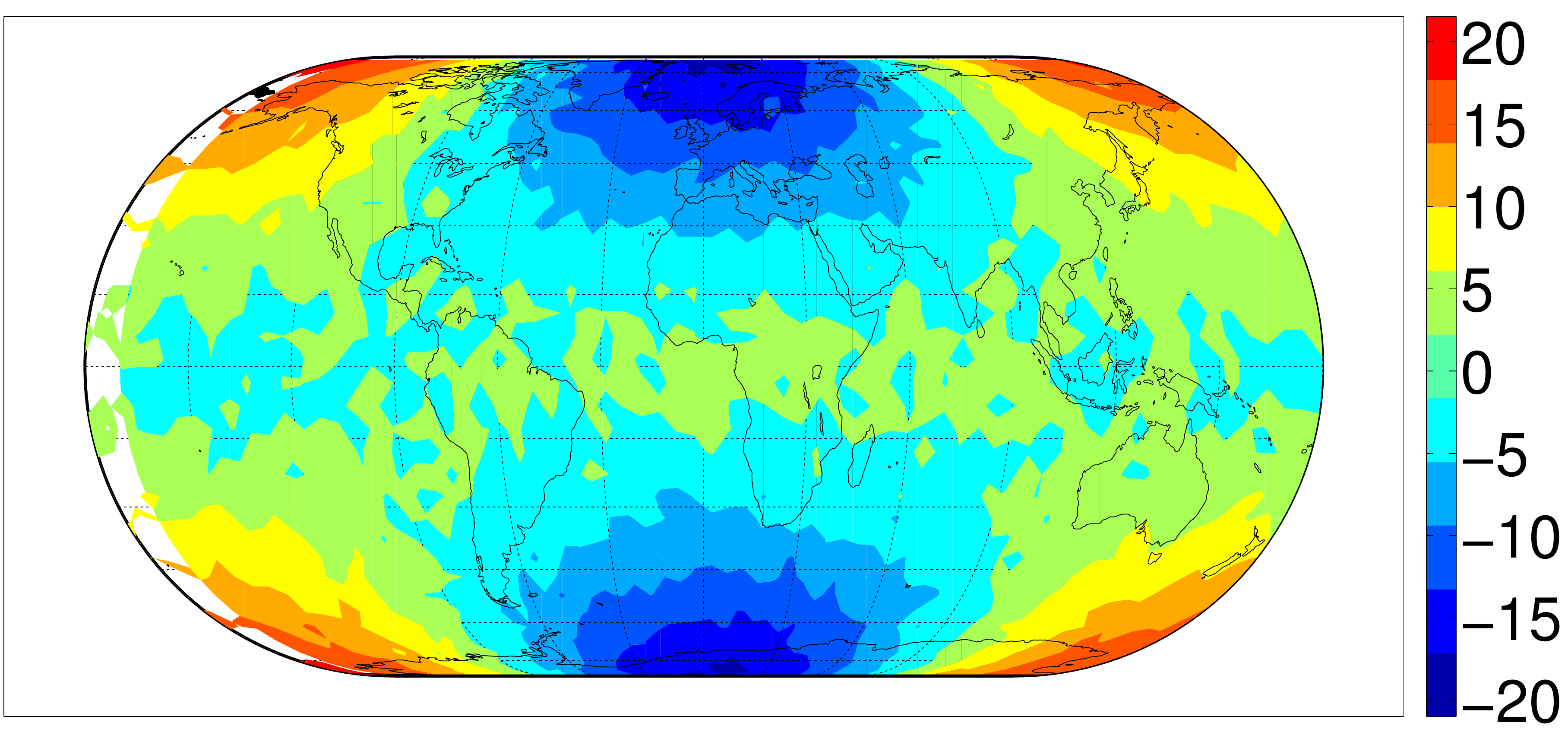}   
 \label{fig:SWE_IC_4DVAR_V_Update_B0_window_2}}
  %
    % caption and label of the whole figure
    \caption{Four dimensional data assimilation results with linear observations. The initial condition solutions at the beginning of the second window are shown. The data assimilation scheme and the state components are indicated under each panel. The assimilation window length is $8$ hours, with hourly observations. The background error covariance matrix $\mathbf{B}_0$ is updated using \eqref{eqn:Linear_Update_of_B0} for both schemes.} 
    \label{fig:SWE_IC_Solutions_Update_B0_Window_2}
    \end{figure*}

\begin{figure*}[htpb!]
    \centering
  \subfloat[Reference solution at the initial time, H component]{%
 \includegraphics[width=0.30\linewidth,height=1.90cm]{SWE_Reference_IC_H_Window_3_Same_B0}   
 \label{fig:SWE_IC_Ref_H_Update_B0_window_3}}
  \hfill
  \subfloat[Reference solution at the initial time, U component]{%
 \includegraphics[width=0.30\linewidth,height=1.90cm]{SWE_Reference_IC_U_Window_3_Same_B0}   
 \label{fig:SWE_IC_Ref_U_Update_B0_window_3}}
  \hfill
  \subfloat[Reference solution at the initial time, V component]{%
 \includegraphics[width=0.30\linewidth,height=1.90cm]{SWE_Reference_IC_V_Window_3_Same_B0}   
 \label{fig:SWE_IC_Ref_V_Update_B0_window_3}}
  \hfill
  \subfloat[HMC smoother analysis at the initial time, H component]{%
 \includegraphics[width=0.30\linewidth,height=1.90cm]{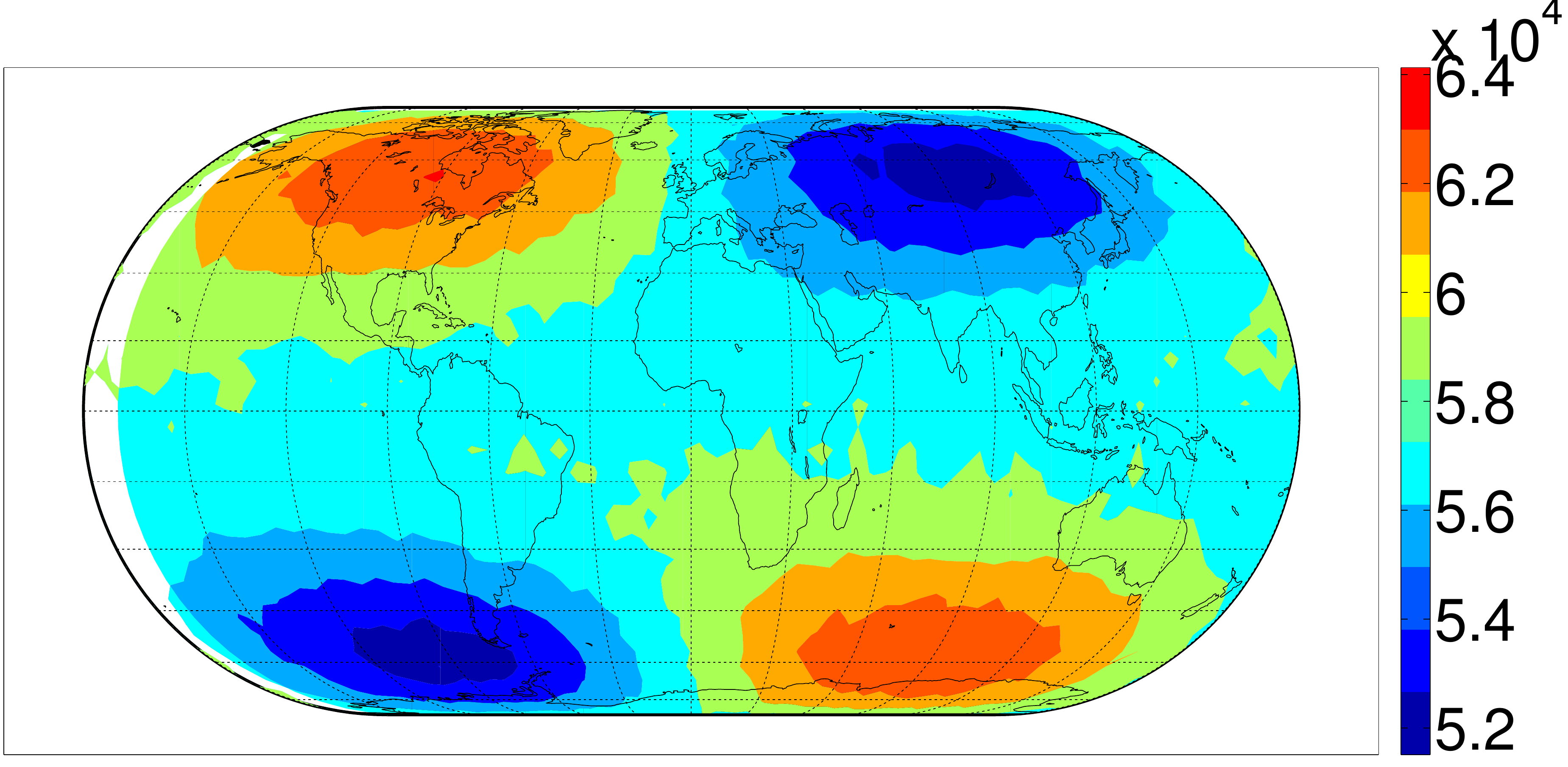}   
 \label{fig:SWE_IC_HMC_H_Update_B0_window_3}}
  \hfill  
  \subfloat[HMC smoother analysis at the initial time, U component]{%
 \includegraphics[width=0.30\linewidth,height=1.90cm]{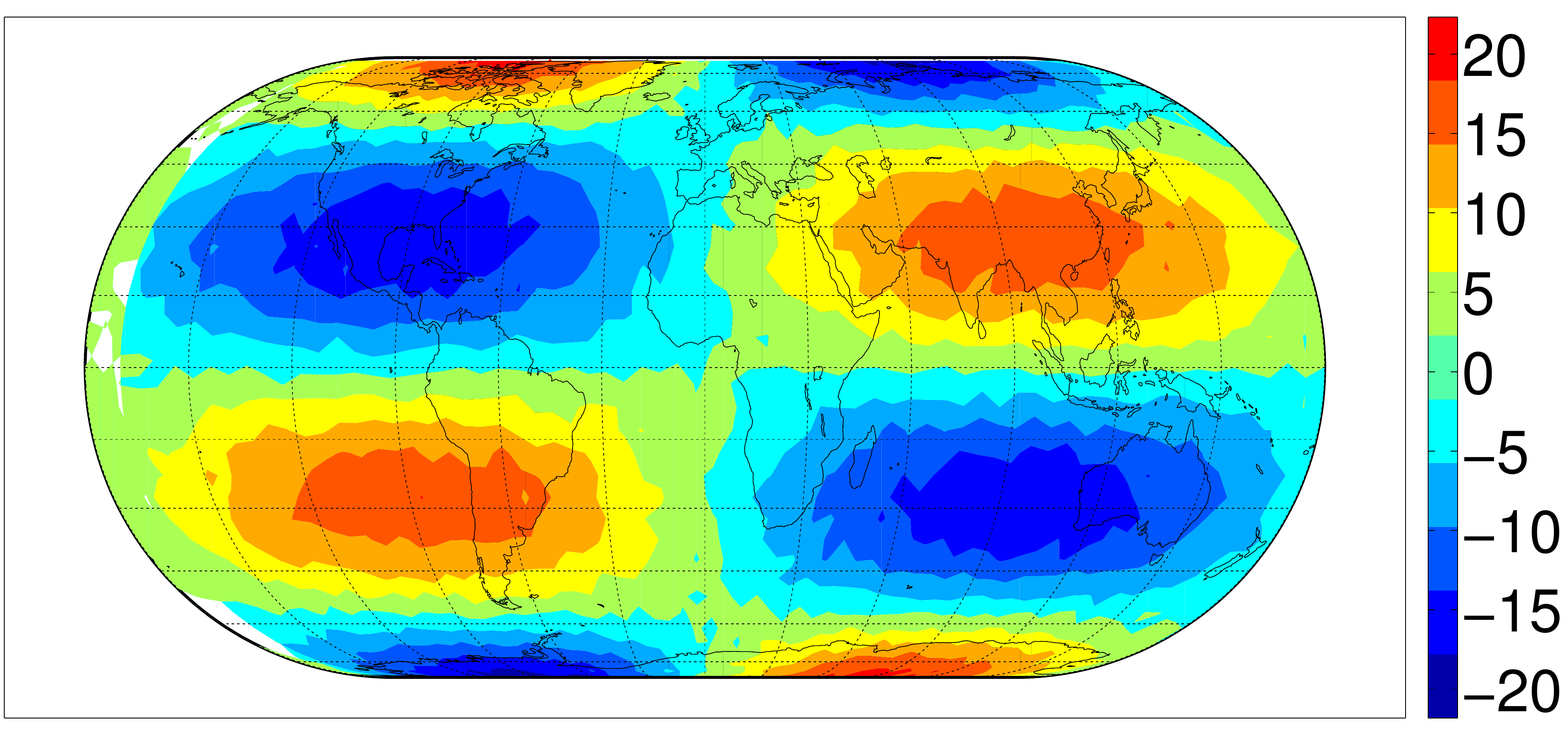}   
 \label{fig:SWE_IC_HMC_U_Update_B0_window_3}}
  \hfill  
  \subfloat[HMC smoother analysis at the initial time, V component]{%
 \includegraphics[width=0.30\linewidth,height=1.90cm]{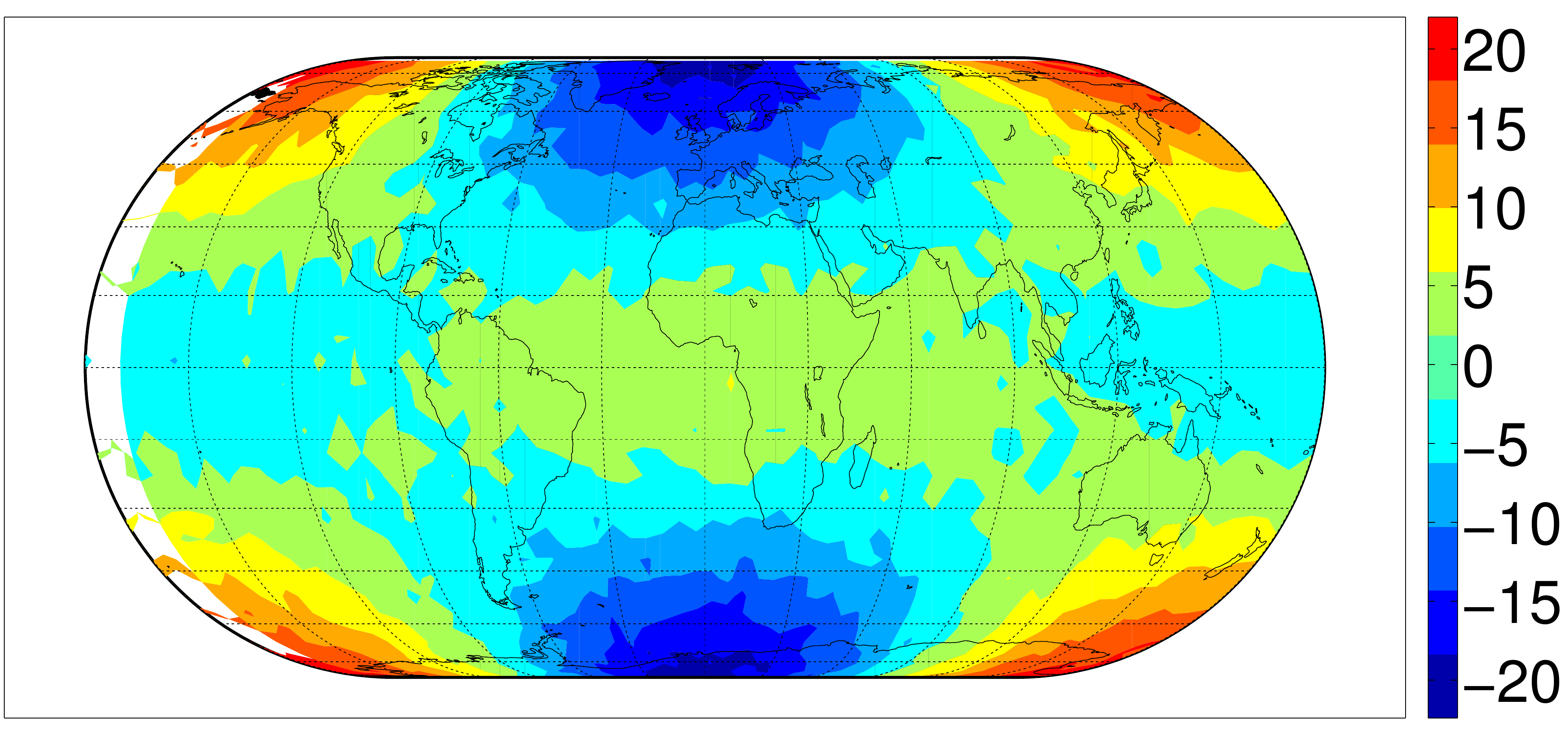}   
 \label{fig:SWE_IC_HMC_v_Update_B0_window_3}}
  \hfill
  \subfloat[4D-Var analysis at the initial time, H component]{%
 \includegraphics[width=0.30\linewidth,height=1.90cm]{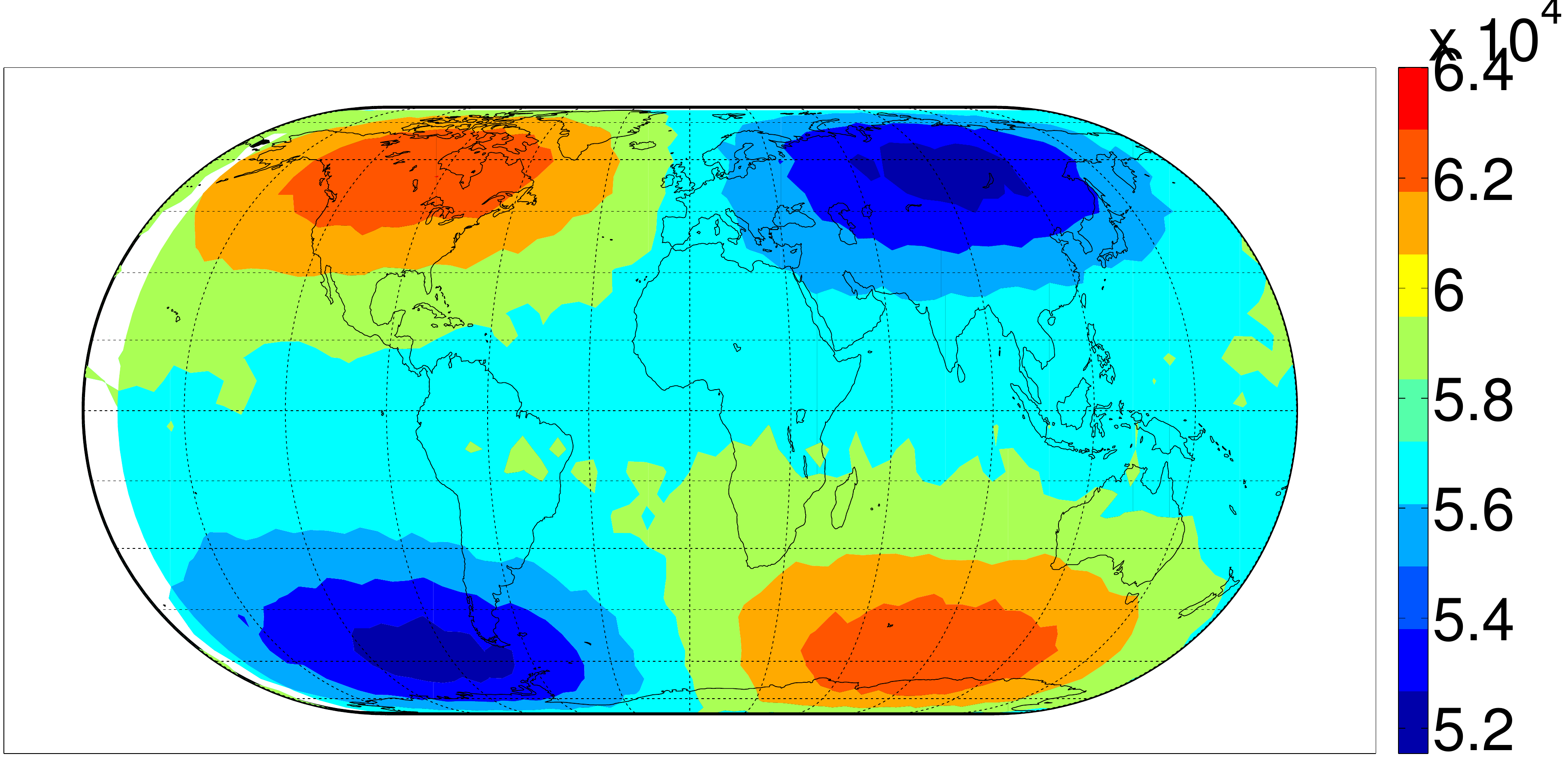}   
 \label{fig:SWE_IC_4DVAR_H_Update_B0_window_3}}
  \hfill
  \subfloat[4D-Var analysis at the initial time, U component]{%
 \includegraphics[width=0.30\linewidth,height=1.90cm]{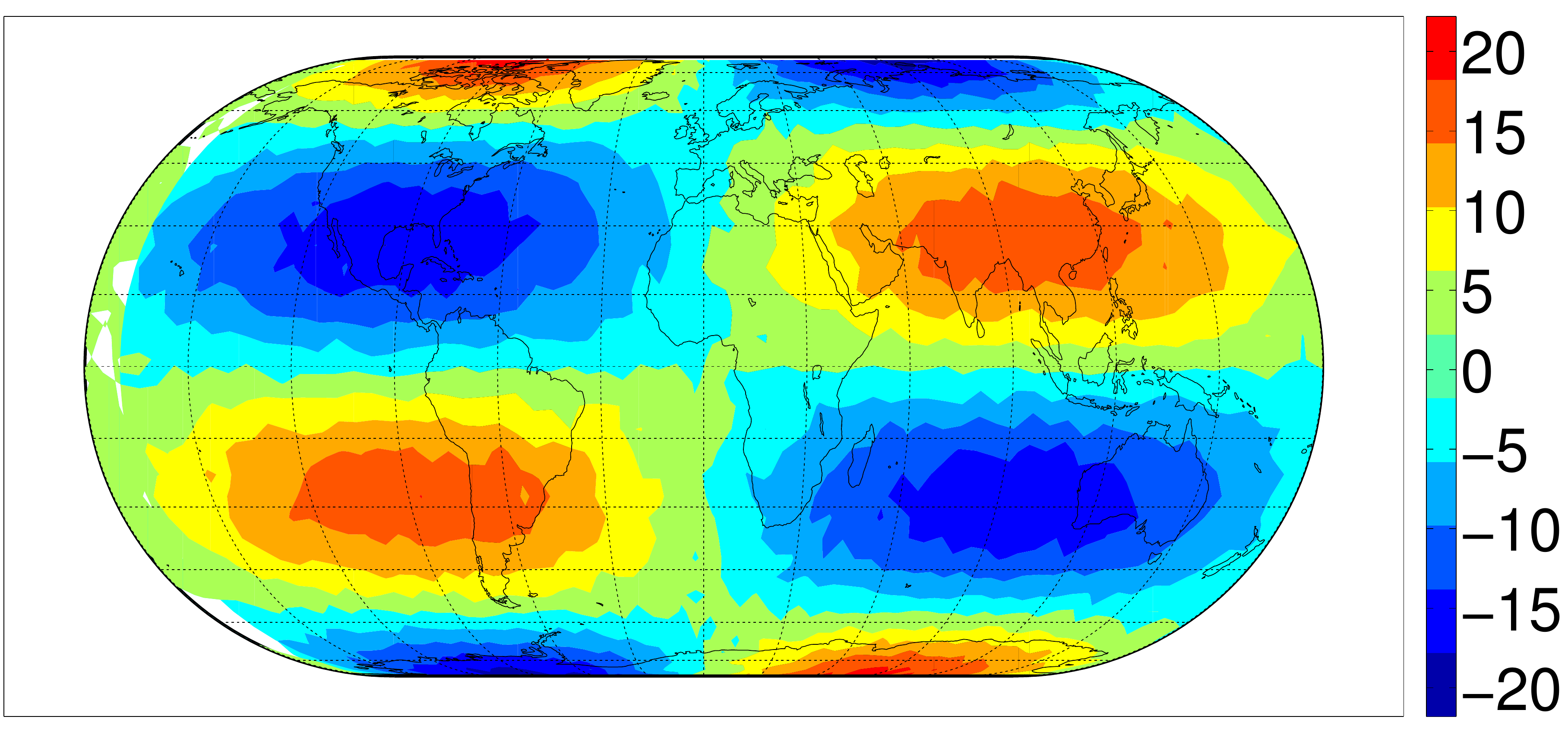}   
 \label{fig:SWE_IC_4DVAR_U_Update_B0_window_3}}
  \hfill  
  \subfloat[4D-Var analysis at the initial time, V component]{%
 \includegraphics[width=0.30\linewidth,height=1.90cm]{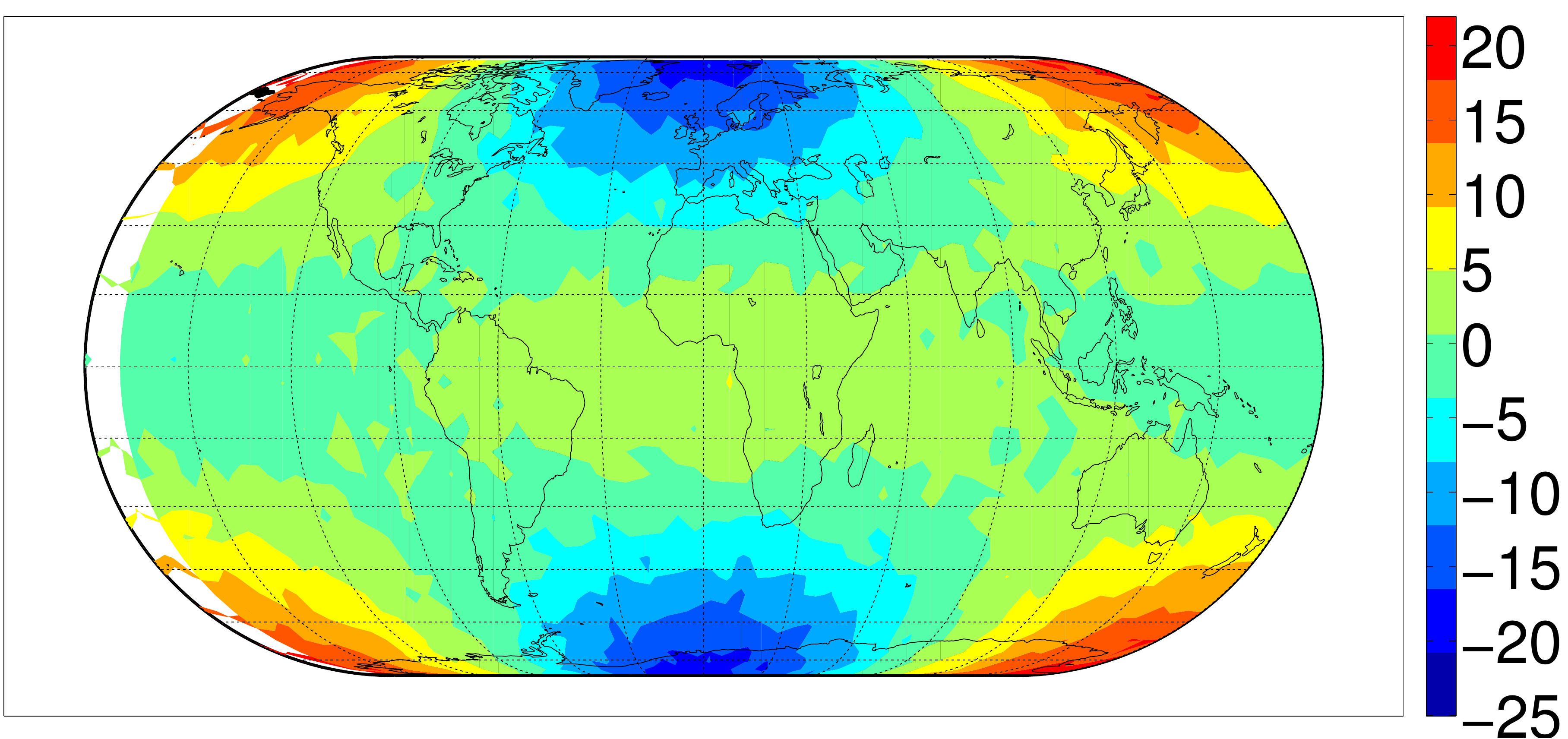}   
 \label{fig:SWE_IC_4DVAR_V_Update_B0_window_3}}
  %
    % caption and label of the whole figure
    \caption{Four dimensional data assimilation results with linear observations. The initial condition solutions at the beginning of the third window are shown. The data assimilation scheme and the state components are indicated under each panel. The assimilation window length is $8$ hours, with hourly observations. The background error covariance matrix $\mathbf{B}_0$ is updated using \eqref{eqn:Linear_Update_of_B0} for both schemes.} 
    \label{fig:SWE_IC_Solutions_Update_B0_Window_3}
    \end{figure*}

    Updating the background error covariance can, in principle, enhance the performance of both the 4D-Var and the HMC smoother. 
    In the case of 4D-Var, updating $\mathbf{B}_0$ results in lower RMSE which indicates that in real applications, the analysis is expected to  be closer to reality.
    In addition to more accurate prior kernel, updating $\mathbf{B}_0$ will result in a better update of the mass matrix $\mathbf{M}$ which in turn is expected to result in better performance of the smoother. 
    In our experiments the update has a small positive impact on the performance of the two data assimilation schemes as explained in Figures \ref{fig:SWE_RMSE}, \ref{fig:SWE_RMSE_Window_2_Zoomed}.  The positive effect here is explained by reduction in the RMS errors.
    The resulting ensemble-based forecast error covariance matrix that is used to update $\mathbf{B}_0$ can be crucial for cases where observations are sparse  or not uniformly distributed over the grid, and therefore well worth the computational overhead of the forward propagation of all the analysis ensemble members to build the full forecast ensemble. 
    Results of the data assimilation process with hybrid (updated) background error covariance matrix on the next two windows are shown in Figures \ref{fig:SWE_IC_Solutions_Update_B0_Window_2} and \ref{fig:SWE_IC_Solutions_Update_B0_Window_3}. 

    %
%    
    
%%%%%%%%%%%%%%%%%%%%%%%%%%%%%%%%%%%%%%%%%%%%%%%
\subsection{Computational considerations}
\label{subsec:computational_considerations}
%%%%%%%%%%%%%%%%%%%%%%%%%%%%%%%%%%%%%%%%%%%%%%%
%
Both 4D-Var and HMC smoother require the same computational infrastructure, namely, an adjoint model that computes the gradient of the 4D-Var cost functional \eqref{eqn:4DVAR_Cost_Functional_HMC}. This gradient calculation is the computational bottleneck for both 4D-Var and HMC smoother. It requires forward propagation of the model, and a backward propagation of the adjoint model. In our shallow water test the cost of one adjoint model run is approximately equivalent to $2.5$ times the cost of one forward model run. This makes the cost of gradient evaluation approximately equal to $3.5$ the cost of a forward model run.

\begin{table*}[htpb!]
\scriptsize
\centering
\caption{Data assimilation scheme cost for the shallow water model with linear observations. Number of function evaluations (forward model runs), 
         and the number of optimization iterstions (for 4D-Var) are listed.
         The cost of one adjoint run is approximately equal to $2.5$ forward model runs in the current settings.
         Total cost is approximated in terms of number of forward model runs. 
         The cost of the HMC sampling smoother is the same for the three assimilation windows.
         }
\scalebox{1.0}{ % use this to scale the table if needed
\begin{tabular}{|p{1.5cm}|p{2.9cm}|p{0.9cm}|p{0.9cm}|p{0.9cm}|p{0.9cm}|p{0.9cm}|p{0.9cm}|}   
  \hline
  \multirow{3}{1.5cm}{\textbf{Data assimilation scheme}}   &  \multirow{3}{*}{\textbf{Cost}}   & \multicolumn{6}{|c|}{\textbf{Assimilation window}}\\ \cline{3-8}
 &&  \multicolumn{2}{|c|}{$(1)$} & \multicolumn{2}{|c|}{$(2)$} & \multicolumn{2}{|c|}{$(3)$} \\  \cline{3-8}  %
&& \multicolumn{2}{|p{0.9cm}|}{Fixed $\Bini$}  & Fixed $\Bini$ & Hybrid $\Bini$ & Fixed $\Bini$ & Hybrid $\Bini$ \\ \hline \hline
   \multirow{3}{1.5cm}{\textbf{4D-Var}}  &  Function evaluations   & \multicolumn{2}{|c|}{151}   &  97   &  101  &  96   &  93\\ \cline{2-8}
&  Number of iterations   & \multicolumn{2}{|c|}{49}    &  47   &  46   &  46   &  45\\ \cline{2-8}
&  Cost in equivalent forward model runs & \multicolumn{2}{|c|}{322.5}  & 261.5 &  262  &  257  &  250.5  \\ \hline  \hline
   \multirow{3}{1.5cm}{\textbf{HMC smoother}}  & Number of proposed states  & \multicolumn{6}{|c|}{530}   \\ \cline{2-8}
&   Cost per proposal  & \multicolumn{6}{|c|}{4.5 }  \\ \cline{2-8}
&  Cost in equivalent forward model runs & \multicolumn{6}{|c|}{2,385}  \\ 
\hline 
\end{tabular} 
}
\label{table:SWE_Computational_Cost}
\end{table*}

The cost of the 4D-Var depends on the number of iterations and function evaluations required by the optimization algorithm. 
On the first window the number of iterations required by the LBFGS optimizer is $49$, with $151$ function evaluations. The total cost of the 4D-Var solution is then $151+49 \times 3.5=322.5$ equivalent forward model runs.

The cost of the HMC smoother depends on the configuration of the chain: the number of burn-in step, the number of dropped states at stationarity,  and step-size settings of the symplectic integrator. The symplectic integrator itself controls the number of adjoint runs to evaluate the gradient  of the cost functional in order to propose a new state to the chain. The size of the desired ensemble controls the length of the Markov chain and consequently the total cost of the analysis step by the HMC smoother. 
The Verlet integrator \eqref{eqn:Verlet}, used in the current experiments, requires a single adjoint run to propose a new state to the chain. 
The acceptance/rejection criterion requires an additional forward run   
to evaluate the loss of energy. This makes the cost of generating a proposal state to the  Markov chain approximately equal to $4.5$ the cost of a model run. 
On all assimilation windows the chosen ensemble size is $100$. The number of burn-in states is $30$, and $4$ states are rejected between consecutive selected samples. The HMC sampling smoother, in this case generates $30+100\times 5=530$ states to collect the analysis ensemble, with a total cost roughly equal to $530\times 4.5=2,385$ forward model runs.

The computational cost of the two data assimilation schemes, 4D-Var and HMC smoother, on each assimilation window are summarized in Table \ref{table:SWE_Computational_Cost}. Notice that the total cost of DA schemes is given in terms of the total number of forward model runs. 
On the first window the cost of the HMC smoother is approximately $9$ times the cost of the 4D-Var scheme. On the next two windows,  the HMC smoother costs roughly $11$ times the cost of the 4D-Var. A cost-reduction in 4D-Var is expected because the forecast state  is closer to the MAP than the case on the first window.
The higher computational cost of the HMC smoother can be handled more efficiently by parallelizing the sampling scheme and the gradient calculations. Another way to reduce the computational cost of the HMC smoother is to replace the burn-in steps with a suboptimal 4D-Var obtained using a small number of iterations.
The computational cost of the proposed sampling smoother can of course by be reduced by decreasing the ensemble size, however this will result in higher sampling error. The impact of the sampling errors can be assessed by the techniques developed in \cite{Rao:2014,rao2014posteriori}. 

The increased cost of the HMC smoother could be acceptable in view of the additional useful information it provides: a sample estimate of the analysis probability distribution (and as immediate consequences an analysis error covariance matrix and a flow-dependent background error covariance matrix for the next cycle). 

%
%%%%%%%%%%%%%%%%%%%%%%%%%%%%%%%%%%%%%%%%%%%%%%%
\section{Conclusion and Future Work}
\label{sec:Conclusion_and_Future_Work}  
%%%%%%%%%%%%%%%%%%%%%%%%%%%%%%%%%%%%%%%%%%%%%%%
 %
A four-dimensional data assimilation smoother is proposed in this paper. The smoother samples from the posterior distribution 
using a Hybrid Monte-Carlo approach. The 4D-Var approach provides a MAP estimate of the true state, but it does not compute a measure of uncertainty of the analysis. The HMC smoother builds an ensemble approximating the posterior PDF. This can be used to estimate the true state together with the uncertainty in  analysis, e.g., by calculating the ensemble mean and ensemble-based analysis error covariance matrix.  Moreover, propagating the analysis ensemble to the beginning of the next assimilation window provides a forecast ensemble that can be used to construct a flow-dependent background covariance matrix for this new window.  Unlike several popular hybrid approaches, the HMC smoother generates an analysis error covariance that is consistent with the analysis state -- because both statistics are produced by one consistent data assimilation scheme.

The HMC smoother  requires an adjoint of the simulation model, and runs on the same computational infrastructure as 4D-Var. The computational cost of the HMC smoother is - as of now - larger than that of 4D-Var. The efficiency issue must be addressed before the HMC smoother becomes fully practical. We are currently investigating several strategies to enhance the performance of the sampling smoother and to reduce its computational cost. Parallelizing the sampling smoother will be considered.  We will also test the HMC smoother on the case of nonlinear observations, imperfect models, and non-Gaussian errors.

%    

%XXXXXXXXXXXXXXXXXXXXXXXXXXXXXXXXXXXXXXXXXXXXXXXXXXXXXXXXXXbibitem version of the referencesXXXXXXXXXXXXXXXXXXXXXXXXXXXXXXXXXXXXXXXXXXXXXXXXXXXXXXXXXXXXXXX
%\FloatBarrier
%\clearpage
\section*{Acknowledgements}
This work was partially supported by awards AFOSR FA9550--12--1--0293--DEF, NSF CCF--1218454, NSF DMS--1419003, and by the Computational Science Laboratory at Virginia Tech.

% \bibliographystyle{spbasic}% basic style, author-year citations
%\bibliographystyle{spmpsci}% mathematics and physical sciences
%\bibliographystyle{spphys} % APS-like style for physics
%
% \bibliographystyle{plain}
%\bibliography{data_assim_general,data_assim_kalman,data_assim_fdvar,data_assim_HMC}
%
%XXXXXXXXXXXXXXXXXXXXXXXXXXXXXXXXXXXXXXXXXXXXXXXXXXXXXXXXXXbibitem version of the referencesXXXXXXXXXXXXXXXXXXXXXXXXXXXXXXXXXXXXXXXXXXXXXXXXXXXXXXXXXXXXXXX
%    
%
%
%
%
\end{document}